\documentclass[12pt,a4paper]{article}
\usepackage{epsfig}
\usepackage{epsfig,amssymb,amsmath,amsthm, color, graphics,hyperref,enumitem}
\usepackage[english]{babel}
\usepackage{multicol, color}
\usepackage{graphics}
\usepackage{hyperref}

\newtheorem{theo}{Theorem}
\newtheorem*{theo2}{Theorem 2'}

\newtheorem{lemm}{Lemma}
\newtheorem{coro}{Corollary}

\theoremstyle{definition}
\newtheorem{rema}{Remark}
\newtheorem{defi}{Definition}
\newtheorem{example}{Example}
\newenvironment{demo}{{\bf Proof: }}{\hfill $\diamond$\medskip}

\def\al{\alpha}
\def\om{\omega}
\def\Om{\Omega}
\def\ga{\gamma}

\def\vp{\varphi}
\def\la{\lambda}

\def\si{\sigma}
\def\Si{\Sigma}
\def\ep{\varepsilon} 

\def\R{\mathbb R}
\def\Z{\mathbb Z}

\def\decomposedscheme{decomposed scheme }
\begin{document}
\sloppy
\date{23 Nov 2017}
\title{A complete topological classification of Morse-Smale diffeomorphisms on surfaces: a kind of kneading theory in dimension two}
\author{V.\,Z.~Grines\footnote{National Research University Higher School of Economics, vgrines@yandex.ru}, O.\,V.~Pochinka\footnote{National Research University Higher School of Economics, olga-pochinka@yandex.ru}, S.\,van~Strien\footnote{Imperial College London, 	
s.van-strien@imperial.ac.uk}}

\maketitle

\begin{abstract} 
{In this paper we give a complete topological classification of orientation preserving  Morse-Smale  diffeomorphisms on orientable closed surfaces. For MS diffeomorphisms with relatively simple behaviour it was known that such a classification 
can be given through a directed graph, a three-colour directed graph or by a certain topological object, called a {\lq}scheme{\rq}. 
Here we will assign to {any} MS  surface diffeomorphism a finite amount of data
which completely determines its topological conjugacy class. 
Moreover,  we show that associated to any abstract version of this data, there exists a unique conjugacy class of MS orientation preserving diffeomorphisms (on some orientation preserving surface).
As a corollary we obtain a different proof that nearby MS diffeomorphisms are topologically conjugate.}
\end{abstract}

{\bf Key words}: Morse-Smale diffeomorphism, topological classification

Bibliography: \ref{last} names.

\tableofcontents

\section{Informal statement of the results}

One of the main objectives in the field of dynamical systems is to obtain a classification in
terms of their dynamics. Such a classification was achieved successfully in the one-dimensional
setting. For example, two circle diffeomorphisms $f,f' : S^1 \to S^1$ are called topologically conjugate if
there exists an orientation preserving homeomorphism $h: S^1 \to S^1$ so that $hf = fh$. In the
1880's, Poincare \cite{Po} showed that if these two diffeomorphisms are topologically conjugate then they have the same rotation number. Moreover, for  two transitive homeomorphisms $f,f'$ the condition that their rotation numbers are the same is both necessary and sufficient for the topological conjugacy. In 1932, Denjoy \cite{De}  improved this result by showing that if the diffeomorphism f is $C^2$ and have no periodic orbits then it is transitive.

If two circle diffeomorphisms have periodic orbits and each of these periodic orbits is hyperbolic (such a diffeomorphism is called {\em Morse-Smale}),  then a necessary and sufficient
condition for these circle diffeomorphisms to be topologically conjugate is that their rotation numbers are the same and that they have the same number of periodic attractors. 
Moreover, if one chooses a rotation compatible permutation on a finite number of points on the circle, then this data corresponds to a Morse-Smale diffeomorphism. 
For {\em non-invertible} Morse-Smale maps of the circle or the interval one has a similar situation: it's so-called  {\em kneading map} (describing itineraries of its turning points) 
is  (essentially) a complete topological invariant and, moreover, each admissible kneading map corresponds to a map of the circle. 

\begin{figure}[h!]
\centering
\includegraphics[width=0.9\textwidth]{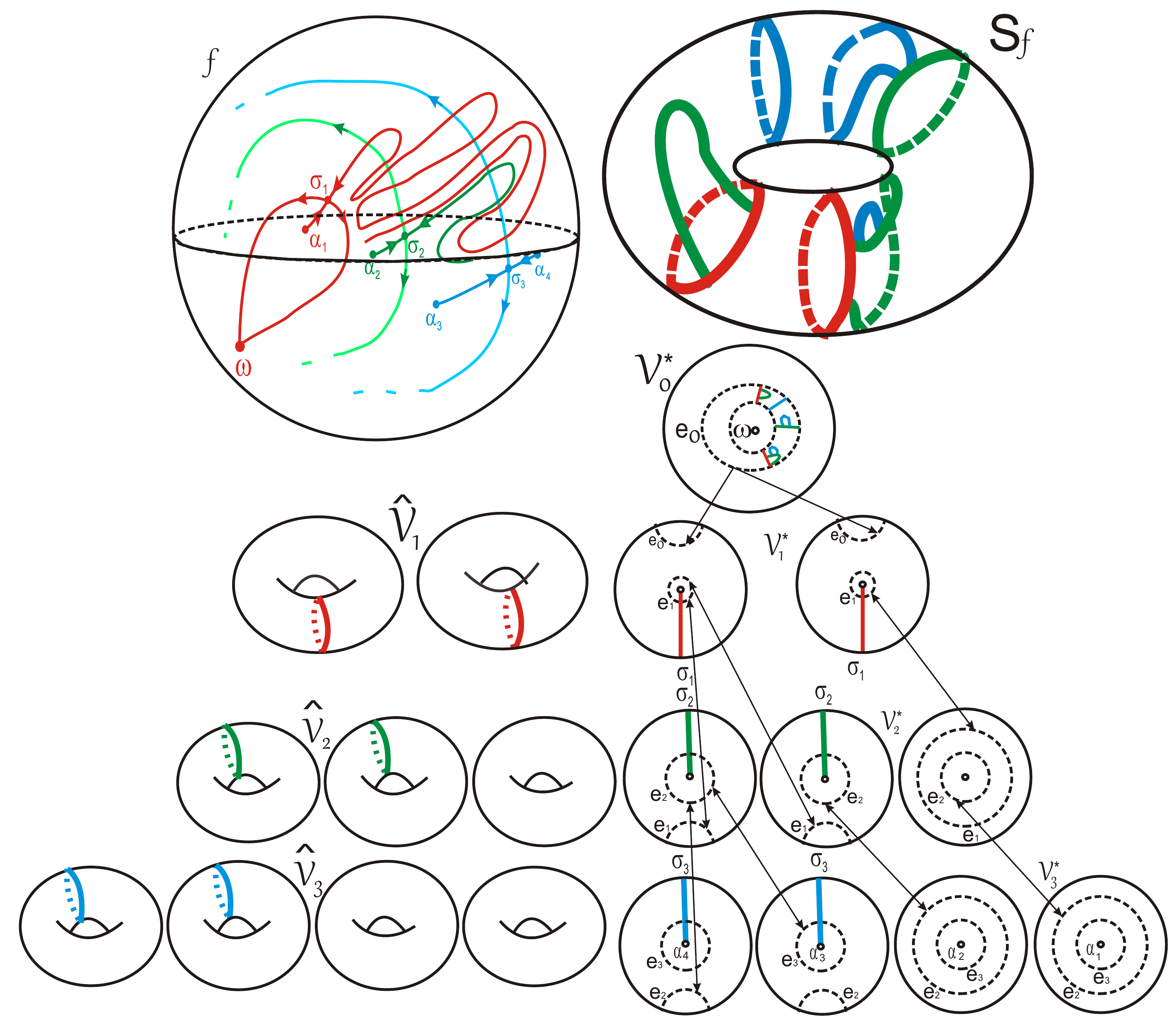}\label{faz+}
\caption{\small Theorem A gives a complete topological invariant for a Morse-Smale diffeomorphism $f$ on the 2-sphere 
(drawn on the top left) by its scheme $S_f$
(on the top right)
and its decomposed scheme (depicted in the remaining figures), while  Theorem B describes all possible Morse-Smale diffeomorphisms. This picture  and the notions of a scheme and a decomposed scheme will be explained in Subsection~\ref{subsec:statementresults}.}
\end{figure}

The {\em aim} of this paper is to establish a corresponding classification in the setting of Morse-Smale diffeomorphisms on {\em closed orientable surfaces}, 
replacing a finite number of points on  a circle by a finite number of annuli on tori.  
In Section~\ref{subsec:statementresults} we will state our results precisely, and define the notion of a {\lq}scheme{\rq} 
and {\lq}decomposed scheme{\rq}, 
but {\em informally speaking} (see Figure~\ref{faz+}) we have

\medskip 

\noindent 
{\bf Theorem A (Classification by finite amount of data)} {\em Let $M$ be a closed orientable surface and $f\colon M\to M$ be an orientation preserving Morse-Smale diffeomorphism. Then one can assign to $f$ a {\bf scheme}  $S_f$ {or a {\bf decomposed scheme}}
 consisting of {\bf  a finite amount data} (given by  a finite union of tori, and the homotopy type of certain annuli in these tori), in such a way that  $f\colon M\to M$ and $f'\colon M'\to M'$ are topologically conjugate if and only if $S_f$ is equivalent to $S_{f'}$. }

\medskip 

\noindent 
{\bf Theorem B (Realisation)} {\em Each  {\bf abstract scheme} $S$ corresponds uniquely to an orientable closed surface $M$ and an  orientation preserving Morse-Smale diffeomorphism $f\colon M\to M$.}
\medskip 

For a formal statement of these results,  see Theorems 2 and 2'  in Subsections~\ref{subsec:statementThm2},\ref{deca} and 
Theorem~\ref{t.realisation} in Subsection \ref{subsection:statementTheorem3}.

\medskip 
Notice that Theorem A implies immediately that Morse-Smale diffeomorphisms on surfaces are structurally stable, providing an {\lq}alternative{\rq} proof to the one given by  Palis in \cite{Pa}, see the corollary below Theorem~\ref{t.cla}.
Of course when two maps are not close to each other, it can be hard to determine whether they are conjugate in the same way 
as it is not immediately obvious whether two knots in $\R^3$ are the same.

Theorem A shows that two Morse-Smale diffeomorphisms are topologically conjugate if and only if their schemes are the same. 
For this reason we say that the scheme is a {\em complete invariant}. Crucially, we show that the scheme requires only finite data.

Theorem B shows that each ``abstract'' admissible  scheme corresponds to an actual Morse-Smale diffeomorphism
in the same way as each rotation number together with some additional information  (about the number of periodic orbits) corresponds to circle diffeomorphism
(uniquely, up to conjugacy), and the same way as every admissible kneading invariant can be ``realised'' within a smooth full family of interval maps, 
see Section II.4 in \cite{MS}.  In the same way Theorem B provides a {\em full  catalogue} of {\em all} Morse-Smale surface diffeomorphisms.

 In this way we suggest that the scheme of a MS surface diffeomorphism should be regarded as the analogue of the rotation number for a MS circle diffeomorphism and of the kneading invariant of a  MS interval map.

Here, as usual, we define a $C^1$ interval or circle map $f$ to be {\em Morse-Smale}, and 
say that $f\in MS([0,1])$ or $f\in MS(S^1)$,  if $f$ has only finitely many critical and periodic points,
if each of its periodic points has multiplier $\lambda\notin  \{0,1,-1\}$ and if no critical point of $f$ 
is eventually mapped onto a critical or periodic point or to an iterate of another critical point
(these latter assumptions are the analogue of transversality of invariant manifolds in the two-dimensional case).  

\begin{rema} A $C^1$ map in $MS([0,1])$ is {\em not} necessarily structurally stable since a nearby diffeomorphism
can have additional  critical points. To obtain structural stability, we need to consider
the space of $C^2$ endomorphisms, and make the additional assumption that 
each critical point is non-degenerate.
\end{rema}

The analogy between the one-dimensional  and two-dimensional case we are referring to 
is summarised in the following table: 

\bigskip 
{\small 
\begin{tabular}{ |  l  |  l | l |    }
\hline  & complete invariant (finite data) & realisation \\
\hline 
MS circle diffeomorphism & rational rotation number + finite data &  no additional conditions \\ 
\hline 
MS interval map  &  kneading invariant + finite data &   admissibility condition  \\ 
\hline 
MS surface diffeomorphism  & scheme  &   admissibility condition\\
\hline 
\end{tabular}
\begin{center}{\small The analogy between the classification of  MS interval and surface maps.}
\end{center}
}
\bigskip 

Indeed, for a MS circle diffeomorphism, the rotation number and the number $\ge 1$
of periodic attractors is a complete invariant; moreover, given to this abstract data corresponds
a Morse-Smale diffeomorphism (unique up to conjugacy).  For piecewise monotone interval maps $f$, 
the analogue of the rotation number is the kneading invariant, defined as follows.
Let $c_1<c_2<\dots<c_d$ be the turning points of $f$, and let $I_0,\dots,I_d$ be 
the components of $I\setminus \{c_1,\dots,c_d\}$ numbered so that $c_i\in \partial I_i\cap \partial I_{i-1}$, $i=1,\dots,d$. 
Then associate to each $c_i$, the sequence defined by $K(c_i)=I_{i_1}I_{i_2}I_{i_3}\dots\in \{I_1,\dots,I_d\}^{\mathbb N_+}$
where $i_k$ is so that $f^k(c_i)\in I_{i_k}$ for $k\ge 1$. The {\em kneading invariant} of $f$ 
is by definition $K(c_1),\dots,K(c_d)$ (there are other more or less equivalent definitions, see \cite{MS}).
Since the kneading invariant does not detect the number of periodic attractors of $f$ additional data is required
to obtain a complete invariant. For this reason we will follow a slightly different approach. 


Let 
 $N\in \mathbb N\cup \{\infty\}$ and define the finite set 
$$P^{N}_f=\bigcup_{c\in C} \bigcup_{k=0}^{N}f^k(c)\cup Per^A_f\mbox{ where }Per^A_f=\{\mbox{attracting periodic points of $f$}\}.$$
The next (essentially well-known) theorem shows that there exists $N<\infty$ so that the conjugacy class of $f$ is fully
determined by the finite set $P^{N}_f$ together with how $f$ acts on this set. 
Much of this can be expressed by the kneading invariant defined above, 
but the set $P^N_f$  is more helpful in this context.

\newtheorem*{theo1D}{Theorem}

\begin{theo1D}[Classification in the one-dimensional case]    For each $f\in MS([0,1])$
there exists $N$ so that $g\in MS([0,1])$ is topological conjugate with $f$ if and only if
there exists an order preserving bijection between 
%
$P_f^N$ and $P_g^N$
so that $h\circ f(x)=g\circ h(x)$ for each $x\in P_f^{N-1}$. 
\end{theo1D}

\begin{proof} 
We say an interval $K$ is a renormalisation interval  of period $n$, if $K,f(K),\dots,f^{n-1}(K)$ have disjoint interiors,
 $f^n(K)\subset K$ and $f^n(\partial K)\subset \partial K$. 
To prove this theorem, we need to 
use the following claim (whose prove can be found in  \cite{MS}[Proposition III.4.2]).
{\bf Claim:} for each $x\in [0,1]$ there exist $n\ge 0$ and  a renormalisation interval $K$ of period $2$, so that 
either  (a) $f^n(x)$ is a fixed point $p$,   (b)  $f^k(x)\to p$ as $k\to \infty$ where $p$ is an attracting 
fixed point or (c) $f^n(x)\in K$.  The assertion  implies in particular that $f\in MS([0,1])$ can only have (finitely many) 
periodic points of periods of the form $\{1,2,\dots,2^b\}$.

For each periodic attractor $p$, take an interval $J=J_p$ with $p\in J\in W^s(p)$ so 
that if $m$ is the period of $p$ then $f^{2m}|J$ is orientation  increasing.
Since $f\in MS([0,1])$  the previous claim implies that each critical point
is in the basin of an attracting periodic point. Hence one  
can choose $N<\infty$ so that if $c,c'\in W^s(p)\cap C(f)$ (where $p$
is a periodic attractor of period $m$)  then there exists $n,n'\le N-2m$ so that $f^n(c),f^{n'}(c')\in J$.
One can arrange is so that if  $f^m|J$ is orientation preserving and  $f^n(c),f^{n'}(c')$ lie
on the same side of $p$ then $f^{n'}(c')$ lies between $f^n(c)$ and $f^{n+m}(c)$
whereas if $f^m|J$ is orientation reversing then $f^n(c')$ lies between $f^n(c)$ and $f^{n+2m}(c)$.

If $g\in MS([0,1])$ and there exists an 
an order preserving bijection $h$ between $P_f^N,P_g^N$
then from the choice of $N$, if  $c,c'\in W^s(p)$ then there exists $n,n'\le N$ so that $f^n(c),f^{n'}(c)$ 
lie in the same fundamental domain of $p$. 
It follows that one can extend $h$ to an order
preserving bijection between  $P^\infty_f,P^\infty_g$. Since  $f,g$ do not have (i) wandering intervals
(this follows from the previous claim), (ii)  intervals consisting of periodic points of constant period, (iii) periodic turning points
which are not attracting (since $f$ is $C^1$ and $f\in MS([0,1])$), it follows by  \cite{MS}[Theorem II.3.1] that $f,g$ are topologically conjugate.
\end{proof}

Conversely, one can derive from  \cite{MS}[Theorem II.5.2], that  given (abstract) finite sets $P'\subset P\subset [0,1]$ and $\pi\colon P'\to P$ with some additional other 
admissibility conditions, 
there exists a $C^\infty$ MS interval map $f$ so that $P^{N-1}_f=P'$ and $P^{N}_f=P$ and so that $f|P'=\pi$. 
So this is the analogue of Theorem B in the interval case. 

%
%
%

\medskip 

Before giving formal statements of theorems  A and B, we will give a historical background of this classification problem, see also the survey \cite{BoGrLa2001} by Bonatti, Grines, Langevin for further details and references.

\bigskip 

{\it Acknowledgements.} This work was supported by the Russian Foundation for Basic Research (project 16-51-10005-Ko\_a), Russian Science Foundation  (project 17-11-01041), the Basic Research Program at the HSE (project 90) in 2017 and the European Union ERC AdG grant 339523 RGDD.

\section{History of the problem of classifying MS diffemorphisms} 

\paragraph{Topological equivalence, roughness and structural stability.}  The concept of {\it roughness} of a dynamical system was conceived in Nizhny Novgorod  (fomely known as Gorky) in 1937.  Andronov and  Pontryagin considered a dynamical system $$\dot{x} = v(x),~~~~~~~~~~~~(1)$$ where $x\in \mathbb R^2$ and $v$  is a $C^1$-vector field   defined in closed compact domain $D\in\mathbb R^2$  bounded  by a smooth  curve  $\partial D$ without selfintersections  transversal to  $v$. They suggested to call  $v$  {\it rough} if for any sufficiently small $C^1$-perturbation of $v$  there exists a homeomorphism $C^0$  {\em close to the identity} which transforms the orbits of the original dynamical system to the orbits of the perturbed system (so the perturbed system is {\it topologically equivalent} to the original one, and the topological equivalence is close to the identity).

In \cite{AP}, they showed that for such a dynamical system  roughness is equivalent to the following two properties: 
\begin{enumerate}[topsep=1pt,itemsep=1pt,parsep=1pt]
\item [1)] the number of the equilibrium points and periodic trajectories is finite and each of these is hyperbolic;
\item [2)]there are no saddle connections (that is stable and unstable separatrices of  any saddle equilibrium state p and of any pairs of different   saddle equilibrium states p and  q have no intersections).
\end{enumerate}
\medskip
The {\it topological classification} (w.r.t. topological equivalence) of structurally stable {\it flows} (dynamical systems with continuous time) on a bounded part of the plane and on the 2-sphere follows from the results by   Leontovich-Andronova and   Mayer obtained in \cite{LeMa1} and \cite{LeMa2}. In fact, in these papers a more general class of dynamical systems was considered.  The classification was based on the ideas of Poincare-Bendixson to pick a set of specially chosen trajectories so that their relative position (the {\it Leontovich-Mayer scheme}) fully determines the qualitative decomposition of the phase space of the dynamical system into the trajectories.

\begin{figure}[h] 
\centering
\includegraphics[width=1.0\textwidth]{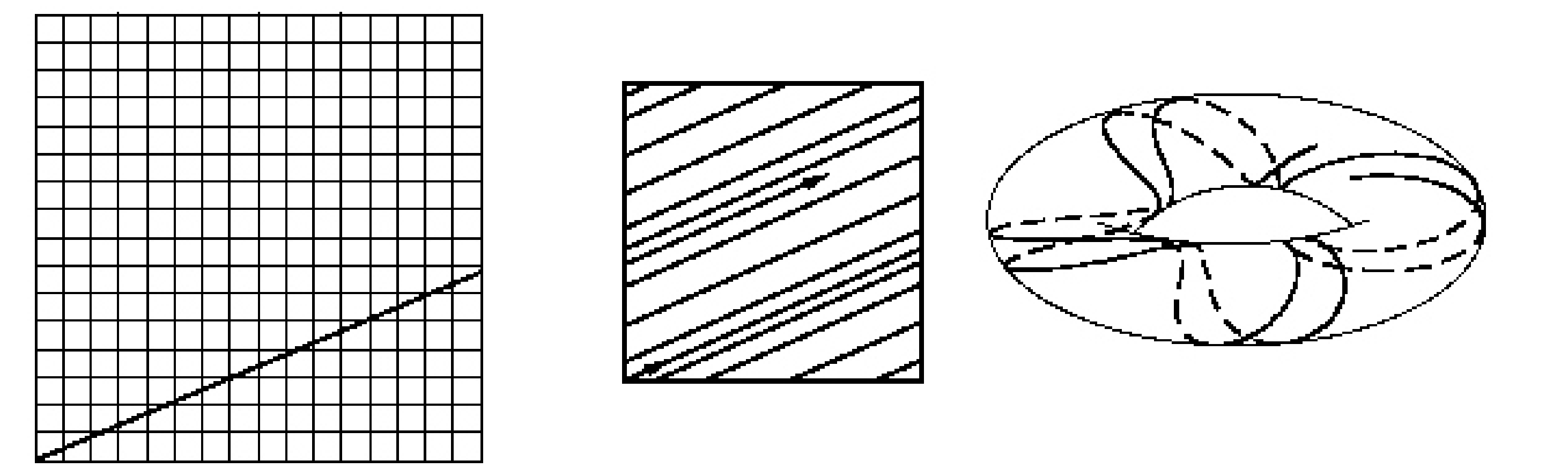}
\caption{\small Irrational winding of the torus}
\label{irr}
\end{figure}

The principal difficulty in generalising this result  to flows on arbitrary orientable surfaces of positive genus is the possibility of new types of trajectories, namely {\it non-closed recurrent trajectories} (see Figure \ref{irr}, where the natural  projection of a straight line with an irrational slope on $\mathbb {R}^2$ to the torus $\mathbb{T}^2$ is shown, which can be a recurrent trajectory of some flow on $\mathbb{T}^2$). The absence of such trajectories for   structurally stable flows without singularities on the 2-torus was first proved by  Mayer in 1939  \cite{Ma}.  Actually, in this  paper he introduced the  notion of roughness for {\em cascades} (i.e., discrete dynamical systems corresponding to $C^1$
{\em diffeomorphisms}) on the {\em circle} and  found necessary and sufficient  conditions for the roughness  these cascades: rough diffeomorphisms of the circle are exactly those that have a finite number of hyperbolic periodic points.

In 1959,  Peixoto \cite{Pe0} introduced the concept of {\it structural stability} of flows to generalize the concept of roughness. A flow $f^t$ is called {\it structurally stable} if for any flow $g^t$ which is $C^1$-sufficiently close to $f^t$, there exists a homeomorphism $h$ sending trajectories of  $g^t$ to trajectories of  $f^t$. The original definition of a rough flow involved the additional requirement that the homeomorphism $h$ is $C^0$-close to the identity map. Peixoto proved that the concepts of roughness and structural stability for flows on surfaces  are equivalent. In 1962,  Peixoto \cite{Pe1,Pe2} proved that the  conditions 1), 2) above plus condition
\begin{enumerate}[topsep=1pt,itemsep=1pt,parsep=1pt]
\item [3)]	all $\omega$- and $\alpha$-limit sets are contained in the union of the equilibrium  points and the periodic trajectories,  
\end{enumerate} 
\noindent
are necessary and sufficient for the structural stability of a flow on any arbitrary  orientable closed (i.e., compact and without boundary) surface and showed that such flows are dense in the space of all $C^1$-flows. 

\paragraph{Morse-Smale systems in arbitrary dimensions.} In 1960,  Smale \cite{S2} introduced a class of diffeomorphisms on manifolds,  generalising the above Andronov-Pontryagin-Peixioto conditions (1), (2), (3)  to the case of diffeomorphism in arbitrary dimensions, requiring that (i) there are at most a finite number of all periodic points, (ii) each periodic orbit is hyperbolic,  and (iii) that each intersection of stable and unstable manifolds is transversal.   This class of systems was also inspired by Smale's earlier work on the Poincar\'e conjecture for manifolds of dimension $\ge 5$,  \cite{Smale60b, Smale61b}, in which he made essential use of Morse theory and vector fields generated by the gradient of Morse functions. Since then  flows and diffeomorphisms with these properties are called  {\it Morse-Smale systems}. Until Smale discovered horseshoe maps, he even assumed that this class of systems
is generic. 

A natural question is the {\em existence} of Morse-Smale systems on closed manifolds.  Smale \cite{Smale61a} proved that any Morse function on the manifold can be approximated by a Morse function whose  gradient vector field is a Morse-Smale flow without periodic orbits. Therefore the time-one map of this flow is a Morse-Smale diffeomorphism. Since Morse functions exist on any closed manifold, it follows that Morse-Smale systems (both flows and diffeomorphisms) exist on any closed manifold. 

In the late 1960's,  Palis \cite{Pa} and  Palis and Smale \cite{PS} proved that  Morse-Smale systems are structurally stable. Therefore, these systems form an open set in the space $C^1$-smooth dynamical systems. From a modern point of view,  Morse-Smale systems on closed manifolds, in the $C^1$ setting, are exactly the structurally stable dynamical systems with zero topological entropy. (This holds because by Ma\~n\'e \cite{Man} any $C^1$ structurally stable system is Axiom A and satisfies the strong transversality condition, and since an Axiom A system with zero topological entropy necessarily has only a finite number of periodic orbits by \cite{Bow}.) From this point of view, they are the simplest structurally stable systems. Anosov \cite{An62, An67} and  Smale \cite{Sm61, Smale65, S3} proved the existence of wide classes of structurally stable dynamical systems with positive topological entropy.

Already in the pioneering work \cite{S2},  Smale established a close relationship between the dynamic characteristics of the Morse-Smale system and the topology of the ambient manifold. 
Later is was found that  the periodic orbits of Morse-Smale flows without equilibrium states, form a rather special set of {\em knots and links}, see Franks \cite{Fr81}, Sullivan \cite{Su}, Wada \cite{Wa}. 
Nevertheless, Morse-Smale vector fields on compact manifolds always have an energy function, as was shown by  Meyer \cite{Mey}.  (By definition, an {\em energy function} for a system is a  Morse-Bott function which decreases along non-periodic orbits and has critical points only at  periodic orbits.)   Pixton \cite{Pi}  showed that this is not true for Morse-Smale diffeomorphisms, namely he proved that for any Morse-Smale diffeomorphism given on a compact  surface there is an energy function and that there is an example of a Morse-€"Smale diffeomorphism on $S^3$ which has no an energy function. Later on, it was shown by Bonatti, Grines, Pixton, Pochinka  \cite{BoGr,Pi,Poch} that there exist Morse-Smale diffeomorphisms in dimension three for which the invariant manifolds of saddle points form wild objects, and hence do not have  an energy function.
Grines, Laudenbach and Pochinka, see \cite{GrLaPo}, established necessary and sufficient conditions  for the existence of an energy function for Morse-€"Smale diffeomorphisms on
 3-manifolds. 

\begin{figure}[h]
\centerline{\includegraphics[height=4cm]{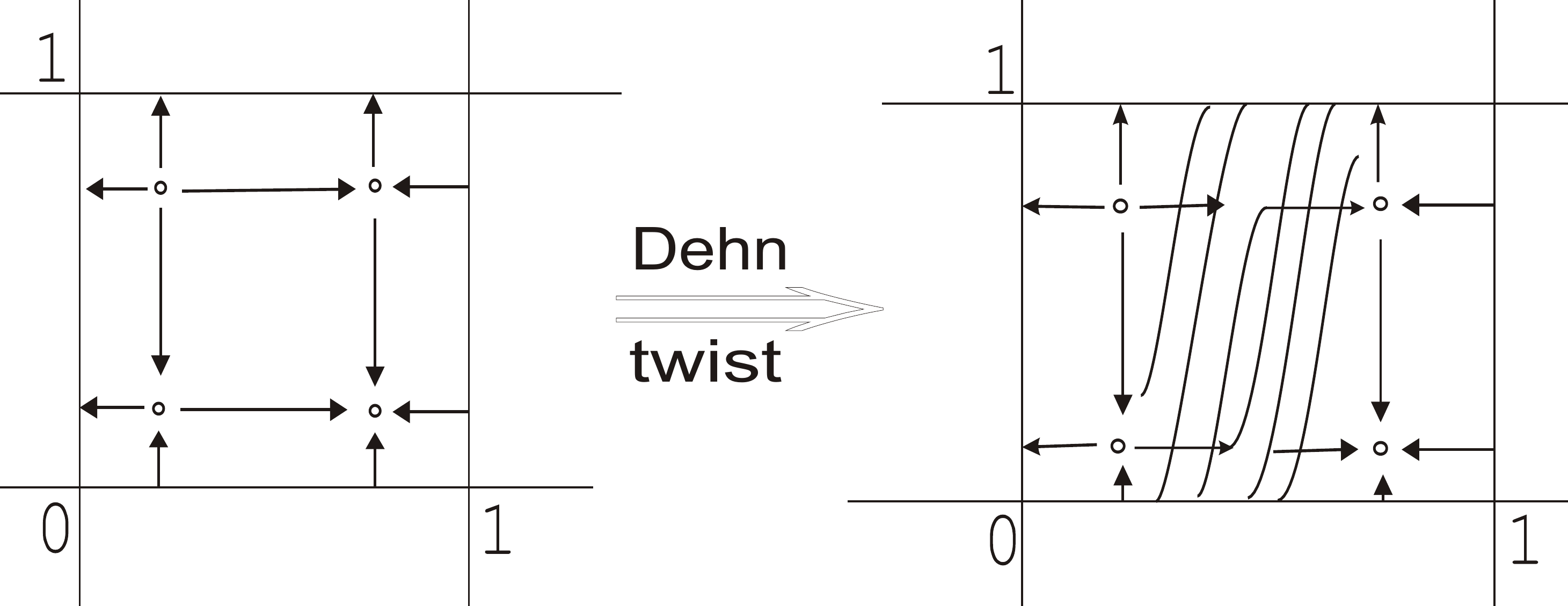}} \caption{{\small A Dehn twist to a Morse-Smale system on $T^2$}} \label{2tor-Den}
\end{figure}

Morse-Smale diffeomorphism which are time-one maps of a flow, obviously induce the identity map in the homology group, which begs the question whether 
Morse-Smale diffeomorphisms  can induce non-trivial isomorphisms in the homology group. The answer is yes, as is easy to see. First, take for example the time-one map of the gradient vector field  on the torus $T^2$ with two saddles and two nodes on the square, see the left side of Figure \ref{2tor-Den}.  Adding a Dehn twist along the closed curve corresponding to $x=\frac{1}{2}$ we get a diffeomorphism whose stable manifold of the first saddle transversely intersects the unstable manifold of the other saddle, see the right side of Figure \ref{2tor-Den}, 
and whose action on the homology group is non-trivial.  However, the action has a special form: 
for arbitrary Morse-Smale diffeomorphism $f: M^d \to M^d$ all eigenvalues of the induced map $ f_ *: H _ * (M ^ d, \mathbb{R}) \to H_*(M ^ d, \mathbb{R}) $  are roots of unity, see works 
\cite{Shub73a,ShubSullivan75} by Shub and Sullivan.

A Dehn twist operation not only leads to a non trivial action in the  homology group but also gives rise to {\it heteroclinic points}, i.e., points of intersection of stable and unstable manifolds of different saddle points. A Morse-Smale diffeomorphism on a closed surface is called {\it gradient-like} if it has no heteroclinic points. Heteroclinic points are an obstruction to the embedding of such a diffeomorphism to a flow.  Palis \cite{Pa} proved that in any neighborhood of the identity map of the surface there is a Morse-Smale diffeomorphism which cannot be embedded in a flow. Moreover he listed necessary conditions for a diffeomorphism to be embedded in a flow, proved that these conditions are sufficient for Morse-Smale diffeomorphisms on surfaces and posed the problem for higher dimensions.  Grines, Gurevich and Pochinka, see  \cite{GrGuPo}, solved this problem 
for Morse-Smale diffeomorphisms on 3-manifolds.

\subsection{Classification of gradient-like diffeomorphisms  on surfaces}\label{subsec:gradient}

In this subsection we describe previous approaches to classify special dynamical systems on surfaces namely  for Morse-Smale flows and for gradient-like diffeomorphisms.

\begin{figure}[h]
\centering
\includegraphics[width=0.55\textwidth]{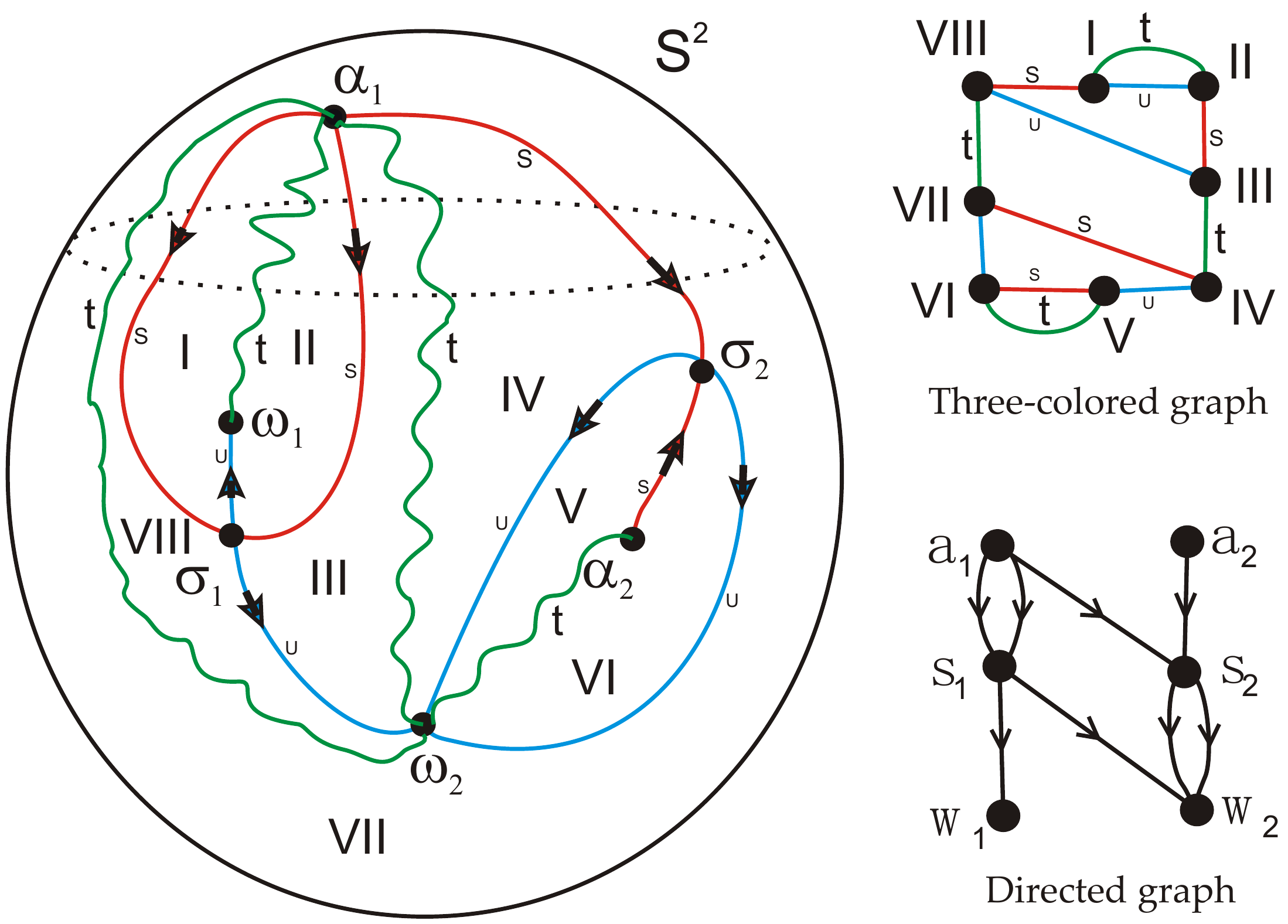}
\caption{\small The directed and three-colour graphs for a gradient-like diffeomorphism on a 2-sphere, see  Subsection~\ref{subsec:gradient}. The regions marked in roman numbers
correpond to vertices in the three color graph, while the stable and unstable manifolds and the curves marked by $t$, correspond to its edges.
The singularities $\alpha_i,\sigma,\omega_i$ correspond to the vertices $a_i,s_i,w_i$ in the graph. }\label{grs}
\end{figure}

\paragraph{A directed graph associated to Morse-Smale flows on surfaces.}  In 1971,  Peixoto \cite{Pe} obtained the classification for  Morse-Smale flows on arbitrary surfaces.
He did this by  generalizing the Leontovich-Mayer scheme for such flows
to a  {\it directed graph}  whose vertices are in one-to-one correspondence with fixed points and closed trajectories of the flow, and whose edges correspond to the connected components of the invariant manifolds of fixed points  and closed trajectories (see Figure \ref{grs}). He proved that the isomorphic class of such directed graph is  a complete topological invariant for Morse-Smale flows on surfaces (where the isomorphisms preserve specially chosen subgraphs).

However, in 1998, Oshemkov and Sharko \cite{OS} pointed out that the Peixoto invariant is unfortunately not complete, 
by giving an example showing that 
an isomorphism of graphs cannot always distinguish between types of decompositions into trajectories for a domain bounded by two periodic orbits of the flow. 
Thus they show that Peixoto's directed graph does not necessarily distinguishes 
non-equivalent systems. They  therefore suggested to use  a {\it three-colour graph}, see Figure \ref{grs}, where vertices correspond to triangular domains and the color ($s,\,t,\,u$) of an edge  corresponds to passing through a side of triangles of the same color. They also showed that this colour graph is a complete invariant. In the next subsection we will explain how such three-colour graphs are obtained.

The construction of the directed and three-colour graphs for Morse-Smale flows is very similar to the construction of similar graphs for gradient-like diffeomorphisms below.  

\paragraph{Directed, equipped and three colour graphs associated to  gradient-like Morse-Smale diffeomorphisms on surfaces.}\label{thr}
As before we say that a Morse-Smale diffeomorphism on a closed surface is called {\it gradient-like} if it has no heteroclinic points.
 In 1985,  Bezdenezhnych and Grines  \cite{BeGr85I,BeGr85II} showed that 
for gradient-like diffeomorphisms on surfaces a directed graph with an automorphism  is again a complete invariant.  In 2014, Grines, Kapkaeva and Pochinka \cite{GrKaPo14},  showed that two gradient-like diffeomorphisms on a surface are topologically conjugate if and only if their three-colour graphs equipped  by periodic automorphisms are isomorphic 
(Grines, Malyshev, Pochinka and Zinina \cite{GrMaPoZi} describe an  efficient algorithm for distinguishing such graphs). Let us describe the above graphs in more detail.

Let $f$ be a gradient-like diffeomorphism of an orientable surface $M^2$. 
Let $\sigma$ be a saddle point of $f$ of period $m_\sigma$. Denote by $\nu_\sigma$ the {\it type of orientation of $\sigma$}, which is $1$ if the diffeomorphism  $f^{m_\sigma}|_{W^u_\sigma}$ preserves orientation and $-1$ otherwise. Let  $l^s_\sigma~(l^u_\sigma)$ be a stable (unstable) separatrix of  $\sigma$, i.e., $l^s_\sigma~(l^u_\sigma)$  is a connected component of $W^s_\sigma\setminus\sigma~(W^u_\sigma\setminus\sigma)$. Since $l^s_\sigma~(l^u_\sigma)$ does not intersect the unstable (stable) manifold of any saddle point, there exists a sink $\omega$ (respectively, a source $\alpha$) such that $cl(l^u_\sigma)=l^u_\sigma\cup\sigma\cup\omega$ (respectively, $cl(l^s_\sigma)=l^s_\sigma\cup\sigma\cup\alpha$); see, for example, Lemma 3.2.1 in the monograph \cite{grin} by Grines, Pochinka. For  $\delta\in\{s, u\}$ we define the direction of the separatrix $l^\delta_\sigma$ to be towards the saddle point if $\delta = s$ and from the saddle point if $\delta= u$.

\begin{figure}[h!]
\centerline{\includegraphics[scale=0.6]{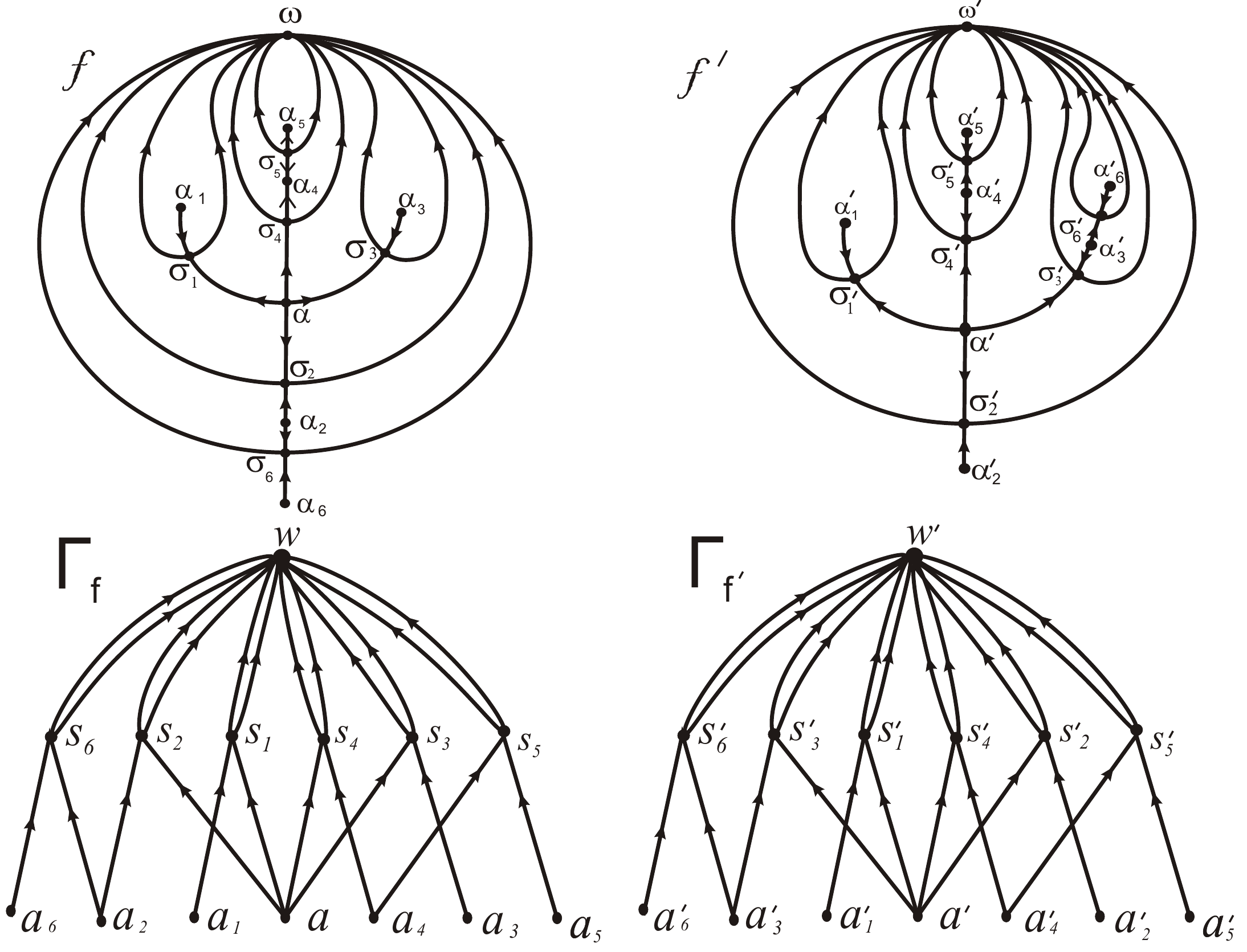}}
\caption{\small The diffeomorphisms $f,f':\mathbb S^{2}\to\mathbb
S^{2}$ have isomorphic graphs but they are not topologically conjugate as their equipped graphs $\Gamma_f^{*}$ are not isomorphic (see the explanations in the text).}\label{Expe}
\end{figure}

We say that a {\it directed graph $\Gamma_f$  is the graph of a diffeomorphism} $f$, see Figure~\ref{grs}, if
\begin{enumerate}[topsep=1pt,itemsep=1pt,parsep=1pt]
\item[1)] the {\em vertices} of the graph $\Gamma_f$ correspond to the periodic points of the non-wandering set $\Omega_f$;  the vertex corresponding to the saddle periodic point $\sigma$ we equip with the value $\nu_\sigma$;
\item[2)] the {\em directed edges }of the graph $\Gamma_f$ correspond to the directed separatrices of the saddle points. 
\end{enumerate}
\noindent 
The diffeomorphism $f$ naturally induces an automorphism $f_*$ of  the graph  $\Gamma_f$. Let $\Gamma_f, ~\Gamma_{f'}$ be the graphs of diffeomorphisms $f,~f'$. 
The existence of an isomorphism between the graphs $\Gamma_f$ and $\Gamma_{f'}$ conjugating the automorphisms $f_*$ and $f'_*$ is necessary for the topological conjugacy of the diffeomorphisms $f,f'$. Unfortunately, in general the existence of an isomorphism of the graphs is not sufficient for the maps $f,f'$ to be conjugacy even if every periodic points is a fixed point and each separatrix is $f$-invariant.
Indeed,  consider diffeomorphisms $f$ and $f^\prime$ with phase portraits shown in Figure \ref{Expe}. Even though these diffeomorphisms have isomorphic graphs, they are not topologically conjugate.
To see this, notice that any conjugating homeomorphism necessarily carries the basin of the sink  $\omega$   of the diffeomorphism $f$ into the basin of the sink $\omega^\prime$ of the diffeomorphism $f^\prime$. However, such a homeomorphism cannot be extended to the entire sphere in such a way that it would carry the invariant manifolds of the saddle points of $f$ into the invariant manifolds of the saddle points of  $f^\prime$.

So the directed graph $\Gamma_f$ of a diffeomorphism $f$ does not determine the topological conjugacy class of $f$, and we are required to add information to $\Gamma_f$.
To obtain a complete classification for gradient-like diffeomorphisms on surfaces,  Bezdenezhnych and Grines  \cite{BeGr85I,BeGr85II} 
introduced in 1985 the notion of {\em equipped graphs}.  Such equipped graphs can be defined as follows. Let $\omega$ be a sink of  $f$ and let $L_\omega$ be the subset of the manifold $M^2$ that consists of the separatrices which have $\omega$ in their closures. 
Then there is a smooth 2-disk $B_\omega$ such that $\omega\in B_\omega$ and each separatrix  $l\subset L_\omega$ intersects $\partial B_\omega$ at an unique point; see, for example,  \cite[Proposition 2.1.3]{grin}. For the vertex $w$ corresponding to the periodic sink point $\omega$,  let $E_w$  denote the set of edges of the directed graph $\Gamma_f$ incident to $w$. Let $N_w$ denote the cardinality of the set $E_w$. We enumerate the edges of the set $E_w$ in the following way. First we pick in the basin of the sink $\omega$ a 2-disk $B_\omega$  and set $c_\omega=\partial B_\omega$. We define a pair of vectors $(\vec{\tau}, \vec n)$ at some point of the curve $c_\omega$ in such a way that the vector $\vec n$ is directed inside the disk $B_\omega$, the vector $\vec{\tau}$ is tangent to the curve $c_\omega$ and induces  a counterclockwise orientation on $c_\omega$ with respect to $B_\omega$ (we call this orientation  positive). Enumerate the edges  $e_{1},\dots,e_{N_w}$ from $E_w$ according to the ordering of the corresponding separatrices as we move along $c_\omega$ starting from some point on $c_\omega$. This enumeration of the edges of the set $E_w$ is said to be {\it compatible} with the embedding of the separatrices.

The graph $\Gamma_f$ is said to be {\it equipped} if each vertex $w$ is  numbered with respect to the numeration of the edges of the set $E_w$ and  the numeration is compatible with the embedding of the separatrices. We denote such a graph by $\Gamma^{*}_f$. For an example, see Figure~\ref{Expe}.

Let  $\Gamma^{*}_f$ and $\Gamma^{*}_{f^{\prime}}$ be equipped graphs of diffeomorphisms $f$ and $f^{\prime}$ respectively and let $\Gamma^{*}_f$ and $\Gamma^{*}_{f^{\prime}}$ be isomorphic by an isomorphism $\xi$. Let a vertex $w$ of the graph $\Gamma^{*}_f$ correspond to a sink and let $w'=\xi(w)$. 
Then the isomorphism $\xi$ induces the permutation $\Theta_{w,w'}$ on $\{1,\dots,N\}$ (where $N=N_w=N_{w'}$) defined by $\Theta_{w,w'}(i)=j\Leftrightarrow \xi(e_i)=e'_j$.

Two equipped graphs $\Gamma^{*}_f$,
$\Gamma^{*}_{f^{\prime}}$ of diffeomorphisms $f$, $f^{\prime}$
are said to be {\it isomorphic} if there exists an isomorphism $\xi$ of the graphs
$\Gamma_f$, $\Gamma_{f'}$ such that

\begin{enumerate}[topsep=1pt,itemsep=1pt,parsep=1pt]
\item[1)] $\xi$ sends the vertices into the vertices and preserves the values of the vertices corresponding to the saddle periodic points; it sends the edges into the edges and preserves their direction;
\item[2)] the permutation $\Theta_{w,w'}$ induced by $\xi$ is a power of a cyclic permutation\footnote{It is directly checkable that the property of the permutation to be a power of a cyclic permutation is independent of the choice of the curves $c_\omega$ and $c_{\omega'}$.} for each vertex $w$ corresponding to a sink;
\item[3)] $f'_*={\xi}f_* \xi^{-1}$.
\end{enumerate}

\noindent
The equipped graph $\Gamma^*_f$ of a diffeomorphism $f$ is again a topological invariant up to isomorphism. Let us show that the equipped graphs $\Gamma_f^*, \Gamma_f'^*$ of $f,f'$ are not isomorphic.   To do this, consider Figure \ref{Expe} once again and suppose that the vertex $w~(w')$ of the graph corresponds to the sink point $\omega~(\omega')$. One can check directly that any isomorphism $\xi$ induces the permutation $\Theta_{w,w'}$ which is not a power of a cyclic permutation and therefore the equipped graphs  $\Gamma^{*}_f$, $\Gamma^{*}_{f^{\prime}}$ are not isomorphic. Generalising this argument, Bezdenezhnych and Grines  \cite{BeGr85I,BeGr85II} were able to show that the equipped graphs are a complete invariant for gradient-like Morse-Smale diffeomorphisms on closed surfaces.

Let us now describe an alternative complete invariant for gradient-like diffeomorphisms,  namely the {\em three-colour graph}. Let, as before,  $f$ be a gradient-like Morse-Smale diffeomorphism on a closed surface $M^2$. The non-wandering set $\Omega_f$ can be represented as  $\Omega_f=\Omega^0_f\cup\Omega^1_f\cup\Omega^2_f$, where $\Omega^0_f$, $\Omega^1_f$, $\Omega^2_f$ denote the set of sinks, saddles, and sources of the diffeomorphism $f$, respectively. Throughout the remainder of this subsection,  we assume that $f$ has at least one saddle point\footnote{If a Morse-Smale diffeomorphism $f:M^n\to M^n$ has no saddle points, then its nonwandering set consists of one source and one sink. All diffeomorphisms ``source-sink'' are topologically conjugate; the proof of this fact is given, for example, in \cite{grin} (Theorem 2.2.1).}.

\begin{figure}[h!]
\centerline{\includegraphics[scale=0.35]{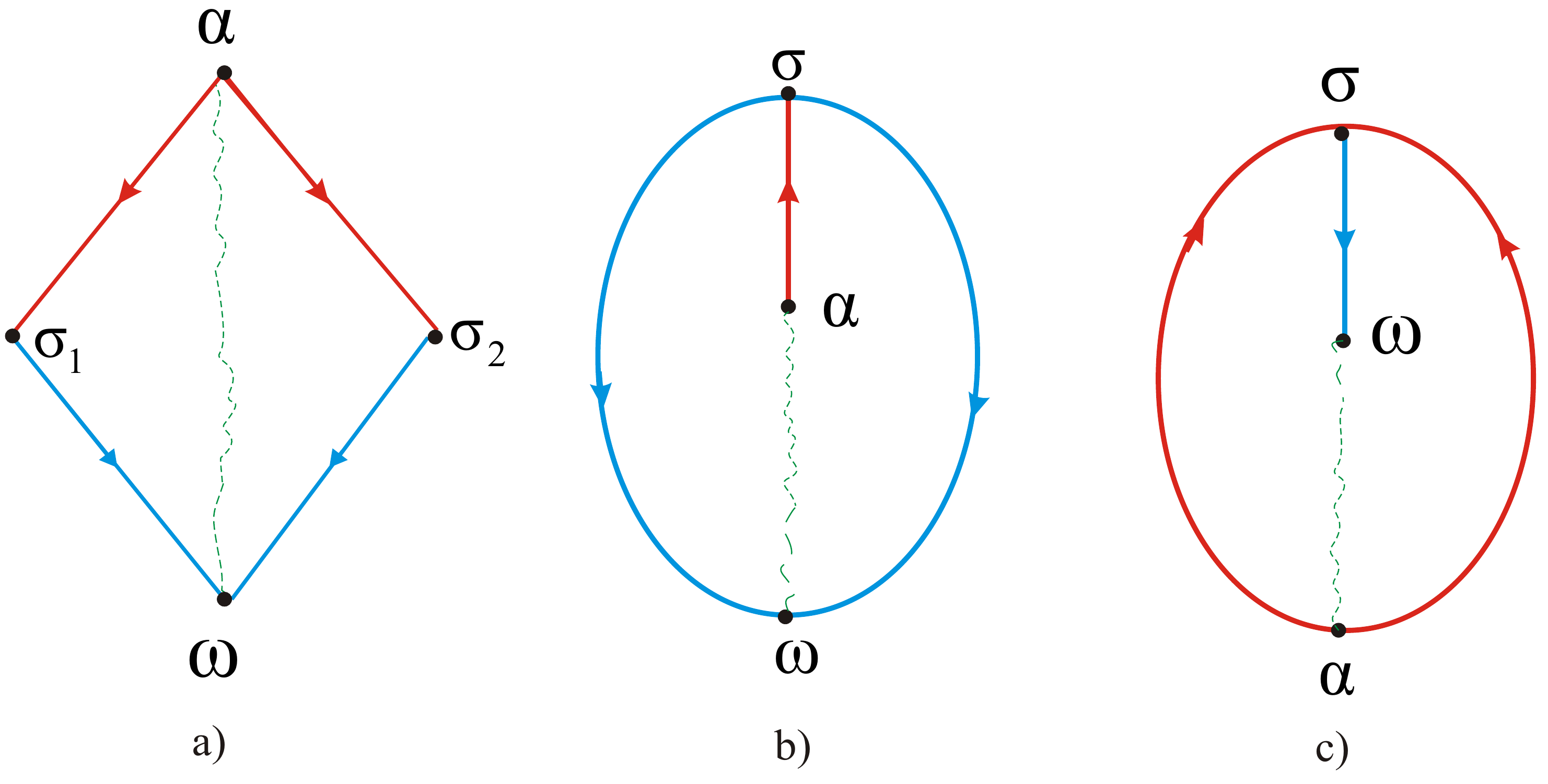}} \caption{\small Types of cells with $t$-curves}\label{ris:division yacheiki}
\end{figure}

We remove from the surface $M^2$ the closure of the union of the stable manifolds and the unstable manifolds of all the saddle points of the diffeomorphism $f$ and let the resulting set be denoted by $\tilde{M}$, that is, $\tilde{M}=M^2\setminus(\Omega_f^0\cup W^u_{\Omega_f^1}\cup W^s_{\Omega_f^1}\cup \Omega_f^2)$. The set $\tilde{M}$ is represented in the form of a union of domains ({\it cells}) homeomorphic to the open two-dimensional disc such that the boundary of each of these cells has one of the forms depicted by boldface lines in Figure \ref{ris:division yacheiki} and it contains exactly one source, one sink, one or two saddle points, and some of their separatrices.

Let $A$ be any cell from the set $\tilde {M}$, and let $\alpha$ and $\omega$ be the source and sink contained in its boundary. A simple curve $\tau\subset A$ whose boundary points are the source $\alpha$ and the sink $\omega$ is called a {\it $t$-curve} (see Figure \ref{ris:division yacheiki}). Let $\mathcal T$ denote a set which is invariant under the diffeomorphism $f$ and which consists of $t$-curves taken one from each cell.

\begin{figure}[h!]
\centerline{\includegraphics[width=6cm,height=4.5cm]{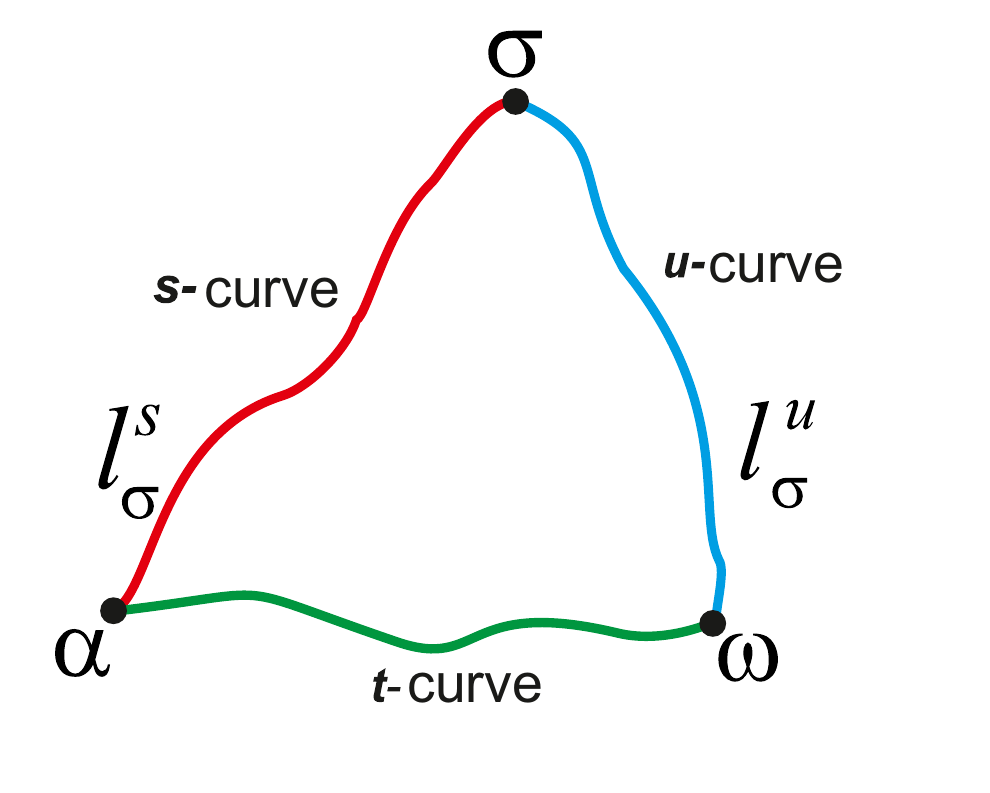}}\vskip -0.5cm \caption{\small Triangular domain}
\label{ris:triangle}
\end{figure}

Any connected component of the set $M_\Delta=\tilde M\setminus\mathcal T$ is called {\it a triangular area}. Let $\Delta_f $ denote the set of all triangular domains of diffeomorphism $f$. The boundary of every triangular domain $\delta \in\Delta_f$ contains three periodic points: a source $\alpha$, a saddle $\sigma$, a sink $\omega$. It contains also the stable separatrix $l^s_{\sigma}$ (called the {\it $s$-curve}) with $\alpha$ and $\sigma$ as boundary points, the unstable separatrix $l^u_{\sigma}$ (called the {\it $u$-curve}) with $\omega$ and $\sigma$ as boundary points and a curve $\tau$ (a $t$-curve) with  $\alpha$ and $\omega$ as boundary points  (see Figure \ref{ris:triangle}). {\it A triangular domain} is bounded by $s$-, $u$- and $t$-curves. We say that two triangular areas {\it have a common side}, if this side belongs to the closures of both domains. {\it The period of the triangular domain} $\delta $ is defined to be the least positive integer $k\in\mathbb{N}$, such that $f^k(\delta) =\delta$.

\begin{figure}[h!]
\centerline{\includegraphics[scale=0.6]{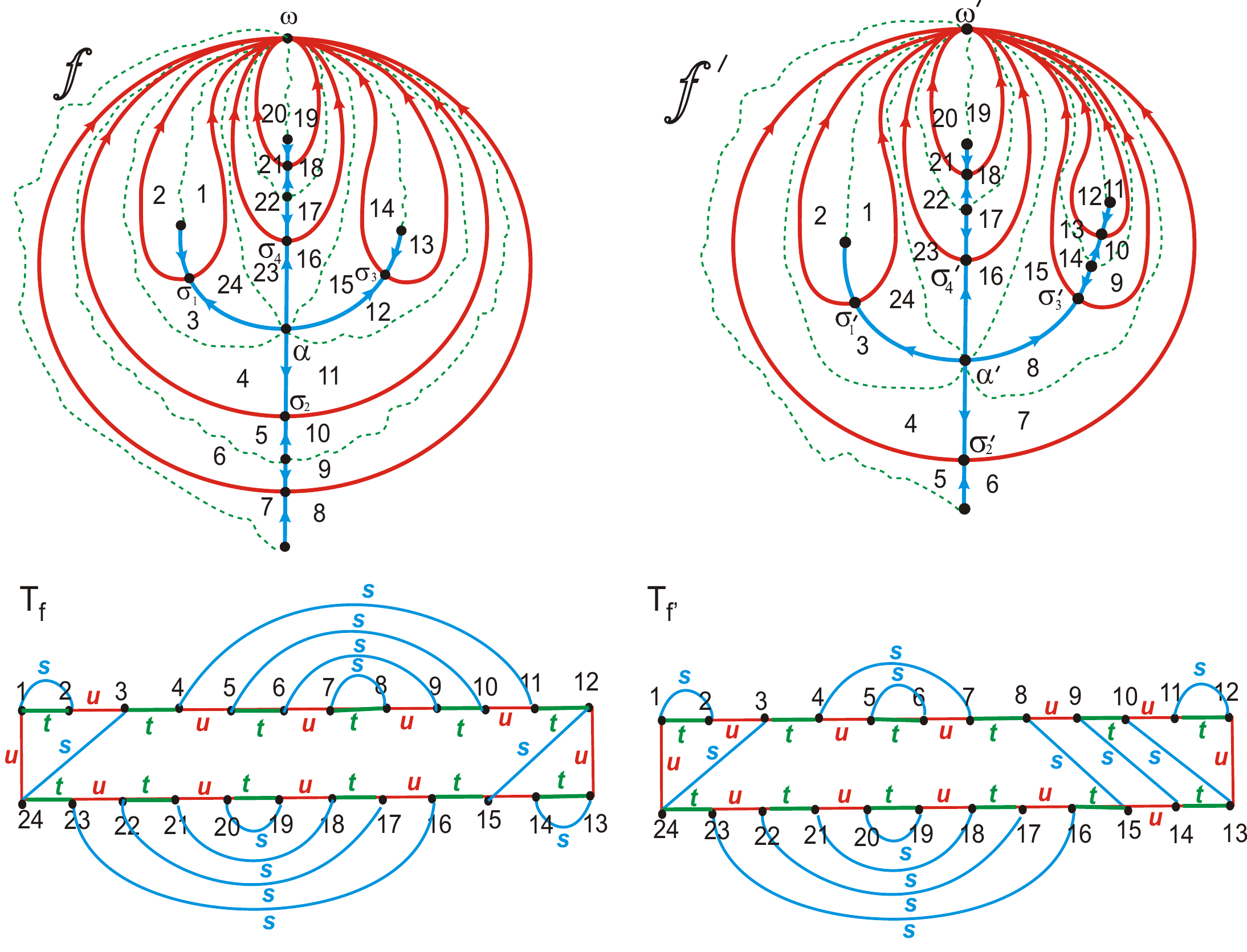}} \caption{\small The non-isomorphic three-colour graphs $T_f,T_{f'}$ associated to the non conjugated gradient-like diffeomorphisms $f,f'$.}
\label{ris:ex}
\end{figure}

We construct a {\em three-color} ($s,t,u$) graph $T_f$, corresponding to a Morse-Smale gradient-like diffeomorphism $f$ as follows (see Figure \ref{ris:ex}):

\begin{enumerate}[topsep=1pt,itemsep=1pt,parsep=1pt]
\item[1)] the vertices of $T_f$ are in a one-to-one correspondence with the  triangular domains of the sets $\Delta$;
\item[2)] two vertices of the graph are incident to an edge of color $s$, $t$ and $u$, if the corresponding triangular domains have a common $s$, $t$ and $u$-curve.
\end{enumerate} 
By construction, three-color graphs obtained from different partitions into triangular domains (depending on the choice of $t$ -curves) are isomorphic.

Let $B_f$ denote the set of the vertices of the graph $T_f$. Since the sides of any triangular domain are assigned different colors, edges of three different colors come together at the vertex corresponding to the triangular domain. Since any side of a triangular domain is adjacent to exactly two different triangular domains, the graph $T_f$ has no cycles of length 1. Thus, the graph $T_f$ satisfies the definition of the three-color graph. Let $\pi_f: \Delta_f \to B_f$ denote a one-to-one map of the set of the triangular domains of the diffeomorphism $f$ into the set of the vertices of the graph $T_f$. The diffeomorphism $f$ induces the automorphism $f_* = \pi_f f \pi_f ^{- 1}$ on the set of vertices  of the graph $T_f$. Let $(T_f,f_*)$ be denote the three colour graph $T_f$ together with the automorphism $f_*$.

Two three-color graphs with automorphisms $(T_f, f_*)$ and $(T_{f'},f'_*)$ of diffeomorphisms $f,f'$ are said to be {\it isomorphic} if there exists a one-to-one correspondence $\xi $ between the sets of their vertices which preserve the relations of incidence and the color, as well as the conjugating automorphisms $f_*$ and $f_*'$ (that is, $f_*'=\xi f_* \xi^{-1}$).

As mentioned, in \cite{GrKaPo14} it is shown that the three-color graph $(T_f, f_*)$ of a diffeomorphism $f$ is a {\em complete} topological invariant up to isomorphism for gradient-like Morse-Smale diffeomorphisms on closed surfaces. 

In this paper we shall associate a different complete invariant to a (general) surface Morse-Smale diffeomorphism, namely a scheme.  For the diffeomorphisms from Figures \ref{Expe} \ref{ris:ex}
the invariants associated to each of the two diffeomorphisms consists of a torus with 10 closed curves, see  Figure \ref{12curves}. Since these curves are ordered differently on the tori for the two diffeomorphisms, it follows at once that such systems are not topologically conjugated.

\begin{figure}[h!]
\centerline{\includegraphics[width=0.65\textwidth]{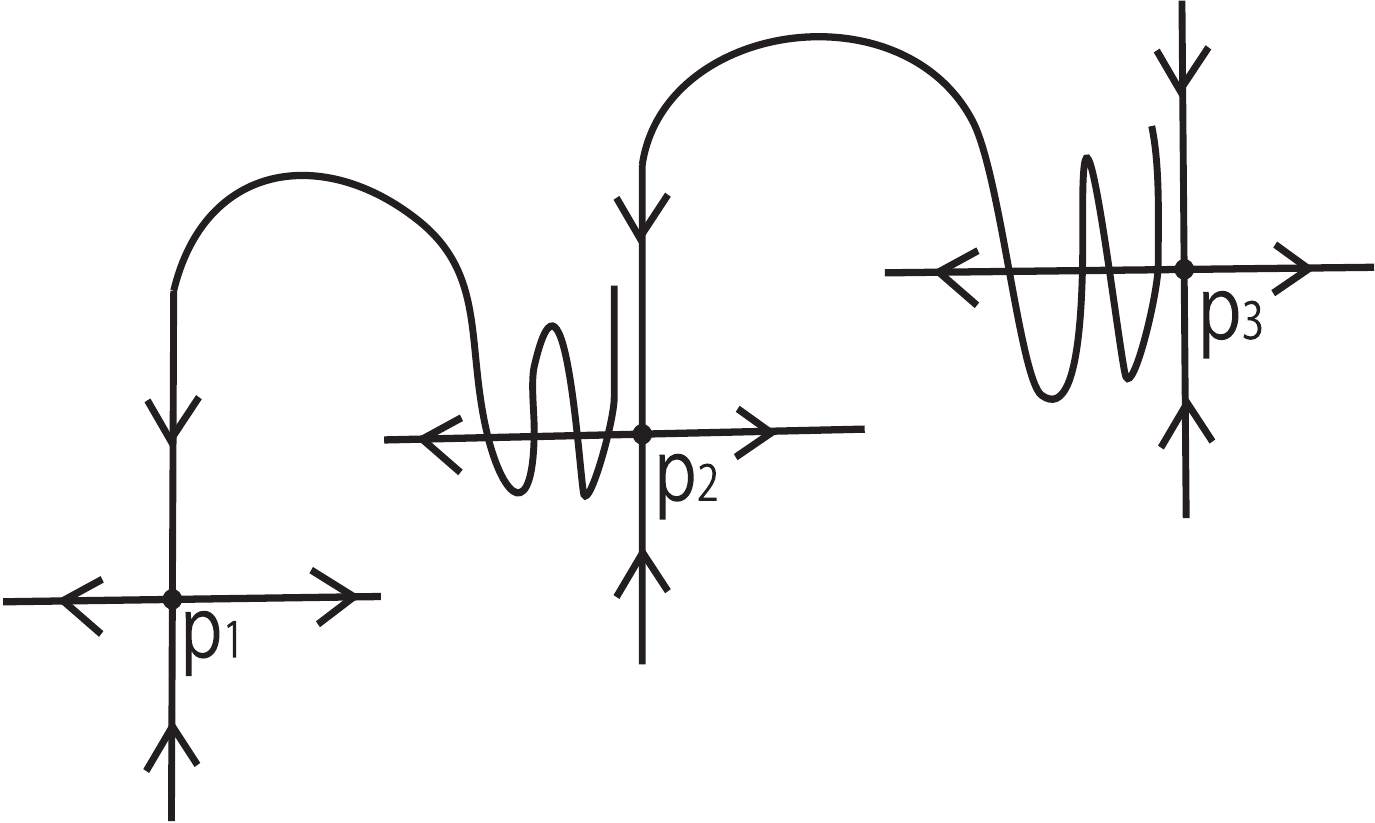}}
\caption{\small A chain of the length 3 \label{chain3}}\centering
\end{figure}

\subsection{Classification of non-gradient-like diffeomorphisms on closed surfaces}

Let us now discuss previous results concerning the classification of general Morse-Smale diffeomorphisms on closed surfaces. 
Let $\mathcal O_i,~\mathcal O_j$ be periodic orbits of Morse-Smale diffeomorphism $f:M^2\to M^2$.  Smale \cite{S3} introduced a {\it partial order relation $\prec$} for the periodic orbits
$$\mathcal O_i\prec\mathcal O_j\iff W^s_{\mathcal O_i}\cap W^u_{\mathcal O_j}\neq\emptyset .$$
A sequence of distinct periodic orbits $\mathcal O_i=\mathcal O_{i_0},\mathcal O_{i_1},\dots,\mathcal O_{i_k}=\mathcal O_{j}$ ($k\geq 1$), such that  $\mathcal O_{i_0}\prec\mathcal O_{i_1}\prec\dots\prec\mathcal O_{i_k}$  is called a {\it chain of length $k$ joining the periodic orbits  $\mathcal O_i$ and $\mathcal O_j$}. The maximal length of the chain joining $\mathcal O_i$ and $\mathcal O_j$ is denoted by $$beh(\mathcal O_j\vert\mathcal O_i)$$ ({\em beh} stands for {\em behaviour}).
where we define  $beh(\mathcal O_j\vert\mathcal O_i)=0$ when $W^u_{\mathcal O_j}\cap W^s_{\mathcal O_i}=\emptyset$. For a subset $P$ of the periodic orbits let us set $beh(\mathcal O_j|P)=\max\limits_{\mathcal O_i\subset P}\{beh(\mathcal O_j|\mathcal O_i)\}$.
Figure \ref{chain3} gives an example where $\mathcal O_1\prec\mathcal O_2\prec\mathcal O_3$ for saddle fixed points $p_1=\mathcal O_1,p_2=\mathcal O_2,p_3=\mathcal O_3$ and $beh(\mathcal O_2|\mathcal O_1)=beh(\mathcal O_3|\mathcal O_2)=1,~beh(\mathcal O_3|\mathcal O_1)=2$. 
Set $$beh(f)=max\{beh(\mathcal O_j\vert\mathcal O_i)\}.$$ 
A Morse-Smale diffeomorphism has $beh(f)=1$ iff it is a so-called {\em sink-source diffeomorphism}. When $beh(f)=2$
then $f$ has no heteroclinic points, and so is {\em gradient-like}. 
A Morse-Smale diffeomorphism $f:M^2\to M^2$ with $beh(f)>2$ has a chain of saddle orbits of the length $beh(f)-2$.

\begin{figure}[h]
\centerline{\includegraphics[height=5cm]{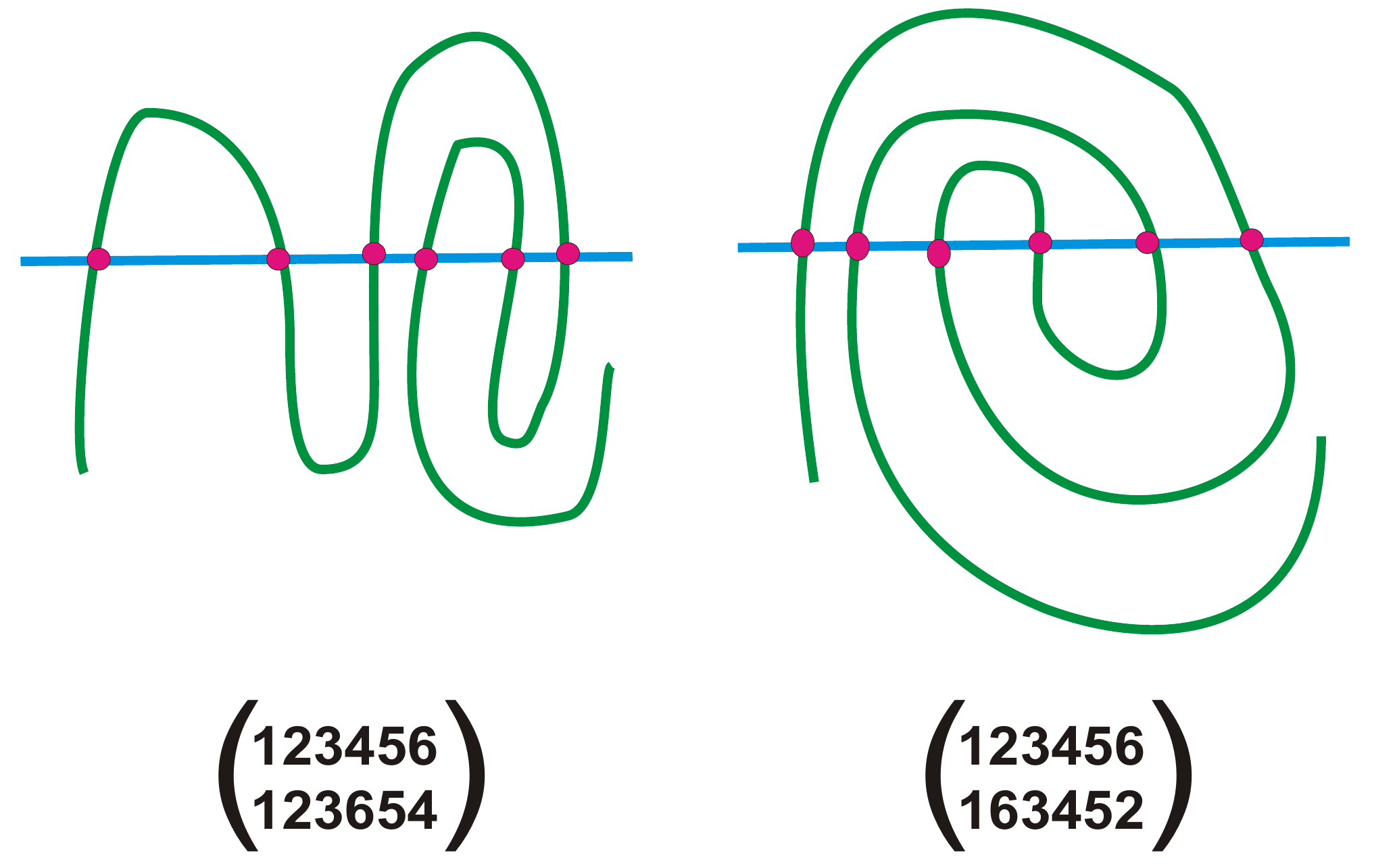}} \caption{\small Heteroclinic substitution}\label{tf20}
\end{figure}

 If $beh(f)=3$ then $f$ has a finite number of heteroclinic orbits.  In 1993,   Grines \cite{Gr97} proved that for such diffeomorphisms an invariant similar to Peixoto's graph carrying additional information on the  {\em heteroclinic substitution}, describing the intersection pattern
 of invariant manifolds as in Figure \ref{tf20}, is sufficient to describe necessary and sufficient conditions for topological conjugacy. 

\begin{figure}
\centerline{\includegraphics[width=11cm,height=7cm]
{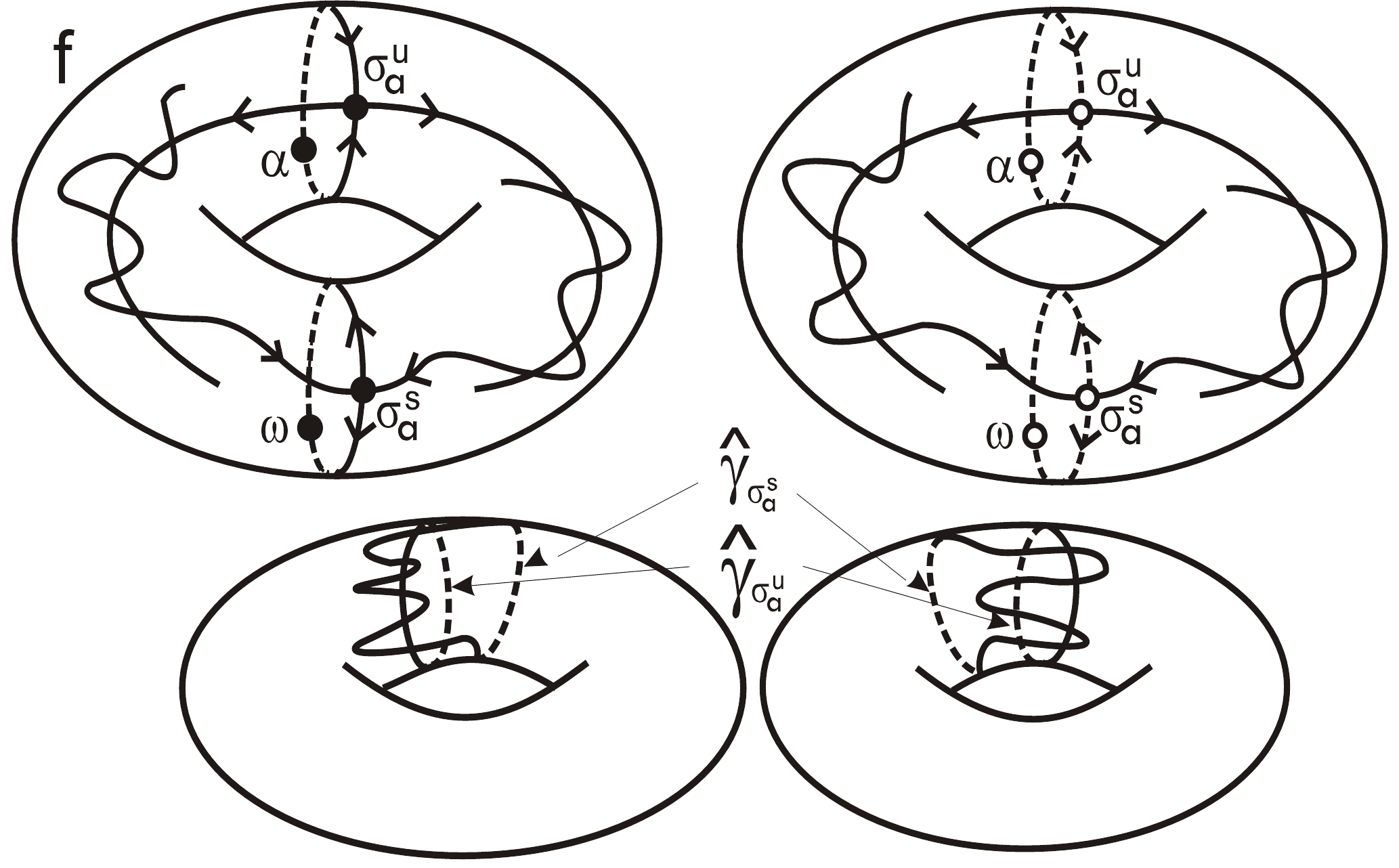}}\caption{\small Scheme of a Morse-Smale diffeomorphism with $beh(f)=3$}\label{sh}
\end{figure}

In 1993, Langevin \cite{La} proposed to consider the {\em orbit space} of the basin of the sink and project to this closed surface the unstable separatries of the saddle points. This approach was generalized and successfully applied by  Bonatti,  Grines,  Medvedev,  Pecou and  Pochinka  in \cite{BoGrMePe}, \cite{BoGrPo} for the topological classification of Morse-Smale  diffeomorphisms $f$ with $beh(f)\leq 3$ on 3-manifolds. In 2010,  Mitryakova and  Pochinka  \cite{MiPo2010} applied this method to the topological classification of Morse-Smale diffeomorphisms $f$ with $beh(f)\leq 3$ on orientable surfaces. Indeed, they constructed a topological invariant (which they called a ``{\em scheme}'') which consists of a finite number of two-dimensional tori (corresponding to the orbit space of the basin of sinks and sources), together with a set of simple closed curves (corresponding to the orbit spaces of separatrices), see  Figure \ref{sh}. They also proved that this invariant is complete when $beh(f)\leq 3$. In 2013, Mitryakova and  Pochinka \cite{MiPo2013}  solved  the  {\em realization problem} for such diffeomorphisms, establishing that each {\em abstract scheme} corresponds to a Morse-Smale diffeomorphism (we will make this notion more precise below). Vlasenko \cite{Vl98} in 1998 proprosed another approach to the topological classification for arbitrary structurally stable diffeomorphisms of orientable surfaces using an equipped oriented graph whose vertices  corresponds to a periodic and heteroclinic point and each directed edge corresponds to a connecting segment of separatrices.
 
In 1998, a different approach was taken   Bonatti and  Langevin \cite{BoLa}, who considered $C^1$-structurally stable diffeomorphisms (Smale diffeomorphisms) of compact surfaces. The main result of that paper consists of a finite combinatorial presentation of the global topological dynamics in terms of the geometrical types of certain Markov partitions of the hyperbolic sets and by gluing the domains along their boundary. One important step of the proof of their theorem consists of a precise analysis of the topological position (the {\lq}intersection pattern{\rq}) of invariant manifolds of the Smale diffeomorphisms.

Let us describe their construction for a Morse-Smale diffeomorphism $f:M^2\to M^2$.  Consider a maximum heteroclinic chain $ O_1 $, $ \ldots $, $ O_h $ of saddle orbits 
for  $ f $ and define $K_ {1h} = \left (\bigcup\limits_ {i = 1} ^ hW ^ u_{O_i} \right) \cap \left (\bigcup\limits_ {i = 1} ^ hW ^ s_{O_i} \right)$. The set $\bigcup\limits_{i = 1} ^ hO_i \cup K_ {1h} $ is called a \textit{saturation} of the chain. So the saturation of the chain is the union of saddle orbits and various intersections of stable and unstable separatrices of saddle periodic points belonging to the chain. Note that the maximum saturation of  heteroclinic chain is an invariant set $K_f$. By definition, separatrices of saddle periodic points intersect transversely for Morse-Smale diffeomorphisms. Using this fact, they proved that one can put a uniform hyperbolic structure
on $K_f$ and that a neighbourhood of $K_f$ has {\it finite topological type}, i.e., is a closed surface with a finite number of holes.

The saturated hyperbolic set can be viewed as a generalization of the notion of a basic set of saddle type. (By definition, a basic set $\Lambda$ is a compact invariant set which carries a hyperbolic structure and so that $f\colon \Lambda\to \Lambda$ is transitive) This makes it possible to apply the Bowen-Sinai technique for constructing Markov partitions for such sets \cite{Sinai68,Bow71}. In other words,  $K_f$ is covered by a special family of curvilinear quadrangles formed by segments of the stable and unstable separatrices, which they case a {\em good Markov partition}. The geometric type of a good Markov partition 
includes a description of the mutual arrangement, orientation and numbering of curvilinear quadrangles, and their images under the action of the diffeomorphism. Two geometric type are called \textit{equivalent} if they are geometric types of some good Markov partitions of the same hyperbolic
saturated set. The main result in \cite{BoLa} in the setting of Morse-Smale systems is the following: let $ K_ {f_1} $, $ K_ {f_2} $ be the saturated hyperbolic sets of Morse-Smale diffeomorphisms $ f_1 $, $ f_2 $ on closed oriented surfaces $ M_1 $, $ M_2 $ respectively, suppose that $ K_ {f_1} $ and $ K_ {f_2} $ have good Markov partitions with equivalent geometric types, then $f_1 $, $ f_2 $ are topologically conjugate on invariant neighborhoods of the sets $K_{f_1}$, $ K_{f_2} $. See Theorem 1.0.3 in \cite{BoLa} for a more precise statement. 

In a  subsequent work, Beguin \cite{Beguin} developed a finite algorithm to decide the equivalence of two realizable geometrical types. In fact, as is shown in \cite{BoLa},  some of the abstract geometrical types do not correspond to any Smale diffeomorphisms on compact surfaces. 

In the next section we will introduce alternative, and in our opinion more natural,  combinatorial objects (called {\lq}schemes{\rq})  which, unlike the previous approaches, not only give a complete classification (as in Theorem A), but also describe which Morse-Smale diffeomorphisms can be realised (as in Theorem B).  

\section{The scheme of a MS surface diffeomorphism is a complete invariant}\label{subsec:statementresults}

The aim of this paper is to give a different classification which is inspired by the one for circle diffeomorphisms. The upshot is that we obtain a {\em  finite amount of data} about a Morse-Smale diffeomorphism which is {\em necessary and sufficient} to determine whether or not it is topologically conjugate to another Morse-Smale diffeomorphism.  One of the main features
of this classification is that, {\em  given this abstract data}, we will always be able {\em construct} a diffeomorphism that realises this data. 

More precisely, in the present paper we consider class $MS(M^2)$ preserving orientation Morse-Smale diffeomorphisms of an orientable surface $M^2$ and give a complete topological classification within this class, and solve the corresponding realization problem, by means of topological invariants similar to those used in  \cite{BoGrMePe,BoGrPo,MiPo2010,MiPo2013}. To make this precise we need to introduce some notions. 

\subsection{Maximal  $u$-compatible system of linearizing neighborhoods:  statement of Theorem~\ref{add}}  \label{subsec:maxdef}
 Let $f\in MS(M^2)$. We will assume that $beh(f)>1$ because all sink-source diffeomorphisms are topologically conjugate. For a periodic point $p$ of $f$ let us denote by $k_p$ the period of $p$ and by $\nu_p$ the orientation type of  $p$, i.e  $\nu_p=+~(-)$ if the map $f^{k_p}|_{W^u_p}$ preserves (reverses) orientation. Denote by $\Sigma$ the set of all saddle points of $f$. In the construction below,  linearizing neighbourhoods of saddles will play a crucial role. It turns out to be convenient to choose  {\lq}canonical{\rq} linearizing neighborhoods. For this aim, for $\nu\in\{+,-\}$, let us introduce a {\it canonical diffeomorphism} $a_\nu:\mathbb R^2\to\mathbb R^2$ by the formula  $$a_\nu(x_1,x_2)=(\nu 2x_1,\frac{\nu x_2}{2}).$$ Notice that the origin $O=(0,0)$ is a saddle point with the unstable manifold $W^u_O=Ox_1$ and the stable manifold $W^s_O=Ox_2$. Define $$N=\{(x_1,x_2)\in\mathbb R^2:|x_1\cdot x_2|\leq 1\}\mbox{ and } N^t=\{(x_1,x_2)\in\mathbb R^2:|x_1\cdot x_2|\leq t\}\mbox{ for } t\in(0,1].$$ 
By construction, $N$ is $a_\nu$-invariant neighborhood of $O$. We will say that $N$ is a {\it canonical neighborhood}. Denote by $F^u$ the one-dimensional foliation which consists of the leaves $\{(x_1,x_2)\in N:x_2=c\}$, $c\in\mathbb R$ and by $F^s$ the one-dimensional foliation which consists of the leaves $\{(x_1,x_2)\in N:x_1=c\}$, $c\in\mathbb R$.  

\begin{figure}[h!]
\centerline{\includegraphics[width=5.5cm,height=5.5cm]{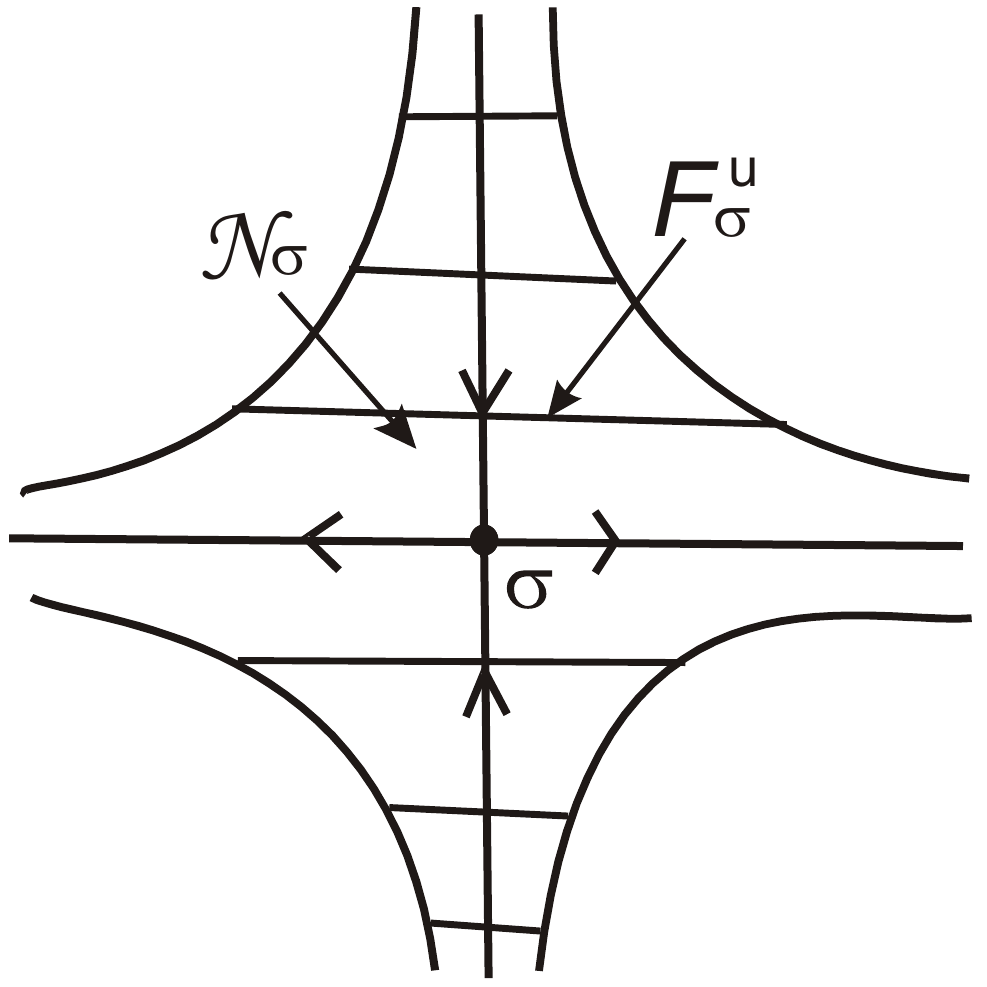}}\caption{\small An $u$-linearizing neighborhood} \label{5}
\end{figure}

\begin{defi}[The $u$-linearizing neighborhood] \label{adop} Let $\sigma$ be a saddle periodic point for $f$. A neighborhood $N_\sigma$ of $\sigma$ together with a one-dimensional foliation ${F}^u_{\sigma}$ containing $W^u_\sigma$ as a leaf, is called \emph{$u$-linearizable} if there is a homeomorphism $\mu_\sigma:N\to{N}_\sigma$ which conjugates the  canonical diffeomorphism $a_{\nu_\sigma}\vert_{{N}}$ to the diffeomorphism  $f^{k_\sigma}\vert_{{N}_\sigma}$ and sends leaves of the foliation ${F}^u$ to leaves of the foliation $F^u_{\sigma}$ (see Figure \ref{5}).
\end{defi}

For every point $x\in N_{\sigma}$ denote by ${F}^u_{\sigma,x}$ the unique leave of the foliation  ${F}^u_{\sigma}$ passing through the point $x$.

\begin{figure}[h!]
\centerline{\includegraphics[width=11cm,height=8cm]{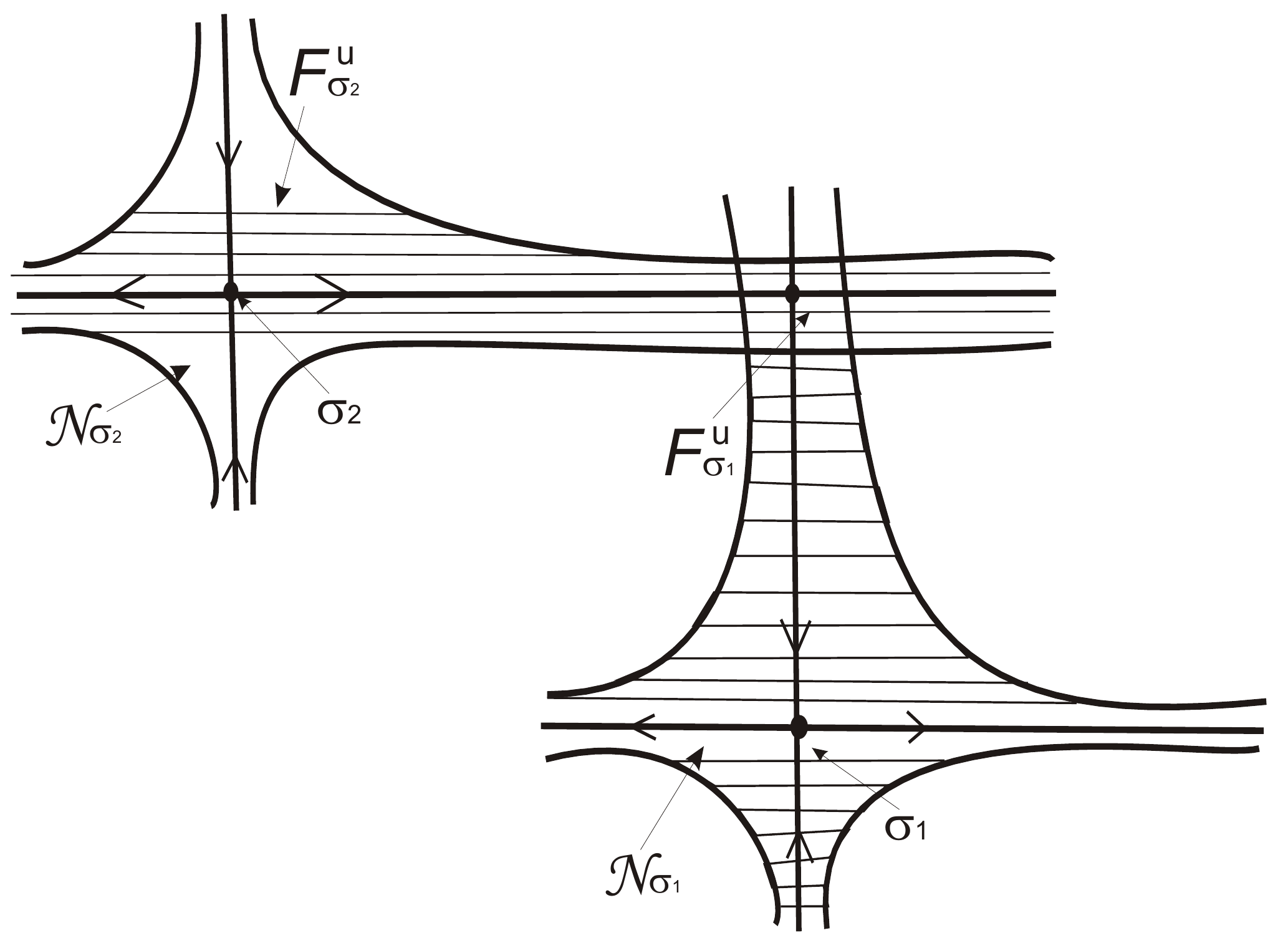}}\caption{\small An $u$-compatible system of neighbourhoods} \label{a}
\end{figure}

\begin{defi}[A {\em $u$-compatible system of neighbourhoods}] \label{dopsystem} An $f$-invariant collection $N_f$ of $u$-linearizable neighborhoods $N_\sigma$ of all saddle points $\sigma\in\Sigma$ is called \emph{$u$-compatible} if the following properties are hold:

1) $\mu_\sigma(\partial{N})$ does not contain heteroclinic points for any $\sigma\in\Sigma$; 

2) if $W^{s}_{{\sigma_1}}\cap W^{u}_{{\sigma_2}}=\emptyset$ and $W^{u}_{{\sigma_1}}\cap W^{s}_{{\sigma_2}}=\emptyset$ for $\sigma_1,\sigma_2\in\Sigma$ then ${N}_{{\sigma_1}}\cap{N}_{{\sigma_2}}=\emptyset$;

3) if $W^s_{\sigma_1}\cap W^u_{\sigma_2}\neq\emptyset$ for $\sigma_1,\sigma_2\in\Sigma$ then $({F}^u_{{\sigma_1},x}\cap {N}_{{\sigma_2}})\subset{F}^u_{{\sigma_2},x}$ for  $x\in(N_{{\sigma_1}}\cap N_{{\sigma_2}})$ $($see Figure \ref{a}$)$.
\end{defi}

It will be proved in Proposition \ref{comp} that there are $u$-compatible neighbourhoods for every diffeomorphism $f\in MS(M^2)$. Indeed these are slight modifications of the admissible system of tubular families constructed by  Palis and  Smale in \cite{Pa} and \cite{PS}. It is easy to see (see, for example Figures \ref{faz+} and \ref{s1}) that in general there is no conjugating homeomorphism which sends an $u$-compatible system of neighborhoods for $f$ to the $u$-compatible system of neighborhoods for $f'$, even when $f$ and $f'$ are topologically conjugated. Therefore we need a more meaningful notion.     

For a saddle point $\sigma$ let us denote by $[a,b]^u_\sigma~([a,b]^s_\sigma)$ the segment of $W^u_\sigma~(W^s_\sigma)$ situated between the points $a,b\in W^u_\sigma~(a,b\in W^s_\sigma)$.

\begin{figure}[h!]
\centerline{\includegraphics[width=17cm,height=4.5cm]{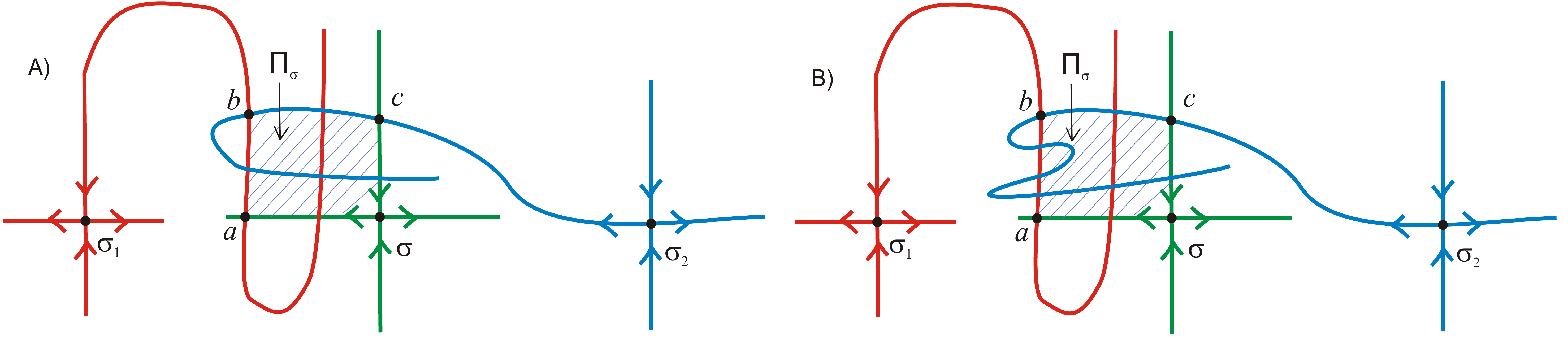}}\caption{\small A) $\Pi_\sigma$ is a heteroclinic rectangle. B) $\Pi_\sigma$ is not a heteroclinic rectangle} \label{rec}
\end{figure}

\begin{defi}[{\em Heteroclinic rectangle}] A closed 2-disc $\Pi_\sigma$ bounded by segments  $[\sigma,a]^u_\sigma,~[a,b]^s_{\sigma_1},~[b,c]^u_{\sigma_2},~[c,\sigma]^s_\sigma,~\sigma_1,\sigma_2\in\Sigma$ and such that $int\,\Pi_\si\cap\Om_f=\emptyset$  is called {\em a heteroclinic rectangle} with respect to $\sigma$ if every connected component of the set $W^s_\Sigma\cap\Pi_\sigma$ intersects every connected component of the set $W^u_\Sigma\cap\Pi_\sigma$ at exactly one point $($see Figure \ref{rec}$)$.  
\end{defi} 

\begin{defi}[The  {\em maximal $u$-compatible system of neighbourhoods}] \label{dopsystem+} A $u$-compatible system of neighbourhoods $N_f$ is called {\em maximal} if every linearizing neighborhood $N_\sigma\in N_f$ contains each heteroclinic rectangle $\Pi_\sigma$. 
\end{defi}

\begin{theo} \label{add} For every diffeomorphism $f\in MS(M^2)$ there is a maximal $u$-compatible system of  neighbourhoods.
\end{theo}

The proof of this theorem will be given in Section~\ref{maxx}. Everywhere below we will denote by $N_\sigma$ a   linearizing neighborhood  of a saddle point $\sigma$ from a maximal $u$-compatible system of neighbourhoods.

\begin{figure}[h!]
\centerline{\includegraphics[width=17cm,height=6.5cm]{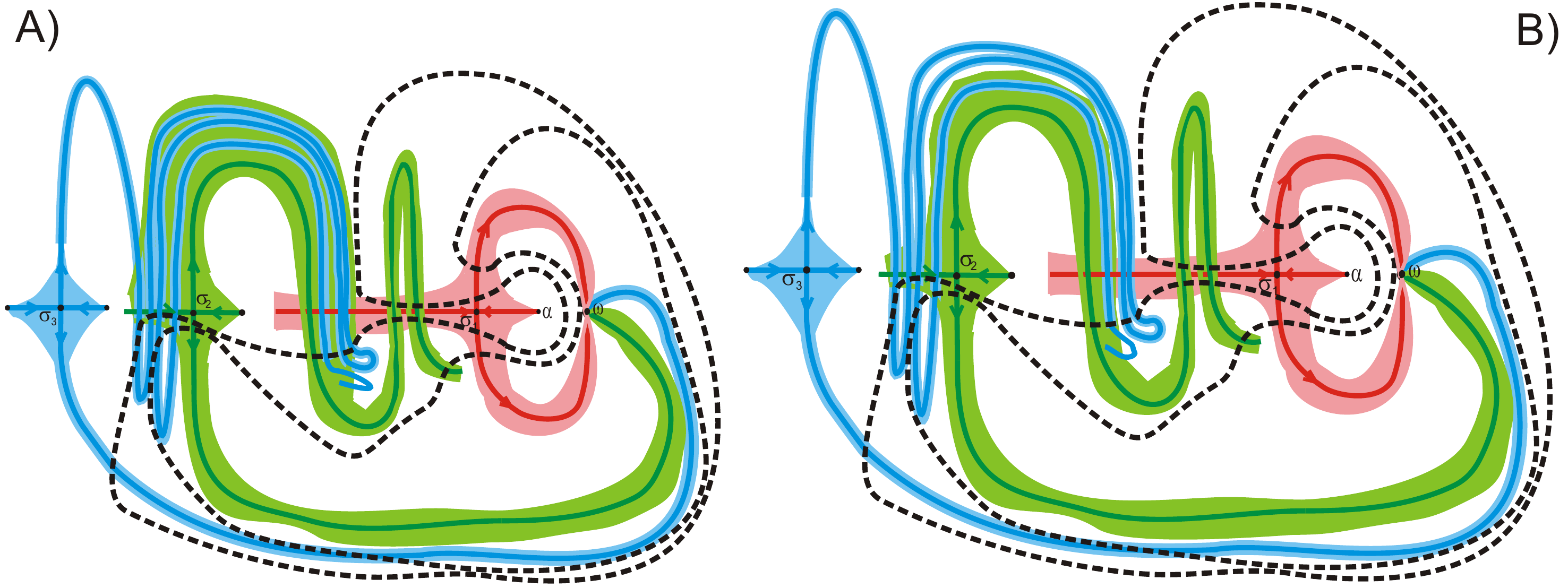}}\caption{\small For a Morse-Smale diffeomorphism on the 2-sphere whose phase portrait is given in Figure \ref{faz+} we represent: A) A maximal $u$-compatible system; B) $u$-compatible system which is non-maximal.} \label{faz}
\end{figure}

\subsection{Associating a  scheme to a Morse-Smale surface diffeomorphism} 
Let $\Sigma_0$ be the set of all sinks of $f$. Let us decompose the set $\Sigma$ of all saddle periodic points of $f$ into subsets $\Sigma_1,\dots,\Sigma_{beh(f)-1}$ inductively as follows: define $\Sigma_i$ to be the set of all saddle points of $f$ such that $beh(\mathcal O_\sigma|\Sigma_{i-1})=1$ for each orbit $\mathcal O_\sigma,~\sigma\in\Sigma_i$.  Let $\Sigma_{beh(f)}$ be the set of all sources.

\paragraph{The quotient space $\mathcal V_f$ of the stable manifold of sinks.}
Let 
\begin{equation}
\mathcal V_f=W^s_{\Sigma_0}\setminus\Sigma_0,~~~\hat{\mathcal V}_f=\mathcal V_f/f.
\label{eq:Vf}
\end{equation}
Since $f$ is a diffeomorphism, the natural projection $p_{_f}:\mathcal V_f\to\hat{\mathcal V}_f$ is a covering map. Every connected component $\hat{V}_j$ of  $\hat{\mathcal V}_f$ is homeomorphic to a 2-torus and corresponds to a unique periodic orbit $\mathcal O_{\omega_j}$ of a sink $\omega_j$. Indeed, the factor space $\hat{\mathcal V}_f=\mathcal V_f/f$ is obtained by taking a fundamental annulus in the basin of each attracting  periodic orbit and identifying its boundaries by $f^{m_j}$ where $m_j$ is the {\em period of this sink}. 

\begin{figure}[h!]\label{factor-tor}
\centerline{\includegraphics[width=14cm,height=6cm]{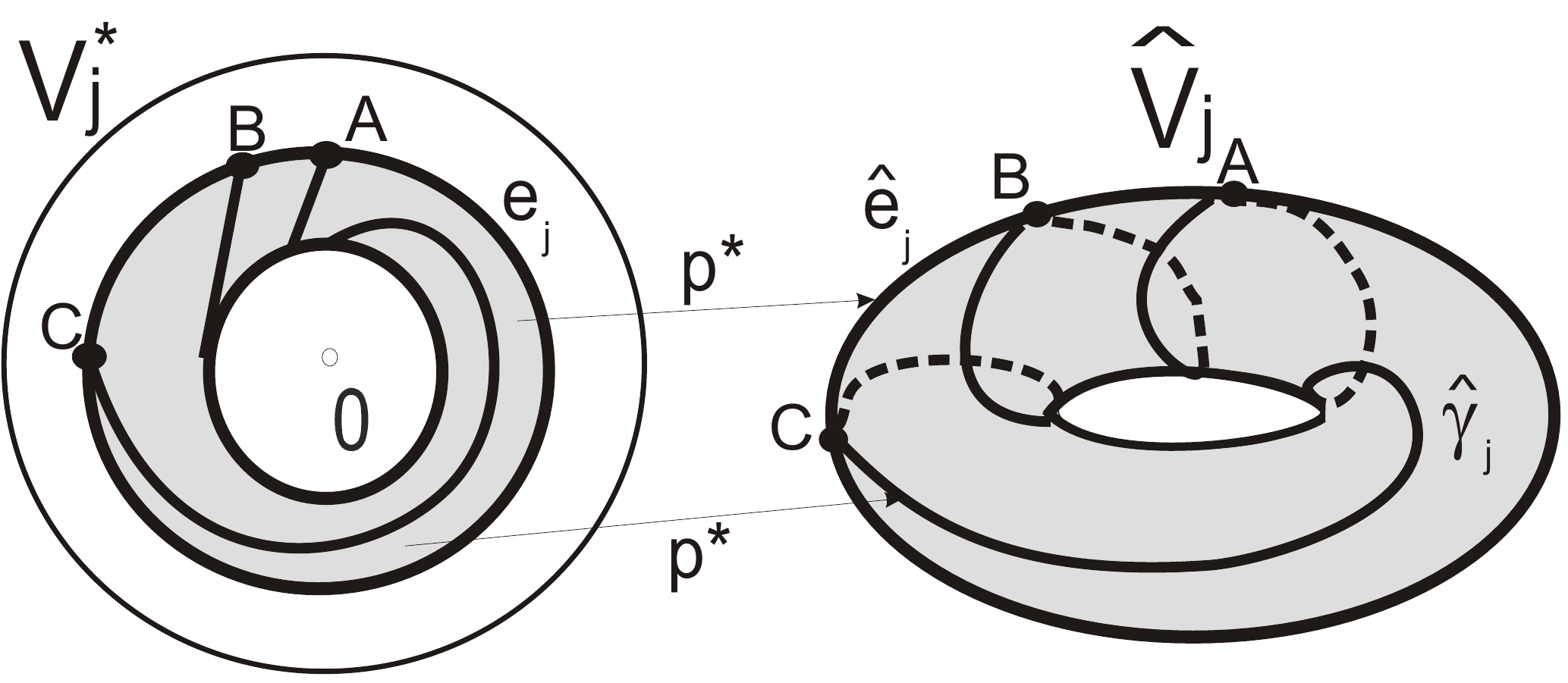}}\caption{\small  $V_j^*$ is represented when $m_j=1$, as is the  equator $\hat e_j$ and some meridian $\hat y_j$ on $\hat V_j$. The curve $e_j\subset V_j^*$ is a connected component of $(p^*)^{-1}(\hat e_j)$. {The curve $\hat\gamma\subset\hat V_j$ winds three times in $\hat V_j$, its preimage $(p^*)^{-1}(\hat\gamma_j)$ consists of three curves passing through the points $A,B,C$, accordingly.}} 
\end{figure}

\paragraph{Equators on the connected components $\hat{V}_j$ of  $\mathcal V_f$.}
It will be helpful  to choose a particular generator of $\pi_1(\hat{V}_j)$ by defining
an epimorphism $$\eta_{\hat{V}_j}:\pi_1(\hat{V}_j)\to m_j\mathbb Z$$ as follows, where $m_j$ 
is the period of the the sink $\omega_j$, as follows. Take the homotopy class $[\hat c]\in\pi_1(\hat{V}_j)$ of a closed curve $\hat c\colon \R/\Z\to \hat{V}_j$. Then $\hat c\colon [0,1]\to \hat{V}_j$ lifts to a curve $c\colon [0,1] \to \mathcal V_f$ connecting a point $x$ with a point $f^{n}(x)$ for some multiple $n\in \Z$ of $m_j$, where $n$ is independent of the lift. So define $\eta_{\hat{V}_j}[\hat c]=n$. A simple closed curve $\hat e_j$ on $\hat V_j$ is called an {\em equator} if $\eta_{{\hat V}_j}[\hat e_j]=0$ and $[\hat e_j]\ne 0$ (see Figure \ref{factor-tor}). Note that the equator $\hat e_j$ is uniquely determined up to homotopy by $\eta_j$. Therefore $\eta_{\hat{V}_j}$ is uniquely  determined by the integer $m_j\ge 1$ and the equator $\hat e_j$ on ${\hat V}_j$. In this way we obtain a unique morphism $$\eta_{_f}:\pi_1(\hat{\mathcal V}_f)\to \mathbb Z$$ so that $\eta_{_f}|\pi_1(\hat{V}_j)=\eta_{\hat{V}_j}$ for each component $\hat{V}_j$ of $\hat{\mathcal V}_f$. 
We will say that a closed curve $\hat\gamma$ {\em winds $n\in \mathbb N$ times 
in $\hat V_j$} if $\eta_{\hat{V}_j}[\hat\gamma]=n \cdot m_j$, see Figure~\ref{factor-tor}.

\paragraph{The covering space $\hat V_j$ of  $\hat{V}_j$.}
Given a component $\hat V_j$ of $\hat{\mathcal V}_f$, the set $p^{-1}(\hat V_j) \subset {\mathcal V}_f$
is equal to $W^s_{\mathcal O_{\omega_j}}\setminus \mathcal O_{\omega_j}$ and therefore homeomorphic to 
$V_j^*:=(\mathbb R^2\setminus O)\times\mathbb Z_{m_{j}}$
where as before $m_j$ is the period of $\omega_j$.  Let $V^*$ be the union of all $V_j^*$ and let $p^*\colon V^* \to \hat{\mathcal V}_f$ be the covering map corresponding to $p_f$. If $\hat e_j$ is an equator on   $\hat V_j$, then $(p^*)^{-1}(\hat e_j)$ is the countable union of simple  closed curves in each of the sets $(\mathbb R^2\setminus O)\times \{k\}  \subset V_j^*$. The complement of these curves form annuli
in $(\mathbb R^2\setminus O)\times \{k\}$, and gluing the two boundaries of such an annulus (according to $f^{m_j}$) we obtain again the torus $\hat V_j$ (see Figure \ref{factor-tor}). For convenience, in the figures below, $(\mathbb R^2\setminus O)$ is drawn as a punctured disc $(\mathbb D^2\setminus O)$. 

\paragraph{The maximal linearizing neighborhoods as subsets of $\mathcal V_f$.} 
For each $i\in\{0,1,\dots,beh(f)\}$, let  $$W^u_i=W^u_{\Sigma_i},~~W^s_i=W^s_{\Sigma_i}.$$ 
For each $i\in\{1,\dots,beh(f)-1\}$, let $$\mathcal N_{i}=\bigcup\limits_{\sigma\in\Sigma_i}N_{{\sigma}}$$ 
be the corresponding linearizing neighbourhoods.

\paragraph{The scheme of a diffeomorphism.} 
For each $i\in\{1,\dots,beh(f)-1\}$, $\sigma\in\Sigma_i$, let 
$$U_\sigma=N_\sigma\setminus \bigcup\limits_{j=1}^{i-1}\mathcal N_j,\,\,\,\hat U_\sigma=p_{_f}(U_\sigma),\,\,\,U^*_\sigma=(p^*)^{-1}(\hat U_\sigma).$$ Let $$\hat{\mathcal U}_i=\bigcup\limits_{\sigma\in\Sigma_i}\hat U_\sigma,\,\hat{\mathcal U}_f=\bigcup\limits_{i=1}^{beh(f)-1}\hat{\mathcal U}_i.$$ 
Let us define $$S_f=(\hat{\mathcal V}_f,\eta_f,\bigcup\limits_{i=1}^{beh(f)-1}\{\hat{U}_\sigma\}_{\sigma\in \Sigma_i}).$$ 

\begin{defi}[The scheme of a diffeomorphism] We call the triple $S_f$
the {\em scheme} associated to the diffeomorphism  $f\in MS(M^2)$. 
\end{defi}

\begin{defi}[Equivalence of schemes] \label{def:equivschemes}
The schemes $S_f\mbox{ and }S_{f'}$ of the two diffeomorphisms  $f, f' \in MS(M^2)$, respectively, are said to be {\em equivalent} if there exist an orientation preserving homeomorphism $\hat\varphi:\hat{\mathcal V}_f\to\hat{\mathcal V}_{f'}$  
such that:

(1) $\eta_{f'}\hat\varphi_*=\eta_{f}$;

(2) $\hat\varphi(\hat{\mathcal U}_f)=\hat{\mathcal U}_{f'}$, moreover for every $i=1,\dots,beh(f)-1$ and every point $\sigma\in\Sigma_i$ there is a point $\sigma'\in\Sigma'_i$ 
such that $\varphi^*(U^*_\sigma)=U^*_{\sigma'}$, where $\varphi^*:{\mathcal V}_f^*\to{\mathcal V}_{f'}^*$ is the lift of $\hat\varphi$.
\label{eqv}
\end{defi}

\begin{rema} It will be proved in Lemma \ref{mama} below that the  equivalence class of a scheme $S_f$  does not depend on a choice of the maximal system of the neighborhoods.
\end{rema}

\begin{rema}
Property (1) in this definition can be restated by requiring that $\varphi$ sends equators of $\hat{\mathcal V}_f$  to equators of $\hat{\mathcal V}_{f'}$, and that the integer  $m_j$ associated 
to a component $\hat{V}_j$ is equal to the integer associated to $\hat\varphi(\hat{V}_j)$. Property (2) in this definition  ensures that the pair of annuli corresponding to some $\hat U_\sigma$ are mapped to a similar pair of annuli for $f'$. The need for the requirement that  $\varphi^*(U^*_\sigma)=U^*_{\sigma'}$ will be clear when considering the diffeomorphisms corresponding to Figures~\ref{11} and \ref{12}. For each of those diffeomorphisms,  the corresponding set $\hat{\mathcal U}_f$ consists of one annulus (which wraps 5 times around the torus). The difference between these diffeomorphisms can only be seen by considering the sets $U^*_\sigma$ in ${\mathcal V}_f^*$.
\end{rema}

\subsection{Schemes are complete invariants: statement of Theorem~\ref{t.cla}}\label{subsec:statementThm2}

\begin{theo} \label{t.cla} Two diffeomorphisms $f,f'\in MS(M^2)$ are topologically conjugate iff their schemes $S_f,S_{f'}$ are equivalent.
\end{theo}

As small perturbations of a diffeomorphism $f\in MS(M^2)$ do not change its periodic dates and the topological structure of the maximal $u$-compatible system of neighborhood, we obtain

\begin{coro}\label{cor:strstab}
 Each diffeomorphisms $f\in MS(M^2)$ is structurally stable (in the $C^1$ topology).
\end{coro}

The previous theorem shows that one can also consider two MS diffeomorphisms which are {\lq}far away from each{\rq}, provided their schemes are equivalent.

In the next section we will also introduce the notion of a {\em decomposed scheme}, which will make it clear
why a scheme is determined by a finite amount of data.  It will also be shown that 
an {\lq}abstract{\rq} decomposed scheme is realizable by a MS diffeomorphisms if and only 
if some simple conditions are satisfied.

\subsection{Examples of schemes associated to gradient MS-diffeomorphisms}

By the construction a scheme of a Morse-Smale diffeomorphism $f$ is completely described by a union of tori $\hat{\mathcal V}_f$, the epimorphism $\eta_{_f}$ and 
$\bigcup\limits_{i=1}^{beh(f)-1}\{\hat{U}_\sigma\}_{\sigma\in \Sigma_i}$ where the sets $\Sigma_1,\dots,\Sigma_{beh(f)-1}$ are finite. As mentioned above, the epimorphism $\eta_{_f}$ is uniquely determined by an equators on each component $\hat V_j$ of $\hat{\mathcal V}_f$ and an integer $m_j\ge 1$. In the figures in this paper,  the equators are depicted as the {\lq}outer boundaries{\rq} of the tori. 

\begin{figure}[h!]
\centerline{\includegraphics[width=14cm,height=9cm]{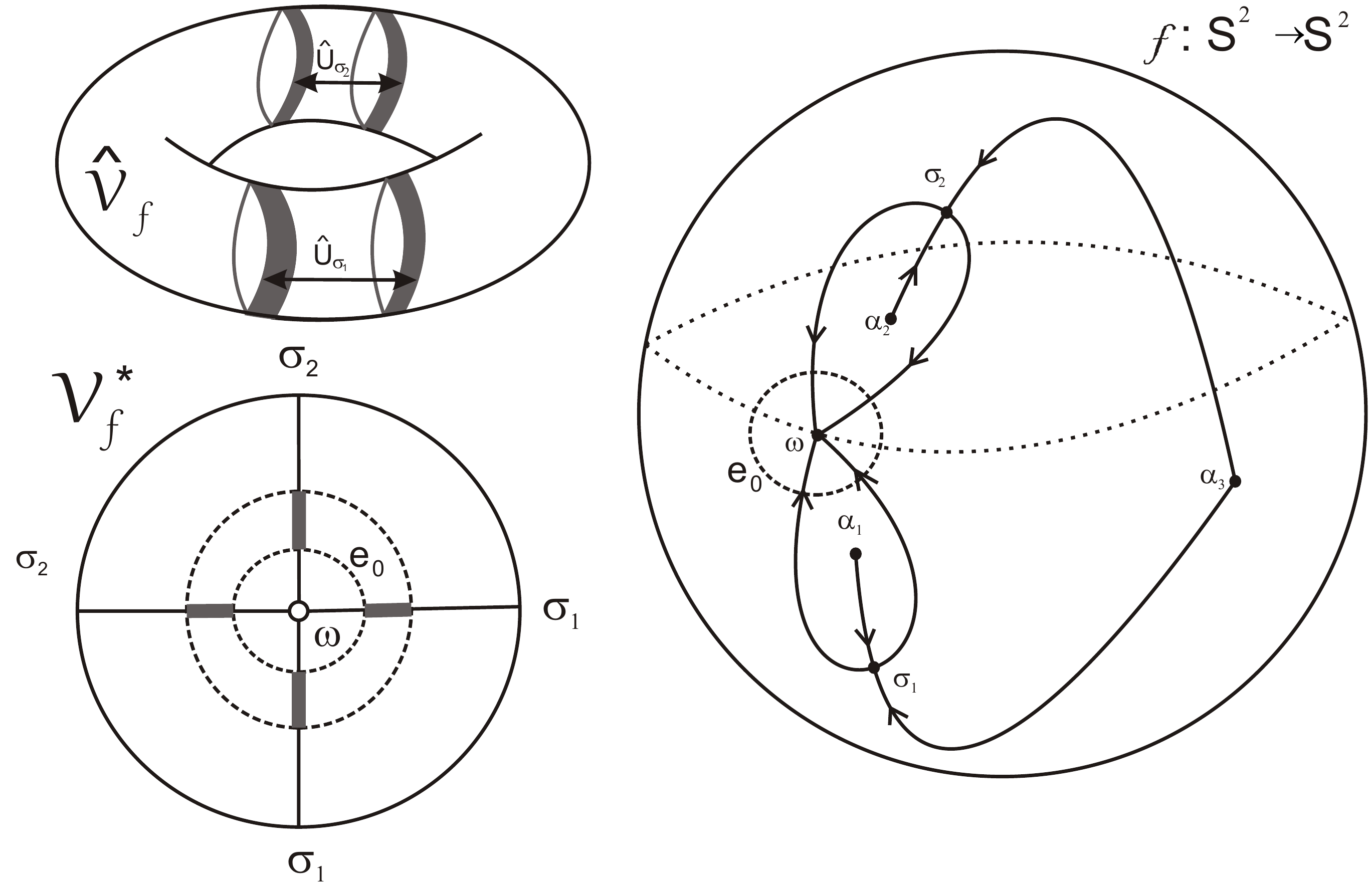}}\caption{\small The scheme $S_f$ of a diffeomorphism $f_S\in MS(\mathbb S^2)$. See Example~\ref{example:++schere} for a further discussion.} \label{++schere}
\end{figure}

\begin{example}[Figure~\ref{++schere}]\label{example:++schere} Consider the diffeomorphism $f\colon S^2\to S^2$  from Figure~\ref{++schere}. Here $\hat{\mathcal V}_f=\mathbb T^2$, $beh(f)=2$, 
$\#\Sigma_0=1$, $m=1$, $\#\Sigma_1=2$ and $\hat{U}_{\sigma_1},\hat{U}_{\sigma_2}$ consist of four annuli on the torus $\hat{\mathcal V}_f$.  The set $\hat V^*$ is represented as $(\mathbb R^2 \setminus O)$.
\end{example}

\begin{figure}[h!]
\centerline{\includegraphics[width=14cm,height=9cm]{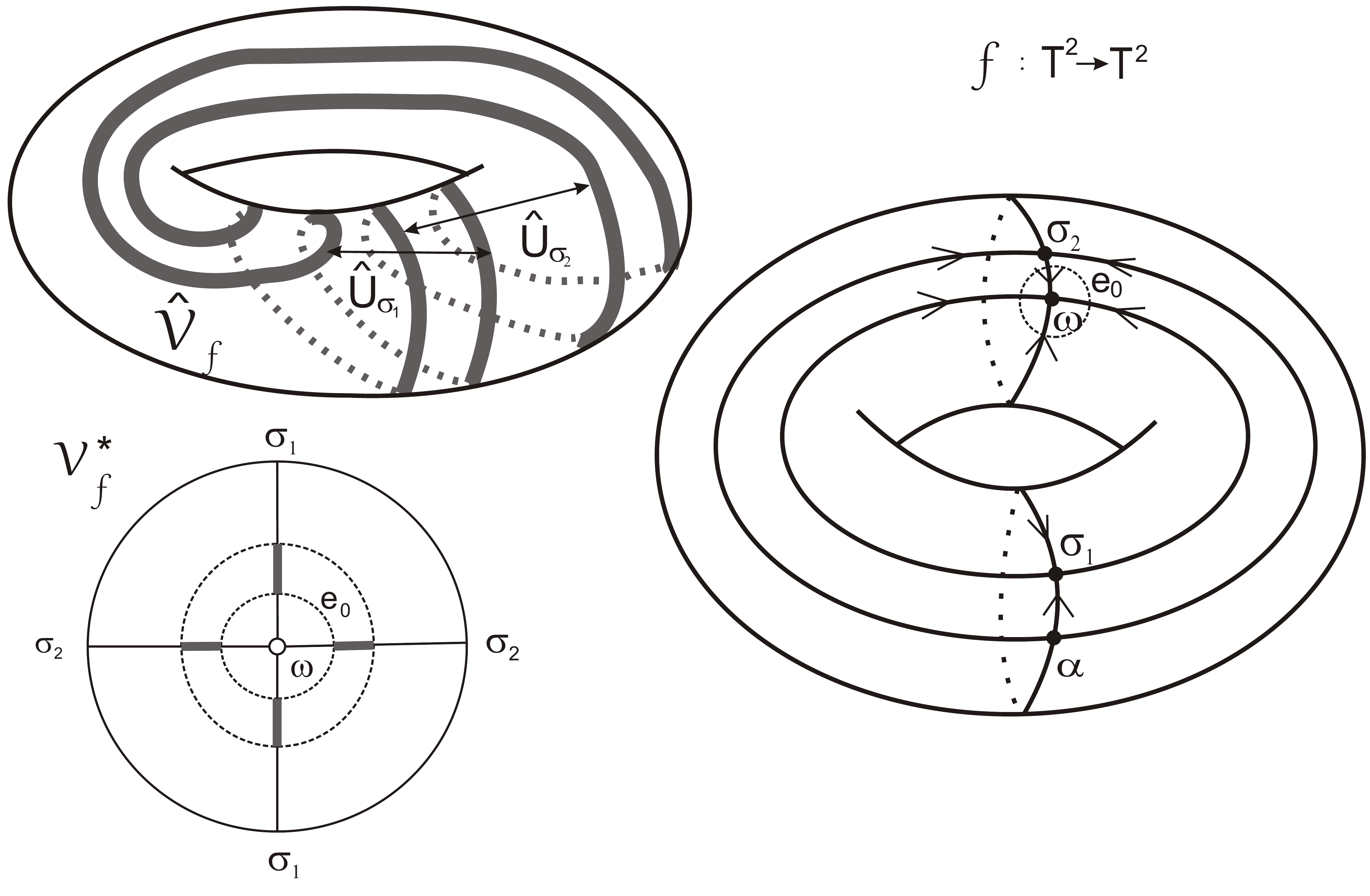}}\caption{\small The scheme $S_f$ associated to a diffeomorphism $f_S\in MS(\mathbb T^2)$ for which the multipliers at   the attractor $\omega$ are negative, and the two separatrices of both saddles are permuted by $f$. See Example~\ref{example:-s13} for a further discussion.} \label{s13}
\end{figure}

\begin{example}[Figure~\ref{s13}]\label{example:-s13} Consider the  diffeomorphism $f\colon S^2\to S^2$  from Figure~\ref{s13} where we take $f$ so that it permutes the  two components 
of $W^u(\sigma_1)\setminus \sigma_1$.  We obtain $\hat{\mathcal V}_f=\mathbb T^2$, $beh(f)=2$, $\#\Sigma_0=1$, $m=1$, $\#\Sigma_1=2$ and  $\hat{U}_{\sigma_i}$, $i=1,2$ is an annulus which winds around twice along the torus $\hat{\mathcal V}_f$.  To see this, notice that $f$ induces an action on ${\mathcal V}_f^*$ which corresponds to the composition of a radial contraction and a half revolution around $0$. Therefore the inner and outer circle drawn in  ${\mathcal V}_f^*$ are identified by a half revolution. 
\end{example}

\begin{figure}[h!]
\centerline{\includegraphics[width=12cm,height=11cm]{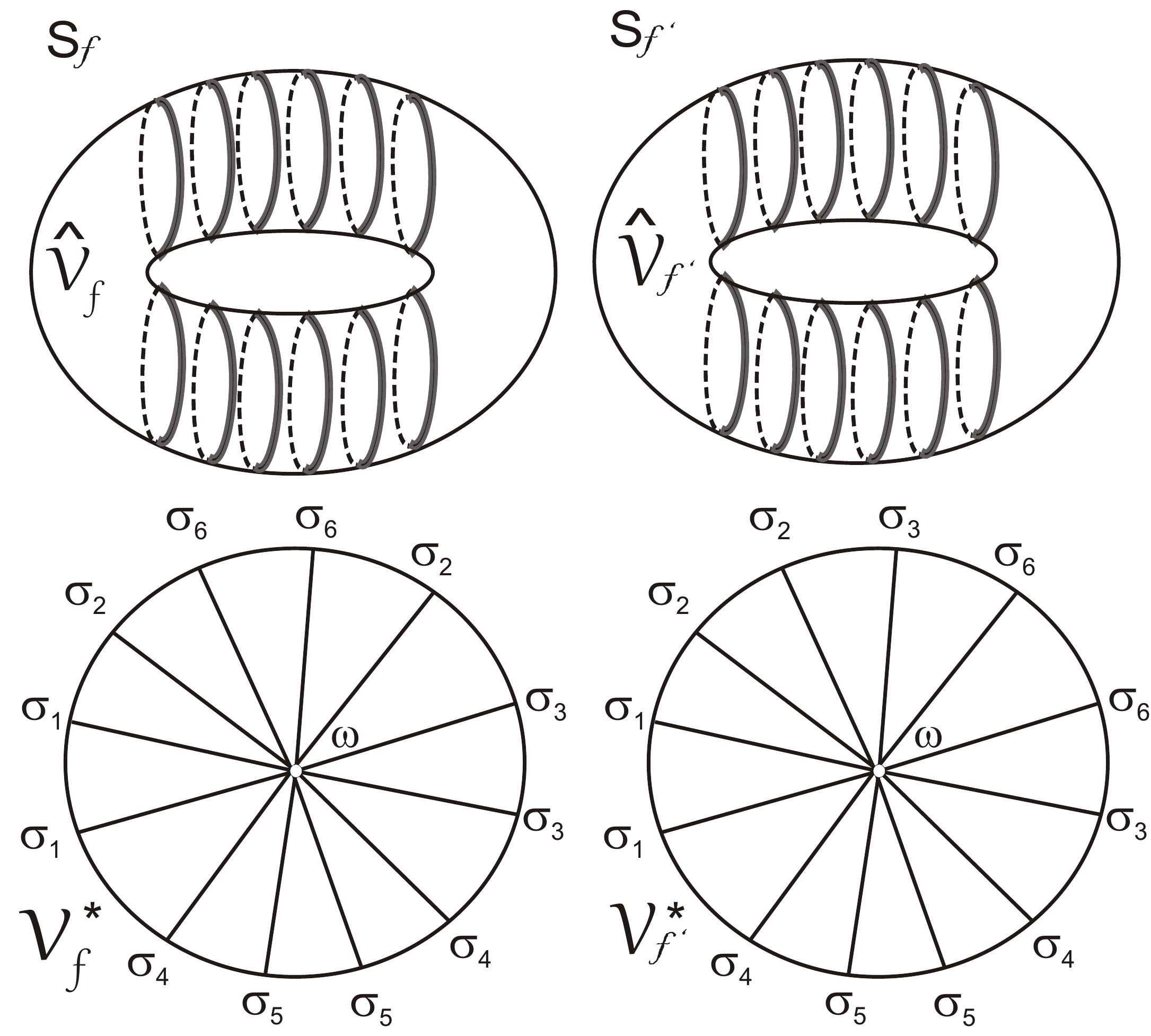}}\caption{The schemes associated to the diffeomorphisms $f,f'$ from Figure \ref{Expe}. For the schemes to be equivalent in the sense of Definition~\ref{def:equivschemes}, each pair of thickened curves corresponding to $\hat U_{\sigma_i}$ for $f$ has to correspond to a pair associated to $f'$.  Since there is no homeomorphism between the two tori preserving the pairs of annuli, by Theorem \ref{t.cla} below, 
the diffeomorphisms $f,f'$ are not topologically conjugated. See Example~\ref{example:ris:ex-12curves} for a further discussion.}
\label{12curves}
\end{figure}

\begin{example}[Figures~\ref{ris:ex},\ref{12curves}]\label{example:ris:ex-12curves}
Consider the diffeomorphisms $f,f'\colon S^2\to S^2$  from Figure~\ref{ris:ex}. The scheme associated to these diffeomorphisms are drawn in   Figure~\ref{12curves} where  $\hat{\mathcal V}_f=\mathbb T^2$, $beh(f)=2$, $\#\Sigma_0=1$, $m=1$, $\#\Sigma_1=5$ and $\hat{U}_{\sigma_i}$, $i=1,\dots,5$  consists of 2 annuli. Note that the difference between the two schemes is the way the pair of components corresponding to each of the five saddles are jointly  embedded in  $\hat{\mathcal V}_f$. 
\end{example}

{\begin{figure}[h!]
\centerline{\includegraphics[width=13cm,height=8cm]{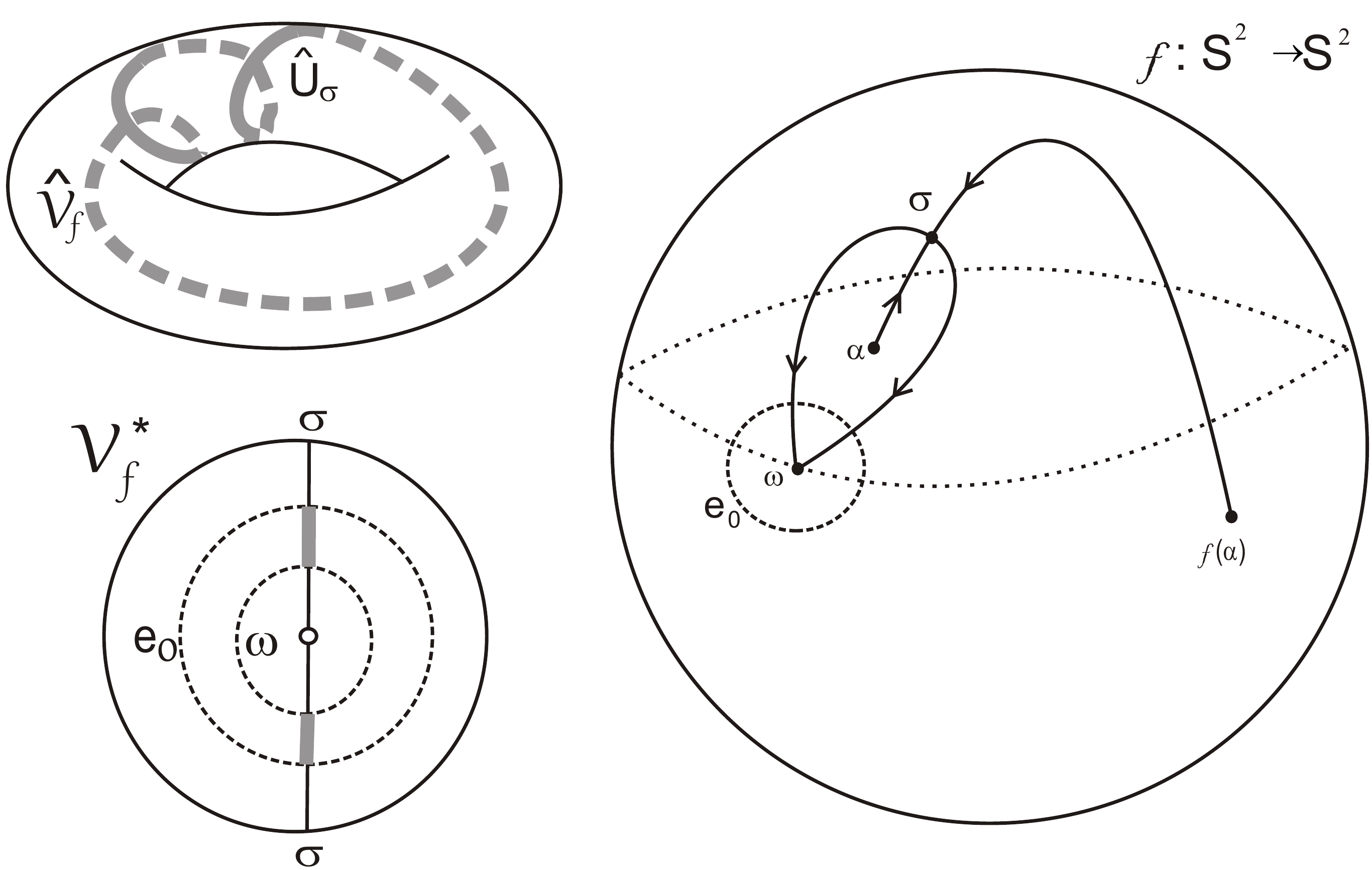}}\caption{\small The scheme $S_f$ associated to a diffeomorphism $f_S\in MS(\mathbb S^2)$ for which the multipliers at the attractor $\omega$ are negative, and the two separatrices of unique saddle are permuted by $f$. See Example~\ref{example:-s} for a further discussion.} \label{-sph}
\end{figure}}

\begin{figure}[h!]
\centerline{\includegraphics[width=16cm,height=6cm]{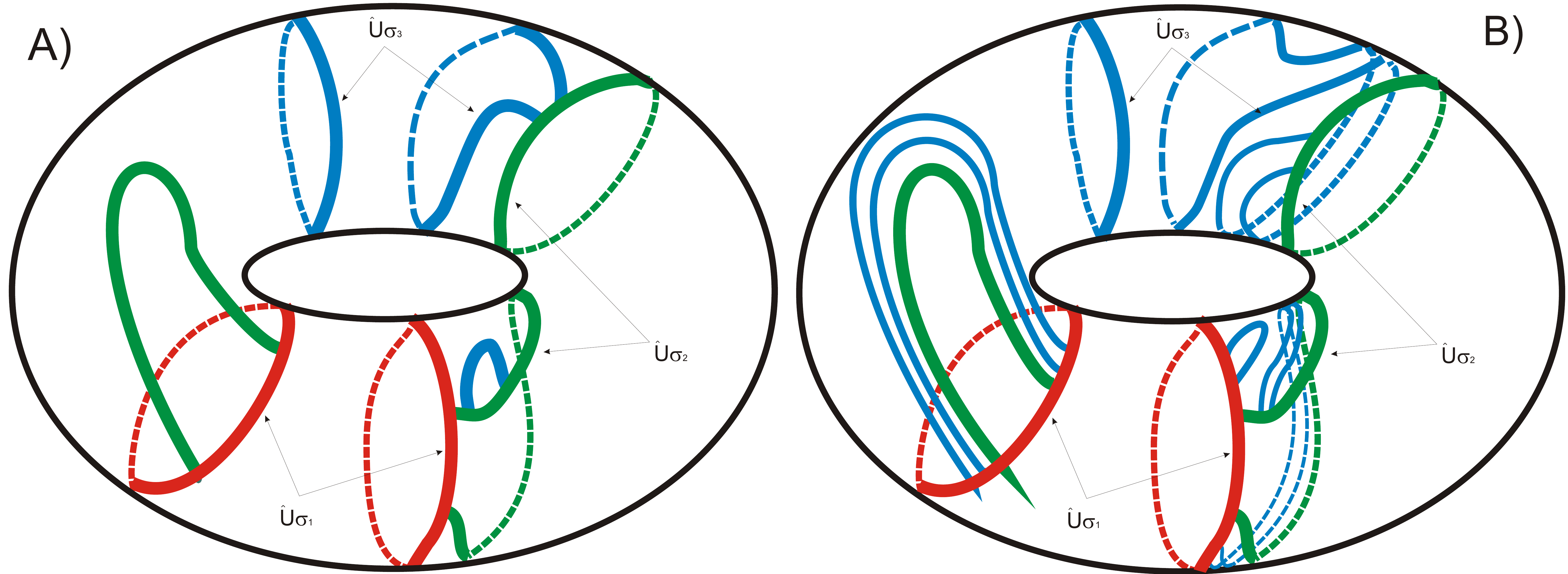}}\caption{\small {A) The scheme associated to the diffeomorphism $f$ from Figures \ref{faz+}A) and \ref{faz} with $beh(f)=4$. Here the torus is obtained by identifying, by $f$, the dashed curves  in Figure~\ref{faz}A) forming the boundary of the fundamental annulus of the sink $\omega$. 
B) The projection to the $\hat V_f$ of non-maximal $u$-compatible neighborhoods from Figures \ref{faz} B) for the same
diffeomorphism. The projections A) and B) are not homeomorphic, showing why it is important to consider 
maximal $u$-compatible neighborhoods.}}
\label{s1}
\end{figure}

{\begin{example}[Figure~\ref{-sph}]\label{example:-s} Consider the  diffeomorphism $f\colon S^2\to S^2$  from Figure~\ref{-sph} where we take $f$ so that it permutes the  two components 
of $W^u(\sigma)\setminus \sigma$.  We obtain $\hat{\mathcal V}_f=\mathbb T^2$, $beh(f)=2$, $\#\Sigma_0=1$, $m=1$, $\#\Sigma_1=1$ and  $\hat{U}_{\sigma}$ is an annulus which winds around twice along the torus $\hat{\mathcal V}_f$ as $\nu_\sigma=-$. To see this, notice that $f$ induces an action on ${\mathcal V}_f^*$ which corresponds to the composition of a radial contraction and a half revolution around $0$. Therefore the inner and outer circle drawn in  ${\mathcal V}_f^*$ are identified by a half revolution.
\end{example}}

\subsection{Examples of schemes associated to non-gradient MS-diffeomorphisms}

Up to now we only considered gradient MS-diffeomorphisms. Let us now consider two non-gradient systems:

\begin{example}[Figure \ref{s14}] The scheme of the MS diffeomorphism $f\colon S^2\to S^2$ from Figure \ref{s14}. Here  $\hat{\mathcal V}_f=\mathbb T^2$, $beh(f)=3$, $\#\Sigma_0=1$, $m=1$, $\#\Sigma_i=1$, $i=1,2$. The sets $\hat{U}_{\sigma_i}$, $i=2,3$ are no longer annuli in $\hat{\mathcal V}_f$. 
\end{example}

\begin{example}[Figures~\ref{faz+},\ref{faz},\ref{s1}]\label{example:s1} The scheme of the MS diffeomorphism $f\colon S^2\to S^2$ from Figures~\ref{faz+} and \ref{faz} is represented in Figure~\ref{s1}. Here  $\hat{\mathcal V}_f=\mathbb T^2$, $beh(f)=4$, $\#\Sigma_0=1$, $m=1$, $\#\Sigma_i=1$, $i=1,2,3$. 
The sets $\hat{U}_{\sigma_i}$, $i=2,3$ are no longer annuli in $\hat{\mathcal V}_f$. \end{example}

As the last example shows, the scheme for non-gradient like MS diffeomorphisms
can in general become quite complicated. Moreover it is much less easy to see how to reconstruct $f$ from the scheme. For this reason we will introduce the notion of a {\em \decomposedscheme} of a MS diffeomorphism in the next section.

\section{The decomposed scheme and realising abstract versions of such schemes by MS diffeomorphisms} \label{reaa}
In this section we will associate a  {\em decomposed scheme} to a MS surface diffeomorphism. 
For a gradient-like MS diffeomorphism this decomposed scheme 
coincides with the scheme defined in the previous section, but for 
non-gradient-like diffeomorphisms it consists of more, but simpler, pieces. 

It turns out that there one can give a few simple rules which determine whether or not an abstract version
of such a decomposed scheme  determines again a MS diffeomorphism.

\subsection{The decomposed scheme of a MS diffeomorphism and the statement of Theorem 2'}
 \label{deca}

For non gradient-like diffeomorphisms, it will be useful to introduce additional factor spaces. Namely, let $A_0=\Omega^0_f$ and for each $i\in\{1,\dots,beh(f)-1\}$ let us define $$A_i=A_0\cup\bigcup\limits_{j=0}^{i}W^u_j,~~\mathcal V_i=W^s_{\Omega_f\cap A_i}\setminus A_i.$$ Observe that $A_i$ is an attractor of $f$ and $f$ acts freely on $\mathcal V_i $. Set $\hat{\mathcal V}_i=\mathcal V_i/f$ and denote the natural projection by $$p_{_{i}}:\mathcal V_{i}\to\hat{\mathcal V}_{i}.$$ Notice that $\hat{\mathcal V}_0=\hat{\mathcal V}_f$ and $p_0=p_{_f}$. It will be proved in Section \ref{maxx} that each connected component $\hat{V}_{i,j}$ of $\hat{\mathcal V}_i$ is a torus,  $p_i$ is a covering map and that $p_i^{-1}(\hat{V}_{i,j})$ is again homeomorphic to ${V}_{i,j}^{*}=(\mathbb R^2\setminus 0)\times \Z_{m_{i,j}}$ where $m_{i,j}$ is the period of  the corresponding components of $A_i$. As before we define a morphism  $\eta_{i}:\pi_1(\hat{\mathcal V}_i)\to \mathbb Z$, an {equator $e_{i,j}$ on ${\hat V}_{i,j}$ and  $\mathcal V^*_i=\bigcup\limits_j {V}_{i,j}^{*}$.}

For each $i\in\{1,\dots,beh(f)-1\}$ let $\mathcal G_i=W^s_{\Sigma_i}\setminus\Sigma_i$ and $\hat{\mathcal G}_i=p_i(\mathcal G_i)$. The {\it decomposed scheme} associated to $f$ is 
$$S_i=(\hat{\mathcal V}_i,\eta_i,\hat{\mathcal G}_{i}),{i=1,\dots,beh(f)-1}.$$ 

\begin{defi}[Equivalence of decomposed schemes] \label{def:equivdeschemes}
Two decomposed schemes $S_i\mbox{ and }S'_{i}$ are {\em equivalent} if there exist orientation preserving homeomorphisms $\hat\varphi_i:\hat{\mathcal V}_i\to\hat{\mathcal V}'_{i}$, $i=1,\dots,beh(f)-1$
such that:\\
(1) $\eta'_{i}\hat\varphi_{i*}=\eta_{i}$;\\
(2) $\hat\varphi(\hat{\mathcal G}_i)=\hat{\mathcal G}'_{i}$, moreover for every point $\sigma\in\Sigma_{i}$ there is a point $\sigma'\in\Sigma'_{i}$ such that $\varphi_i^*(W^{s*}_\sigma)=W^{s*}_{\sigma'}$, where $\varphi_i^*:{\mathcal V}_i^*\to{\mathcal V}_{i}^{\prime *}$ is the lift of $\hat\varphi_i$.
\label{eqvde}
\end{defi}

Let $\mathcal{N}_f=\{{N}_\si,\si\in\Si\}$ be a maximal $u$-compatible system of neighborhoods for a diffeomorphism $f\in MS(M^2)$, $\mathcal{N}_i=\bigcup\limits_{\si\in\Si_i}{N}_\si$ for $i\in\{1,\dots,beh(f)-1\}$. 

Analogously to Theorem~\ref{t.cla} we have 

\bigskip

\noindent 
\begin{theo2}\label{thm2'} 
{\em Two diffeomorphisms $f,f'\in MS(M^2)$ are topologically conjugate iff their decomposed schemes $S_f,S_{f'}$ are equivalent.}
\end{theo2} 

\medskip

We also will define an abstract version of a decomposed scheme. Theorem \ref{t.realisation}
states that an {\em abstract} decomposed scheme is equivalent to a decomposed scheme 
of a  diffeomorphism $f\in MS(M^2)$ if and only if Definitions~\ref{def:allowed} and \ref{opr}
are satisfied.

\subsection{Cut \& paste operations and examples of decomposed schemes}\label{cp}
\begin{figure}[h!]\centerline{\includegraphics[width=17cm,height=9cm]{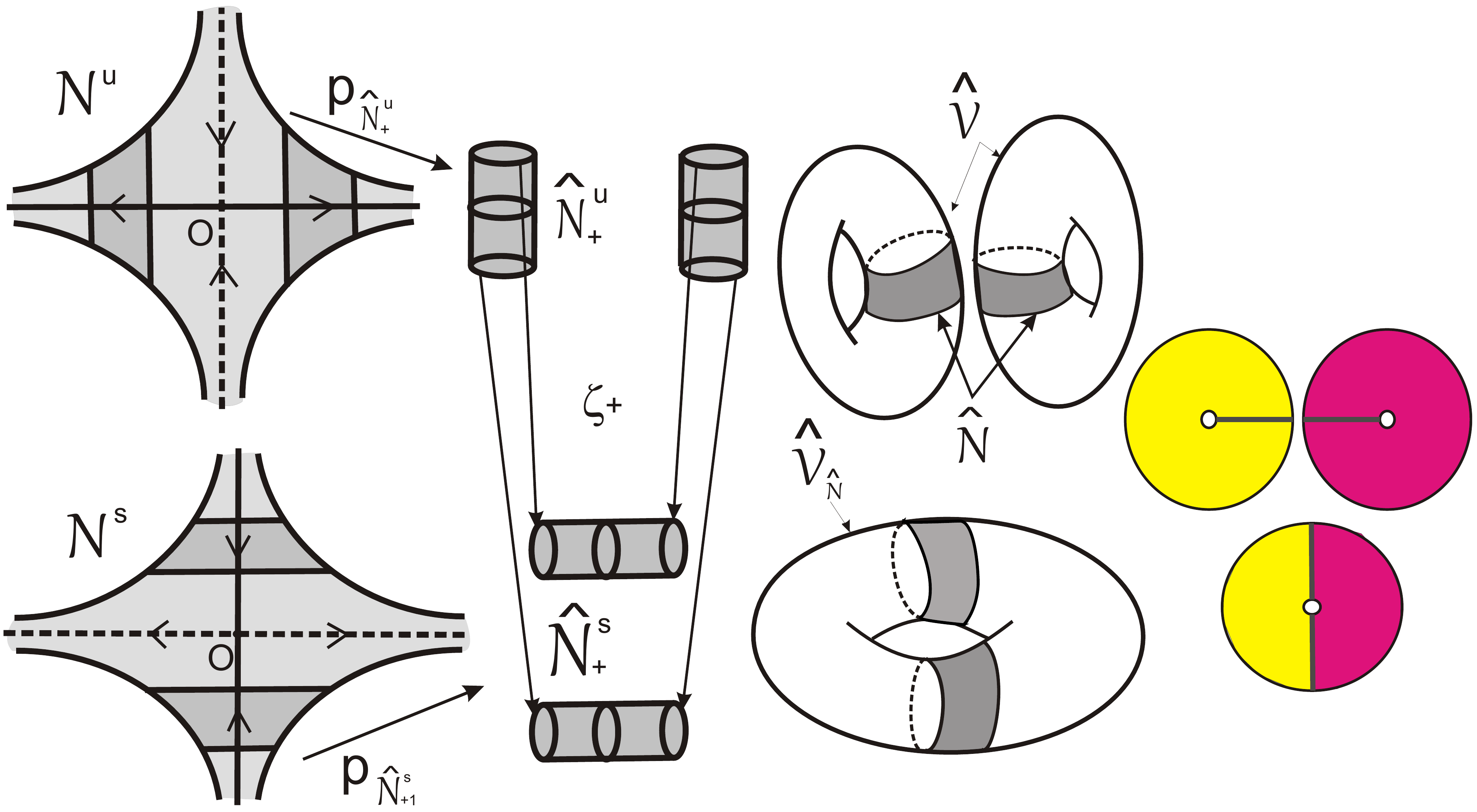}}\caption{\small Cut \& paste operations}\label{perec}
\end{figure}
\paragraph{Cut \& paste operations.}
For every diffeomorphism $f\in MS(M^2)$ the set $\hat{\mathcal V}_{i-1}$ can be obtained from
$\hat{\mathcal V}_{i}$  by a cut \& paste operation. Indeed, 
notice that the manifolds $\hat{\mathcal V}_{i}\setminus p_i(W^s_{i})$ and $\hat{\mathcal V}_{i-1}\setminus p_{i-1}(W^u_{i})$ are homeomorphic by the homeomorphism $\phi_i=p_{_{i}}p^{-1}_{_{i-1}}$. Also the manifolds $\hat{\mathcal V}_{i}\setminus p_i({\mathcal N}_{i})$ and $\hat{\mathcal V}_{i-1}\setminus p_{i-1}({\mathcal N}_{i})$ are homeomorphic by $\phi_i$. Each connected component $\hat {\mathcal N}_{\si}$ of the set $\hat {\mathcal N}_{i}$ is homeomorphic to $\hat N^u_{\nu_\si}$ by means $\mu_{_{\hat N_\si}}=p_{_{\hat N^u_{\nu_\si}}}\mu_{\si}p^{-1}_{_{i-1}}$ and each connected component $\hat {\mathcal N}^s_{\si}$ of the set $\hat {\mathcal N}^s_{i}$ is homeomorphic to $\hat N^s_{\nu_\si}$ by means $\mu_{_{\hat N^s_\si}}=p_{_{\hat {N}^s_{\nu_\si}}}\mu_{\si}p^{-1}_{_{i}}$. Thus $\hat{\mathcal V}_{i-1}$ formally can be obtained from $\hat{\mathcal V}_{i}$ by a regluing of annuli $\hat{\mathcal V}_{i}\setminus\hat{\mathcal G}_i$ or regluing of sectors ${\mathcal V}^*_{i}\setminus{\mathcal G}
_i$ (see Figure \ref{perec}).

For example for a gradient-like case we have only two spaces $\mathcal V_0$ which is a union of  basins of the sink points and $\mathcal V_1$ which is a union of  basins of the source points. The transition from $\mathcal V_1$ to $\mathcal V_0$ consists of a consequent execution of cut \& paste operations along every saddle orbits.  How we will see below the transition along one saddle orbit  is of three types: 

\begin{figure}[h!]
\centerline{\includegraphics[width=13cm,height=8cm]{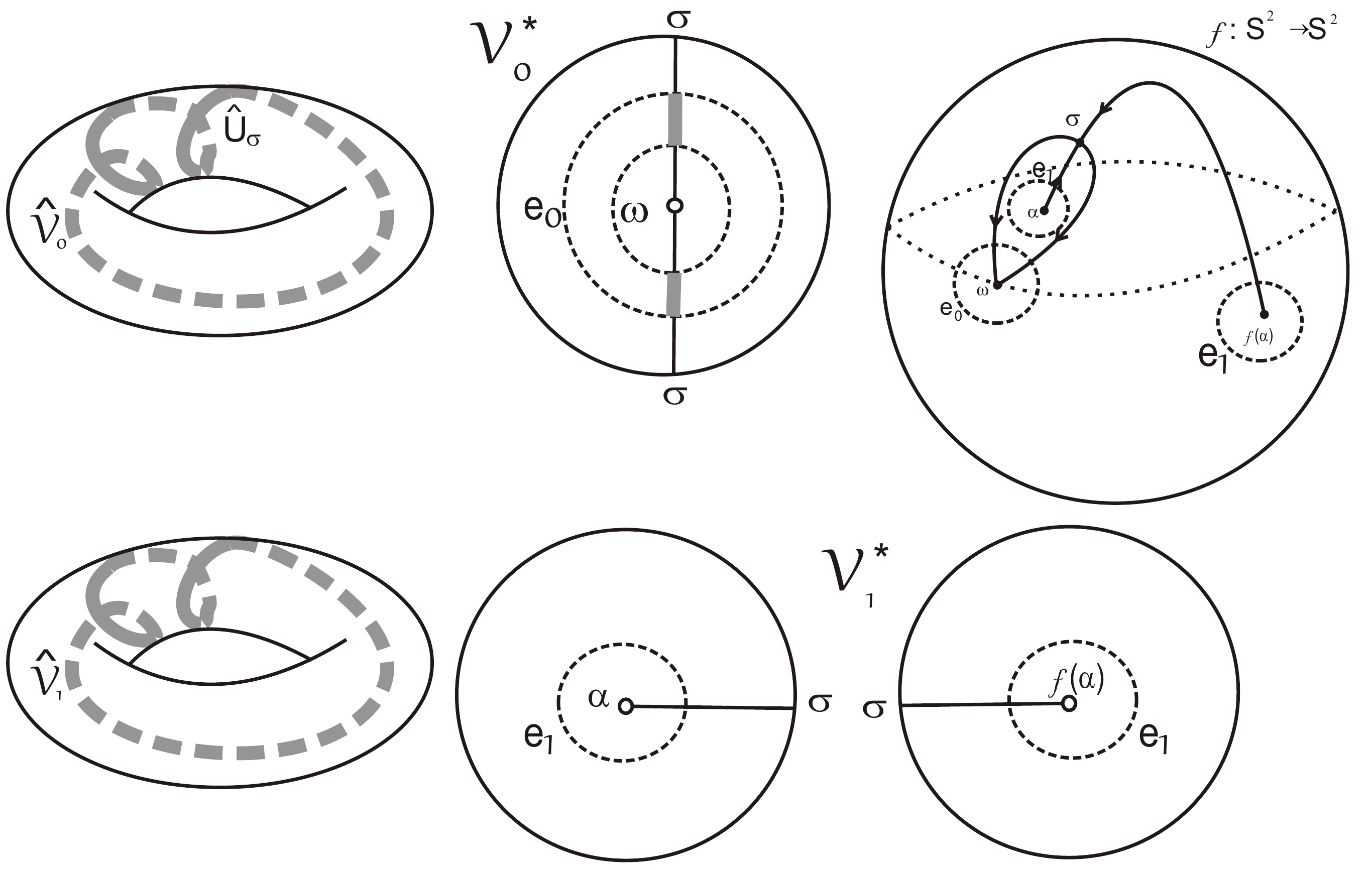}}\caption{\small The transition of the type 1 from $\mathcal V_1$ to $\mathcal V_0$ for the diffeomorphism $f_S\in MS(\mathbb S^2)$ from Figure \ref{-sph} for which the multipliers at the saddle point $\sigma$ are negative, and the two separatrices of unique saddle are permuted by $f$} \label{-sph+}
\end{figure}  

\begin{figure}[h!]
\centerline{\includegraphics[width=13cm,height=8cm]{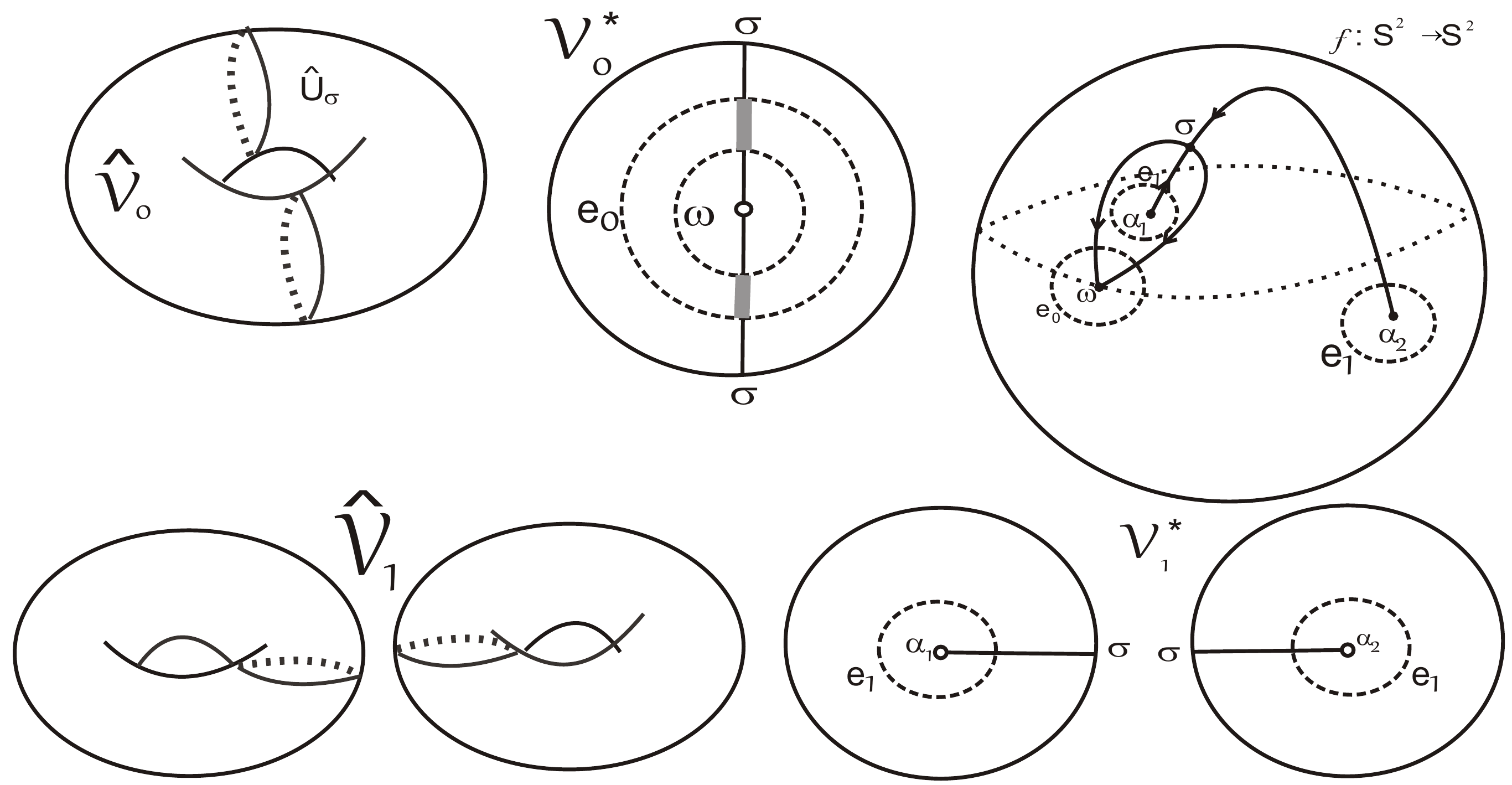}}\caption{\small The transition of the type 2 from $\mathcal V_1$ to $\mathcal V_0$ for the diffeomorphism $f_S\in MS(\mathbb S^2)$  for which the multipliers at the saddle point $\sigma$ are positive, and the two separatrices to a basin of the same sink point $\omega$} \label{-sph-2}
\end{figure}

\begin{figure}[h!]
\centerline{\includegraphics[width=13cm,height=7.5cm]{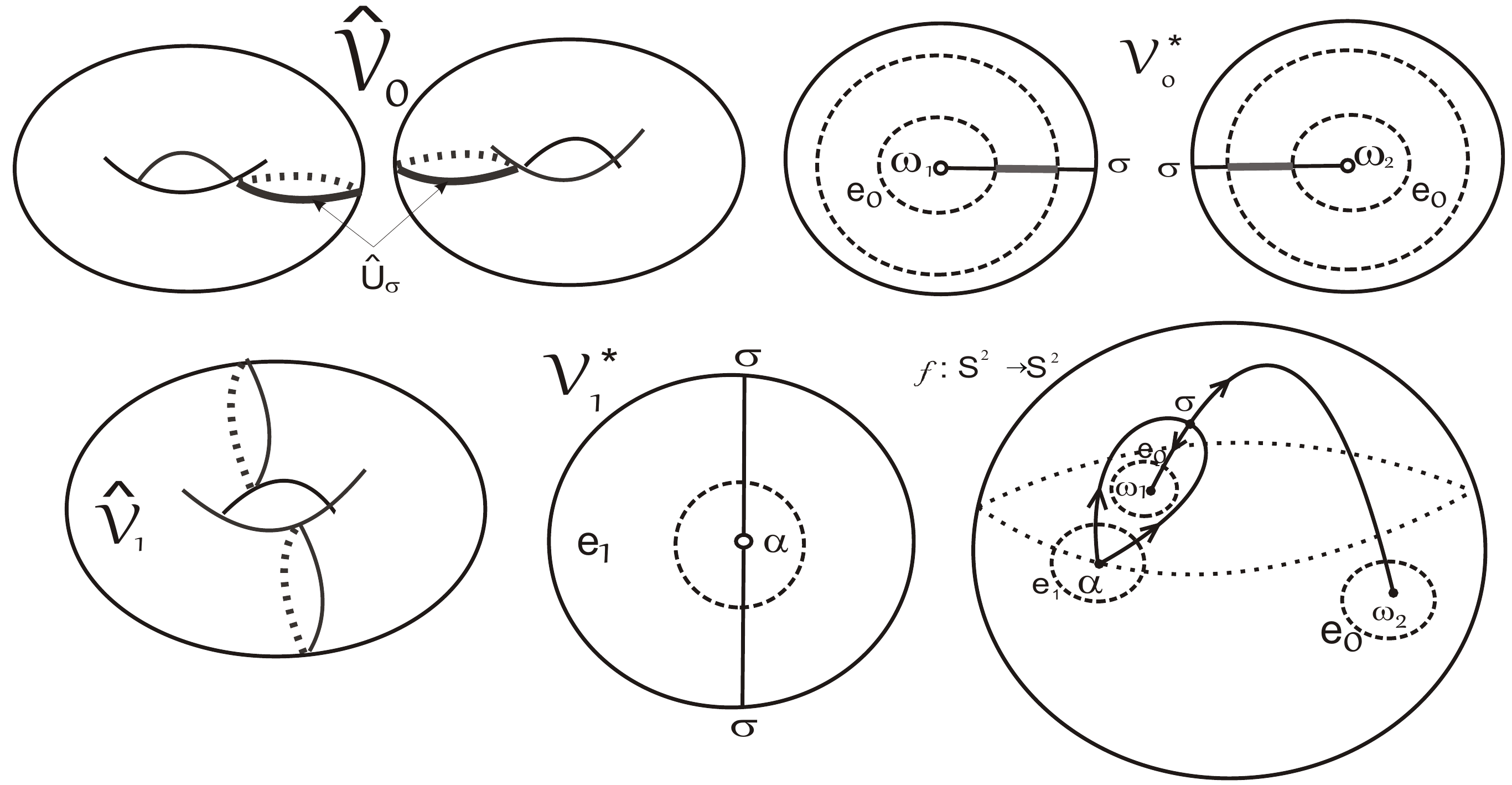}}\caption{\small The transition of the type 3 from $\mathcal V_1$ to $\mathcal V_0$ for the diffeomorphism $f_S\in MS(\mathbb S^2)$ for which the multipliers at the saddle point $\sigma$ are positive, and the two separatrices to basins of different sink points $\omega_1,\omega_2$} \label{-sph-3}
\end{figure}

\begin{enumerate} 
\item The stable separatrices of a saddle periodic orbit with {\em negative multipliers} enter the basin of a periodic attractor, see Figure \ref{-sph+}.
\item The stable separatrices of a saddle periodic orbit with {\em positive multipliers} enter  the basin of a periodic attractor, see  Figure \ref{-sph-2}.
\item The stable separatrices of a saddle periodic orbit with {\em positive  multipliers}  belong to the basins of different periodic attractors, see Figure \ref{-sph-3}.  
\end{enumerate}

\begin{figure}[h!]
\centerline{\includegraphics[width=16cm,height=7cm]{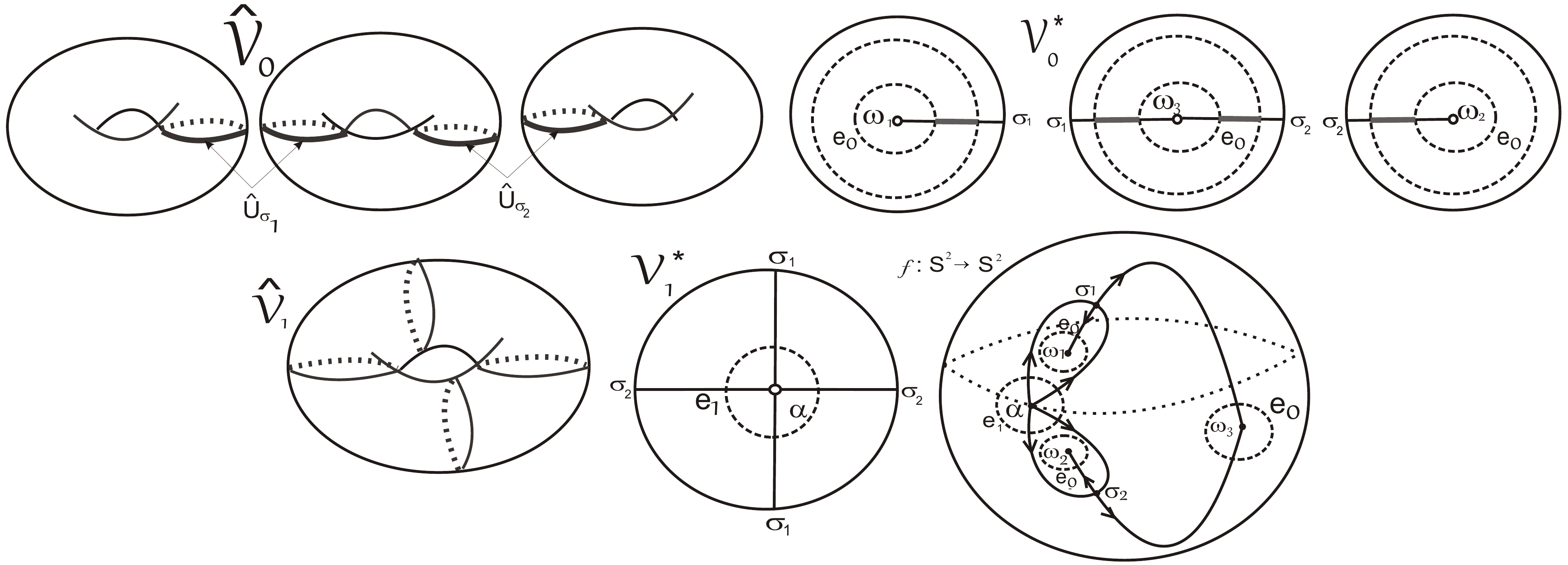}}\caption{\small The transition from $\mathcal V_1$ to $\mathcal V_0$ for the diffeomorphism $f_S\in MS(\mathbb S^2)$ for which there are several saddles with the separatrices entering the basins of sinks.}\label{two}
\end{figure}

If there are several saddles with the separatrices entering the basins of sinks then we have to cut and paste several things at once, see Figure \ref{two}. 

\begin{figure}[h!]
\centerline{\includegraphics[width=16cm,height=10cm]{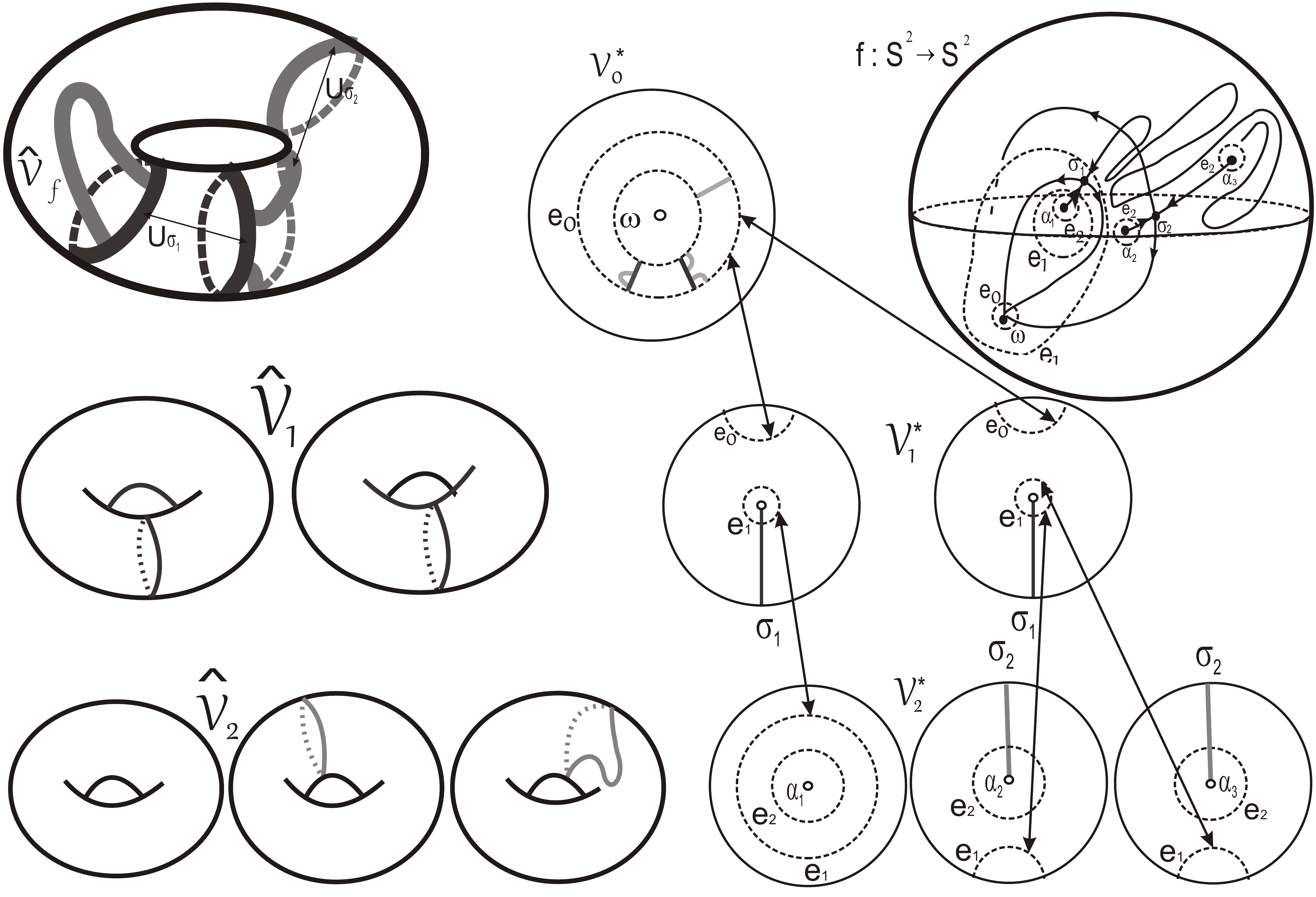}}\caption{\small On the top left the scheme $S_f$ is represented of the diffeomorphism $f_S\in MS(\mathbb S^2)$ on the top right. 
The corresponding decomposed scheme consists of two objects $\hat V_1,\hat V_2$ with circles on them.  
Here $\hat{\mathcal V}_1~(\hat{\mathcal V}_0)$ is obtained by cutting and gluing along the circles in the torus $\hat{\mathcal V}_2~(\hat{\mathcal V}_1)$. The sets $\mathcal V_0^*,\mathcal V_1^*,\mathcal V^*_2$  are explained Example~\ref{example:s14},  as is the construction how  these can be used to get back the original diffeomorphism $f$. Namelly we have to identify the parts of the boundary in $\mathcal V_1^*,\mathcal V^*_2$ as suggested by the labelling of saddle separatrices and after that take the connected sum of $\mathcal V_0^*,\mathcal V_1^*,\mathcal V^*_2$ as suggested by the labelling of equators.}\label{s14}
\end{figure}

For a non gradient-like case there are more than one successive transitions (see Figure \ref{s14}) but every transition along one saddle orbit is of one of these three types. 
 
\begin{example}[Figure~\ref{s14}]\label{example:s14} Consider the MS-diffeomorphism $f$ shown in the 
top-right of Figure~\ref{s14} with $beh(f)=3$. 
\begin{itemize}
\item The set $\mathcal V_0=W^s_\omega\setminus \omega$ is homeomorphic to $\mathbb R^2\setminus O$.  We draw two equators $e_0$ in ${\mathcal V}_0^*$ which bound an annulus. The torus $\hat{\mathcal V}_0$ is  obtained by identifying the boundary curves of this annulus. We also add two curves $e_1$ both on the left and right side of  ${\mathcal V}_0^*$, corresponding to curves $e_1$ in $\mathbb S^2$ (and corresponding to equators on the tori $\mathcal V_1$). 

\item The attractor $A_1=W^u_{\sigma_1}\cup\omega$ is a circle, and its basin $\mathcal V_1$ consists of two components, each topologically a punctured disc, and so each contains two equators $e_1$. 

\item The attractor $A_2=W^u_{\sigma_2}\cup W^u_{\sigma_1}\cup\omega$, whose basin consists of three components surrounding $\alpha_1,\alpha_2,\alpha_3$, each topologically a punctured disc, and so each contains three equators $e_2$. 

\item Now it is easy to reconstruct the original diffeomorphism from the decomposed scheme as follows: we identify the parts of the boundary in $\mathcal V_1^*,\mathcal V^*_2$ as suggested by the labelling of saddle separatrices and after that take the connected sum of $\mathcal V_0^*,\mathcal V_1^*,\mathcal V^*_2$ as suggested by the labelling of equators.
\end{itemize}
\end{example} 

\begin{figure}[h!]
\centerline{\includegraphics[width=16cm,height=12cm]{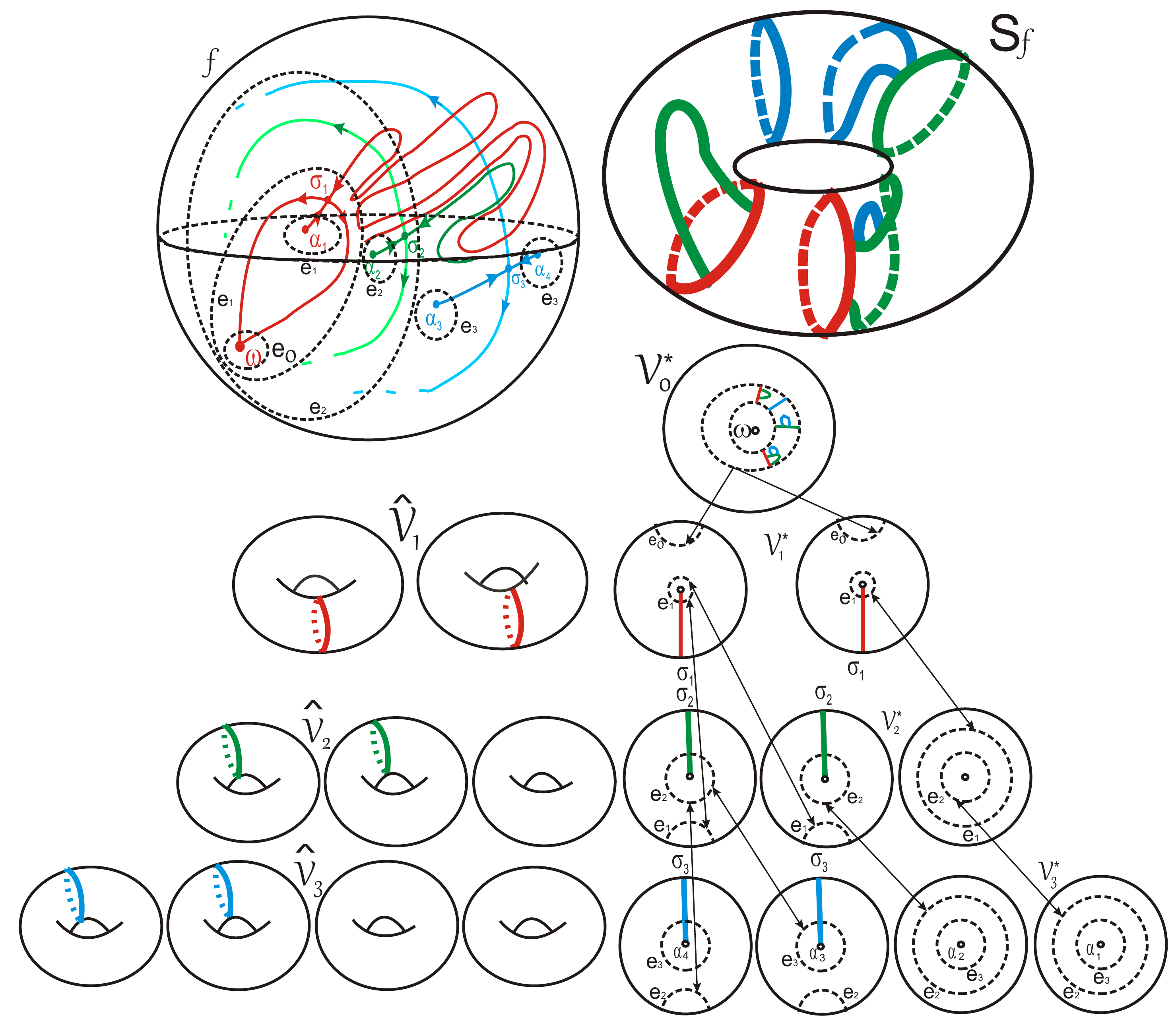}}\caption{\small On the top right the scheme $S_f$ is represented of the diffeomorphism $f_S\in MS(\mathbb S^2)$ on the top left. 
The corresponding decomposed scheme consists of two objects $\hat V_1,\hat V_2,\hat V_3$ with circles on them.  
Here $\hat{\mathcal V}_2~(\hat{\mathcal V}_1,\hat{\mathcal V}_0)$ is obtained by cutting and gluing along the circles in the torus $\hat{\mathcal V}_3~(\hat{\mathcal V}_2,\hat{\mathcal V}_1)$. The sets $\mathcal V_0^*,\mathcal V_1^*,\mathcal V^*_2,\mathcal V^*_3$  are explained Example~\ref{08},  as is the construction how  these can be used to get back the original diffeomorphism $f$. Namelly we have to identify the parts of the boundary in $\mathcal V_1^*,\mathcal V^*_2,\mathcal V^*_3$ as suggested by the labelling of saddle separatrices and after that take the connected sum of $\mathcal V_0^*,\mathcal V_1^*,\mathcal V^*_2,\mathcal V^*_3$ as suggested by the labelling of equators.}\label{beh4}
\end{figure}

\begin{example}[The decomposed scheme associated to the diffeomorphism from 
Figures \ref{faz+}A) and \ref{faz} with $beh(f)=4$.]\label{08} Consider the MS-diffeomorphism $f$ shown in the 
top-left of Figure~\ref{beh4} with $beh(f)=4$. 
\begin{itemize}
\item The set $\mathcal V_0=W^s_\omega\setminus \omega$ is homeomorphic to $\mathbb R^2\setminus O$.  We draw two equators $e_0$ in ${\mathcal V}_0^*$ which bound an annulus. The torus $\hat{\mathcal V}_0$ is  obtained by identifying the boundary curves of this annulus. We also add two curves $e_1$ both on the left and right side of  ${\mathcal V}_0^*$, corresponding to curves $e_1$ in $\mathbb S^2$ (and corresponding to equators on the tori $\mathcal V_1$). 

\item The attractor $A_1=W^u_{\sigma_1}\cup\omega$ is a circle, and its basin $\mathcal V_1$ consists of two components, each topologically a punctured disc, and so each contains two equators $e_1$. 

\item The attractor $A_2=W^u_{\sigma_2}\cup W^u_{\sigma_1}\cup\omega$, whose basin consists of three components surrounding $\alpha_1,\alpha_2,\alpha_3$, each topologically a punctured disc, and so each contains three equators $e_2$. 

\item The attractor $A_3=W^u_{\sigma_3}\cup W^u_{\sigma_2}\cup W^u_{\sigma_1}\cup\omega$, whose basin consists of four components surrounding $\alpha_1,\alpha_2,\alpha_3,\alpha_4$, each topologically a punctured disc, and so each contains three equators $e_3$. 

\item Now it is easy to reconstruct the original diffeomorphism from the decomposed scheme as follows: we identify the parts of the boundary in $\mathcal V_1^*,\mathcal V^*_2,\mathcal V^*_3$ as suggested by the labelling of saddle separatrices and after that take the connected sum of $\mathcal V_0^*,\mathcal V_1^*,\mathcal V^*_2,\mathcal V^*_3$ as suggested by the labelling of equators.
\end{itemize}
\end{example}

In the next subsection we will discuss this cut \& paste operation in the setting of some model objects, 
and show that  the decomposed scheme associated to a MS surface diffeomorphism satisfies the compatibility and realizability properties from Definitions~\ref{def:allowed} and \ref{opr} (see 
Lemma~\ref{lem:real-vs-abstract} below).

\subsection{Abstract decomposed schemes defined through model objects}  \label{mode} 
Let $m\ge 1$ be an integer and $V_m=\mathbb S^1\times\mathbb R^+\times\mathbb Z_{m}$. Thus $V_m$ is a model for the basin of a periodic attractor (equivalently the basin of dual repeller) of {\em period} $m$. Let $k\in\mathbb N$, an integer $n\geq 0$ so that $n=0$ if $k=1$ and otherwise $n\in\{1,\dots,k-1\}$ so that $n$ and $mk$ are coprime. Here $k$ models the number of saddle stable separatrices in each connected component of $V_m$ and $\frac{n}{k}$ represents their {\lq}rotation number{\rq}, i.e.,  how the diffeomorphism permutes these separatrices. As a local modal for the diffeomorphism
on the basin we take  the contraction $\phi_{m,k,n}:V_m\to V_m$ given by the formula: 
$$\phi_{m,k,n}(z,r,l)=(e^{\frac{2\pi n}{mk}i}z,\frac{r}{2^m},l+1(mod~m)).$$ Let $t\in[0,\frac1k)$, $j\in\{0,\dots,mk-1\}$, $\gamma^j=\phi^j_{m,k,n}(e^{i2\pi t}\times\mathbb R^+\times\{0\})$ and $$\gamma=\bigcup\limits_{j=0}^{mk-1}\gamma^j.$$ Thus $t$ defines the angle of a ray $\gamma^0$ in $V_m$ and $\gamma$ is $\phi_{m,k,n}$-invariant union of rays, containing $\gamma^0$, which we will refer to 
as a {\lq}frame{\rq}. So $\gamma$ models a saddle stable separatrix 
of period $mk$. Notice that $\hat V_m=V_m/\phi_{m,k,n}$ is a torus. Denote by  $p_{m,k,n}:V_m\to\hat V_m$ the natural projection. The set $\hat\gamma=p_{m,k,n}(\gamma)$ is a knot on $\hat V_m$. Let $e=\mathbb S^1\times\{1\}\times\{0\}$ and $\hat e=p_{m,k,n}(e)$. Denote by $\eta_{_{\hat V_m}}:\pi_1(\hat V_m)\to m\mathbb Z$ an epimorphism given by conditions: $\eta_{_{\hat V_m}}([\hat e])=0$ and $\eta_{_{\hat V_m}}([\hat \gamma])=km$.   
\begin{figure}[h!] 
\centerline{\includegraphics[width=17cm,height=6.5cm]{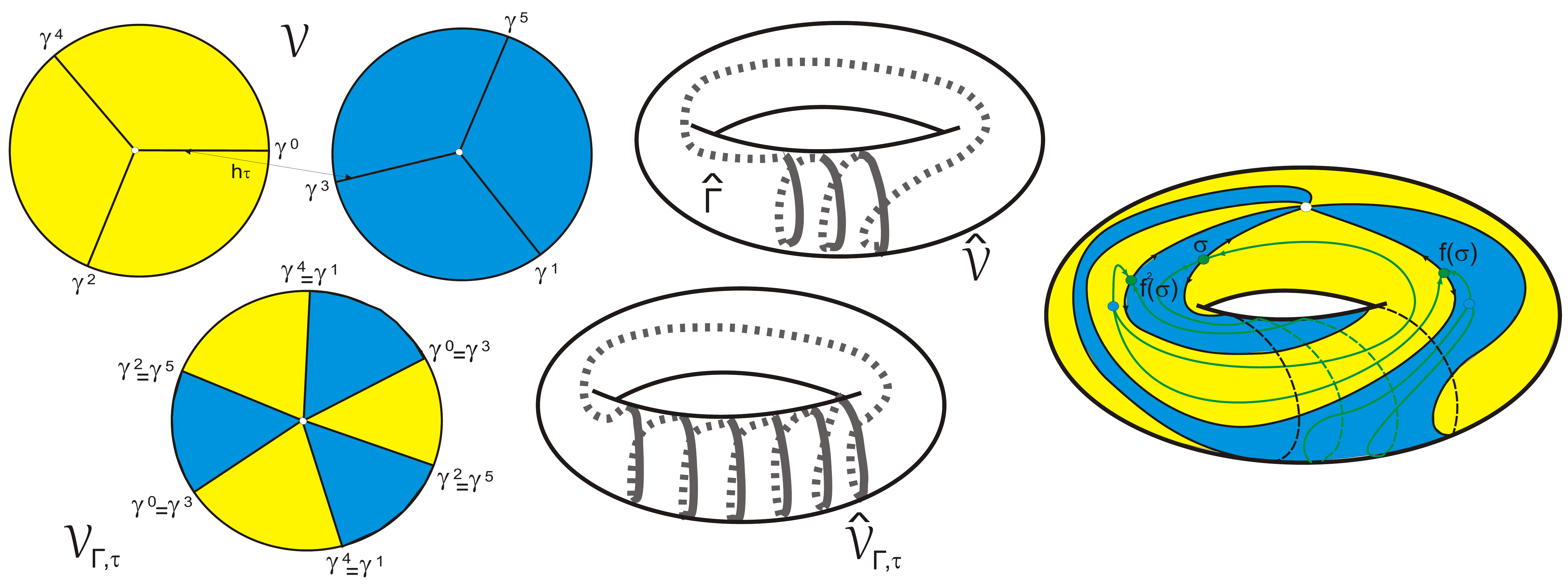}}\caption{\small The {\em first} type of ${\mathcal V}$, $\Gamma$, $h_\tau$ for $m=2,k=3,n=1,\tau=3$. Here $\hat{\mathcal V}_{\Gamma,\tau}$ is a unique torus with $m_{\Gamma,\tau}=1$ and a knot with 6 rotations as a projection of frame. On the left we can see realization of this part of dynamics on torus. Here blue points are a new attractor, which consists of two components of the period two and blue point with the union of green arcs is an old attractor of the period one.} \label{3NP}
\end{figure}

\begin{figure}[h!] 
\centerline{\includegraphics[width=16cm,height=10cm]{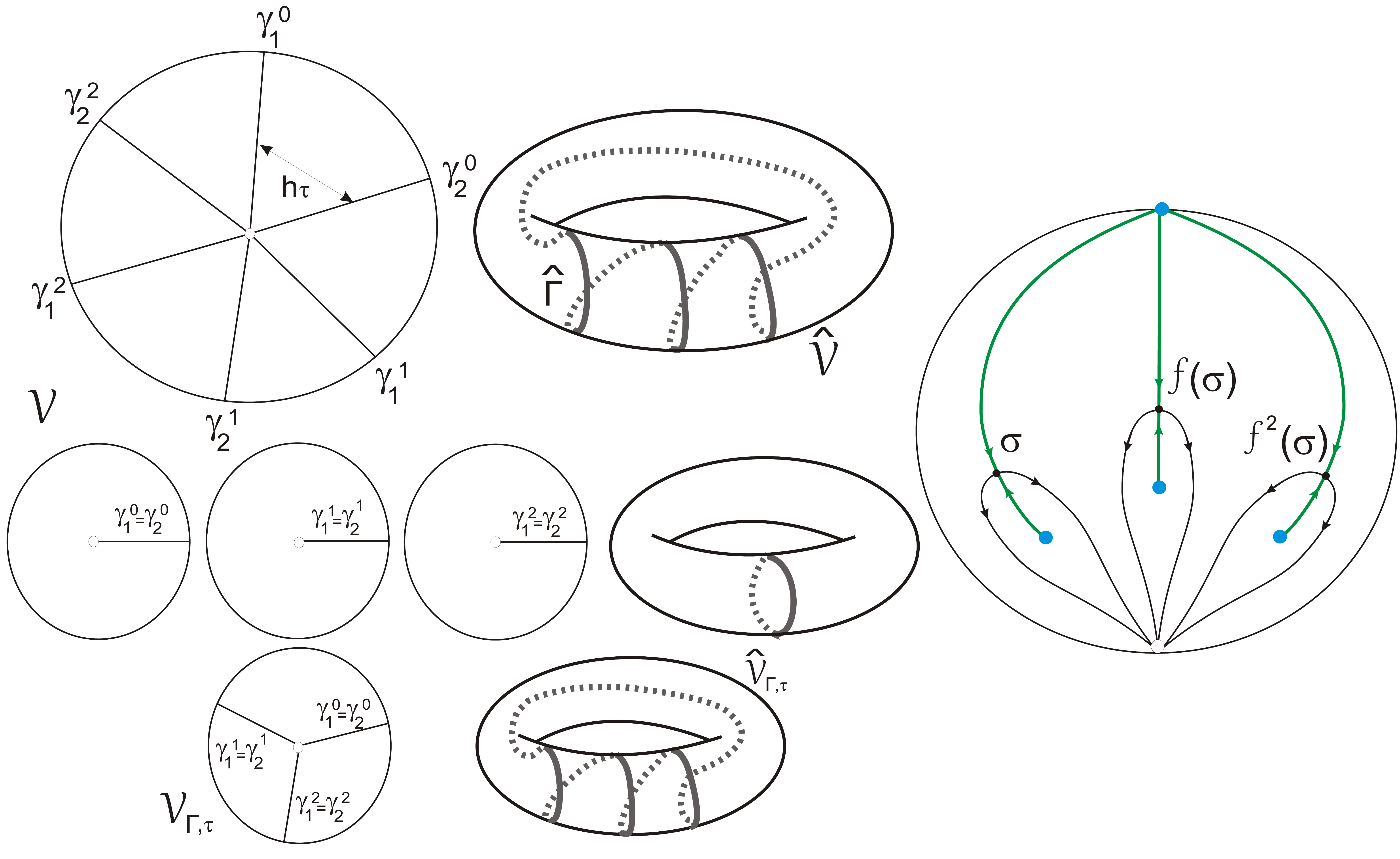}}\caption{\small The {\em second} type of ${\mathcal V}$, $\Gamma$, $h_\tau$ for $m=1,k=3,n=1,\tau=0$. Here $\hat{\mathcal V}_{\Gamma,\tau}$ consists of two tori with $m_{\Gamma,\tau,1}=3,\,m_{\Gamma,\tau,2}=1$ and two knots with 1 and 3 rotations as a projection of frame. On the left we can see realization of this part of dynamics on 2-sphere. Here blue points are a new attractor, which consists of one fixed component and other component of the period three and blue point with the union of three green arcs is an old attractor of the period one.} \label{2NP}
\end{figure}

\begin{figure}[h!] 
\centerline{\includegraphics[width=14cm,height=14cm]{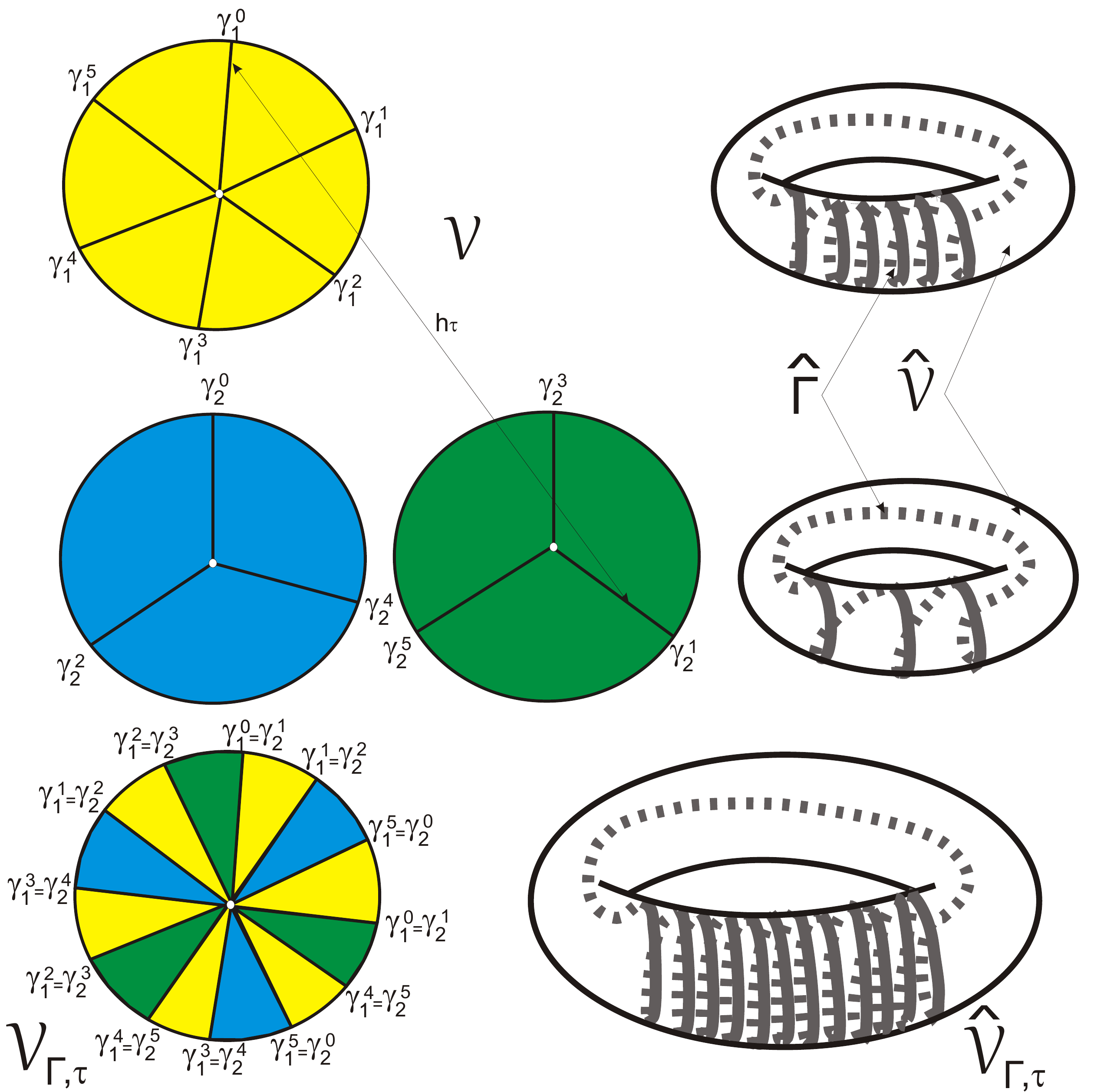}}\caption{\small The {\em third} type of ${\mathcal V}$, $\Gamma$, $h_\tau$ for $m_1=1,k_1=6,m_2=2,k_2=3,n_1=n_2=1,\tau=5$. Here $\hat{\mathcal V}_{\Gamma,\tau}$ is a unique torus with $m_{\Gamma,\tau}=1$ and a knot with 12 rotations as a projection of frame. The dynamics can represented on the surfase of the genus 2.}\label{1NP}
\end{figure}

\paragraph{Three types of collections of model objects.}
Let ${\mathcal V},\Gamma,\tau, h_\tau$ be a collection of one of the following three types:
\begin{enumerate}
\item ${\mathcal V}$ is $V_m$ for some $m$, $\Gamma$ is $\phi_{m,k,n}$-invariant frame $\gamma$ for some $k,n$ with
$km$ is  even,  $\tau=\frac{mk}{2}$, $h_\tau$ identifies the rays $\gamma^j$ and $\gamma^{\tau+j}$ for every $j\in\{0,\dots,\frac{km}{2}-1\}$, so that $h_\tau$ is equal to identity with respect to the coordinate in $\mathbb R^+$, i.e., $h(z, r,l)=(?,r,?)$ for 
all $(z,r,l)\in V_m$.  This case models the situation where the stable separatrices of a saddle periodic orbit with {\em negative multipliers} enter the basin of a periodic attractor of period $m$ so that $\gamma^j$ and $\gamma^{\tau+j}$ correspond to  the stable separatrices of one saddle point, see Figure \ref{3NP}.
\item ${\mathcal V}$ is $V_m$ for some $m$, $\Gamma$ is a {\em pair} of  $\phi_{m,k,n}$-invariant pairwise disjoint frames $\gamma_1,\gamma_2$ for some $k,n$, $\tau\in\{0,\dots,km-1\}$, $h_\tau$ identifies the rays $\gamma^j_1$ and $\gamma^{\tau+j}_2$ for every $j\in\{0,\dots,km-1\}$, where $h_\tau$  identity with respect to the coordinate in $\mathbb R^+$. This case models stable separatrices of a saddle periodic orbit with {\em positive multipliers} enter  the basin of a periodic attractor of period $m$ so that $\gamma^j_1$ and $\gamma^{\tau+j}_2$ correspond to the separatrices of one saddle point, see  Figure \ref{2NP}.
\item ${\mathcal V}$ is a disjoint union of $V_{m_1}$ and $V_{m_2}$ for some $m_1,m_2$, $\Gamma$ is a pair of $\phi_{m_1,k_1,n_1}$- and $\phi_{m_2,k_2,n_2}$-invariant pairwise disjoint frames $\gamma_1\subset V_{m_1},\gamma_2\subset V_{m_2}$ for some $k_1,n_1,k_2,n_2$ with $m_1k_1=m_2k_2$, $\tau\in\{0,\dots,k_1m_1-1\}$, $h_\tau$ identifies the rays $\gamma_1^j$ and $\gamma_2^{\tau+j}$ for every $j\in\{0,\dots,k_1m_1-1\}$ identity with respect to the coordinate in $\mathbb R^+$. This case models stable separatrices of a saddle periodic orbit with {\em positive  multipliers} so that $\gamma^j_1$ and $\gamma^{\tau+j}_2$ correspond to two separatrices of one saddle belonging to the basins of different periodic attractors, see Figure \ref{1NP}. 
\end{enumerate}

For every type denote by $\hat{\mathcal V}$ the corresponding space orbit, by $p_{_{\hat{\mathcal V}}}:\mathcal V\to\hat{\mathcal V}$ the natural projection, by $\hat\Gamma$ the projection of $\Gamma$ and by $\eta_{_{\hat{\mathcal V}}}$ a map composed by the corresponding epimorphisms. 

\paragraph{Cut \& paste operations to obtain ${\mathcal V}_{\Gamma,\tau}$.}
Denote by ${\mathcal V}_{\Gamma,\tau}$ the space obtained by regluing the closures of sectors in $\mathcal V\setminus\Gamma$ along the boundary with respect to $h_{\tau}$ (see Figures \ref{3NP}-\ref{1NP} in the bottom). It follows from the definition of $h_{\tau}$ that every connected component of the set ${\mathcal V}_{\Gamma,\tau}$ is homeomorphic to $\mathbb S^1\times\mathbb R^+$. As the homeomorphism $h_\tau$ commutes with $\phi_{m,k,n}$ (conjugates $\phi_{m_1,k_1,n_1}$ with $\phi_{m_2,k_2,n_2}$)  there is a diffeomorphism $\phi_{\Gamma,\tau}:{\mathcal V}_{\Gamma,\tau}\to{\mathcal V}_{\Gamma,\tau}$ which permutes the connected components.

Let $\hat{\mathcal V}_{\hat\Gamma,\tau}={\mathcal V}_{\Gamma,\tau}/\phi_{\Gamma,\tau}$ 
and $p_{_{\hat{\mathcal V}_{\Gamma,\tau}}}:{\mathcal V}_{\Gamma,\tau}\to\hat{\mathcal V}_{\hat\Gamma,\tau}$ be the natural projection. By  construction $\hat{\mathcal V}_{\Gamma,\tau}$ is obtained by regluing the closures of annuli in $\hat{\mathcal V}\setminus\hat\Gamma$ along the boundary with respect to the projection of $h_{\tau}$. The set $\hat{\mathcal V}_{\Gamma,\tau}$ consists of one torus (resp. is a union of two tori) for the first and third types (resp. the second type). It means that $\phi_{\Gamma,\tau}$ forms the unique orbit (resp. two orbits) from the connected components of the set ${\mathcal V}_{\Gamma,\tau}$. Denote by $m_{\Gamma,\tau}$ its period (resp. $m_{\Gamma,\tau,1},m_{\Gamma,\tau,2}$ their periods) and by $\eta_{_{\hat{\mathcal V}_{\hat\Gamma,\tau}}}$ the corresponding epimorphisms. Thus $\hat{\mathcal V}_{\Gamma,\tau}$ models the basin of a new attractor which is the initial attractor with the stable separatrices of the saddle orbit and which has one periodic component of the period $m_{\Gamma,\tau}$ (two  periodic components of the periods $m_{\Gamma,\tau,1}$, $m_{\Gamma,\tau,2}$). 

We will say that $\hat{\mathcal V}_{\hat\Gamma,\tau}$ is {\it a result of  regluing $\hat{\mathcal V}$ along $\hat\Gamma$ with the parameter $\tau$}. Thus $\hat{\mathcal V}_{\hat\Gamma,\tau}$ is a torus for the first and third types and consists of two tori for the second type.

\paragraph{Cut \& paste operations to obtain ${\mathcal V}_{\hat {\mathcal G},\mathcal T}$.}

Let $\hat{\mathcal V}$ be pairwise disjoint tori $\hat V_1,\dots,\hat V_l$ with a collection $\eta_{_{\hat{\mathcal V}}}$ of epimorphisms $\eta_1:\pi(\hat V_1)\to m_1\Z,\dots,\eta_l:\pi(\hat V_l)\to m_l\Z$, $\hat\Gamma_1,\dots,\hat\Gamma_l\subset\hat{\mathcal V}$ be pairwise disjoint sets of the described above types with parameters $\tau_1,\dots,\tau_l$ and $$\hat{\mathcal G}=\bigcup\limits_{i=1}^l\hat\Gamma_i,\quad \mathcal T=\{\tau_1,\dots,\tau_l\}.$$ We  say that the manifold $\hat{\mathcal V}_{\hat{\mathcal G},\mathcal T}$  {\it is a result of the regluing of $\hat{\mathcal V}$ along $\hat{\mathcal G}$ with the parameters $\mathcal T$} if we execute the regluing along  $\hat\Gamma_i$ with  parameter $\tau_i$ for every $i$. Denote by $p_{_{\hat{\mathcal V}_{\hat{\mathcal G},\mathcal T}}}:\hat{\mathcal V}\to\hat{\mathcal V}_{\hat{\mathcal G},\mathcal T}$ the natural projection. Then $\hat{\mathcal V}_{\hat{\mathcal G},\mathcal T}$  again consists of a finite number tori on which the regluing induces a collection of epimorphisms $\eta_{_{\hat{\mathcal V}_{\hat{\mathcal G}}}}$. Notice that the result does not depend on the order in which one does the gluing. 

\begin{figure}[h!] 
\centerline{\includegraphics[width=12cm,height=6.5cm]{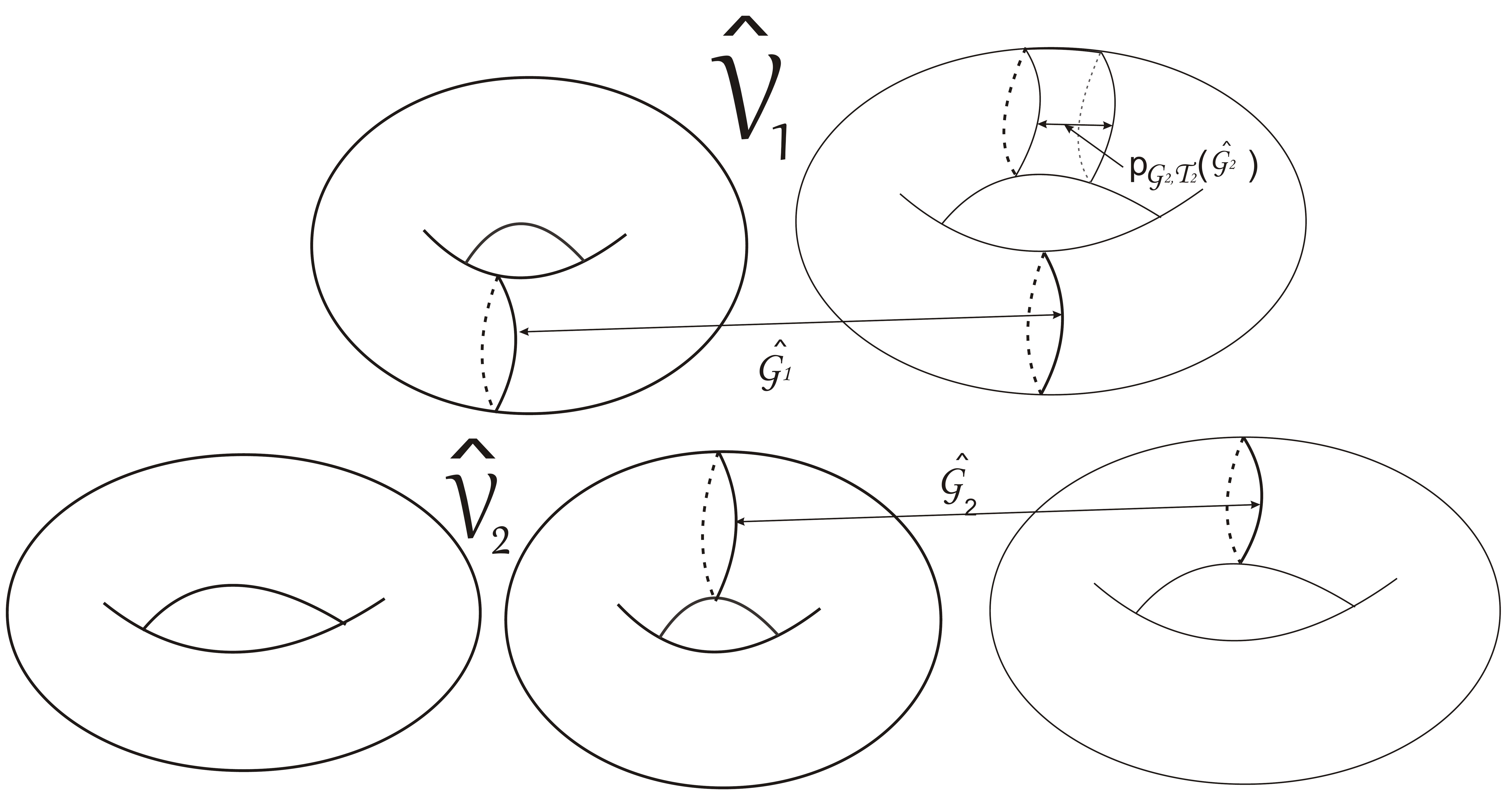}}\caption{\small {The collection $(\hat{\mathcal V}_1,\eta_{_{\hat{\mathcal V}_1}},\hat{\mathcal G}_0,\mathcal T_0),(\hat{\mathcal V}_2,\eta_{_{\hat{\mathcal V}_2}},\hat{\mathcal G}_1,\mathcal T_1)$ is not decomposed because $\hat{\mathcal G}_{1}$ has empty intersection with $p_{_{\hat{\mathcal G}_2,\mathcal T_2}}(\hat{\mathcal G}_2)$.}} \label{1contr}
\end{figure}

\begin{figure}[h!] 
\centerline{\includegraphics[width=12cm,height=9cm]{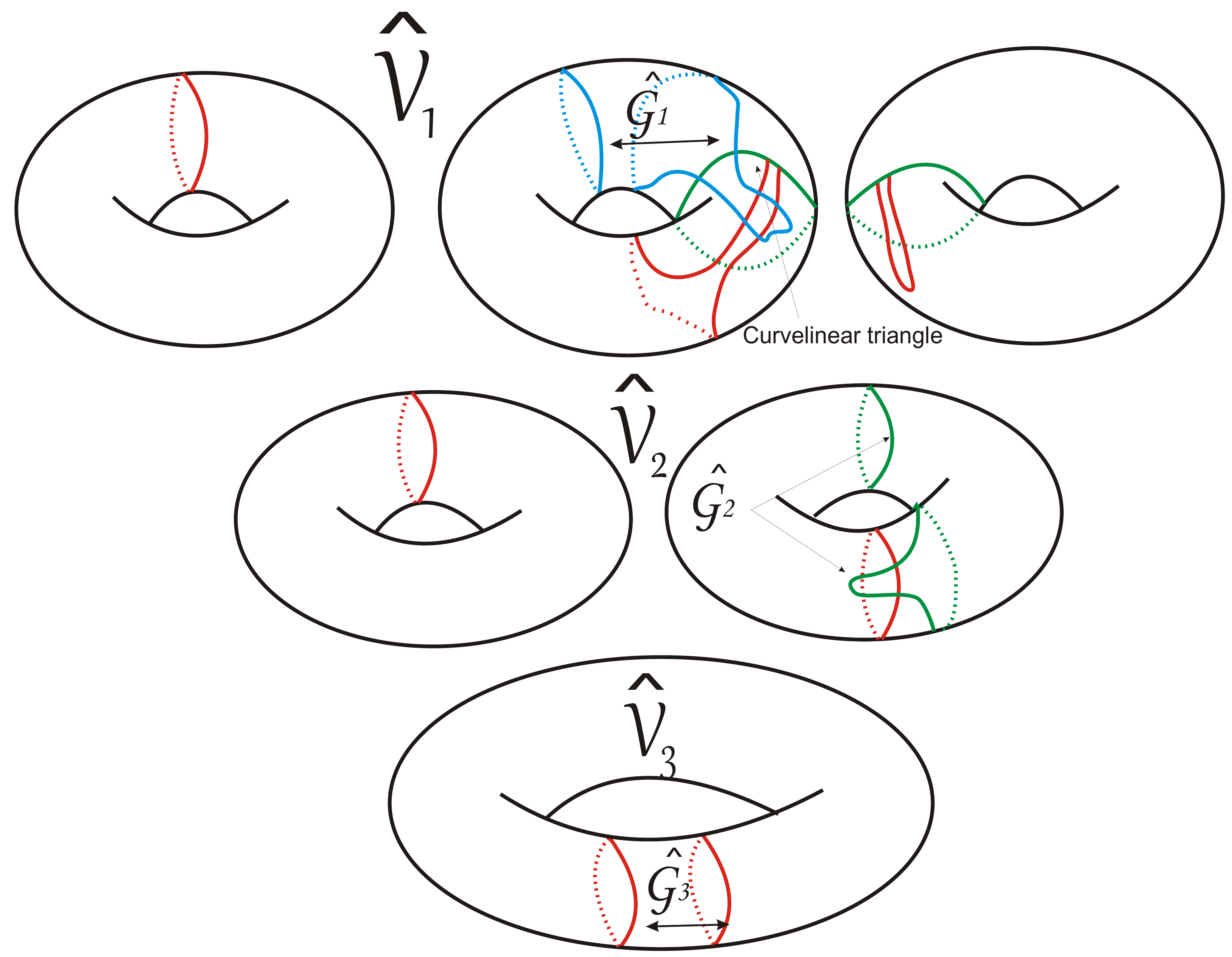}}\caption{\small Here there  is a curvelinear triangle and so  this situation
does not represent a abstract decomposed scheme according to  Definition~\ref{def:allowed}.} \label{3contr}
\end{figure}

\paragraph{The notion of an abstract decomposed scheme.}
\begin{defi}\label{def:allowed}
We say that a sequence of collections $(\hat{\mathcal V}_1,\eta_{_{\hat{\mathcal V}_1}},\hat{\mathcal G}_1,\mathcal T_1),\dots,(\hat{\mathcal V}_{n},\eta_{_{\hat{\mathcal V}_{n}}},\hat{\mathcal G}_{n},\mathcal T_{n})$ is an  {\em abstract  decomposed scheme} if for every $i\in\{2,\dots,n\}$ we have:     
\begin{itemize}[topsep=1pt,itemsep=1pt,parsep=1pt]
\item[1)] $({{\hat{\mathcal {V}_i}}},\eta_{_{\hat{\mathcal {V}_i}}})_{{\hat{\mathcal {G}_i},\mathcal T_i}}=(\hat{\mathcal V}_{i-1},\eta_{_{\hat{\mathcal V}_{i-1}}})$; 
\item[2)]  $\hat{\mathcal G}_{i-1}$ is transversal to $\bigcup\limits_{j=0}^{n-i}p_{_{\hat{\mathcal G}_{i},\mathcal T_{i}}}\circ\dots\circ p_{_{\hat{\mathcal G}_{i+j},\mathcal T_{i+j}}}(\hat{\mathcal G}_{i+j})$ and 
every component $\hat\Gamma_{i-1,\ell}$ of $\hat{\mathcal G}_{i-1}$ has nonempty intersection with $p_{_{\hat{\mathcal G}_i,\mathcal T_i}}(\hat{\mathcal G}_i)$ ({see Figure \ref{1contr} where this condition is  not satisfied});
\item[3)] the set $\hat{\mathcal V}_{i-1}\setminus int\left(\hat{\mathcal G}_{i}\cup\bigcup\limits_{j=0}^{n-i}p_{_{\hat{\mathcal G}_{i-1},\mathcal T_{i}}}\circ\dots\circ p_{_{\hat{\mathcal G}_{i+j},\mathcal T_{i+j}}}(\hat{\mathcal G}_{i+j})\right)$ does not contain curvilinear triangles as the connected components ({see Figure \ref{3contr} where this condition is failed}). 
\end{itemize}
\end{defi}  

\medskip

Notice that that for arbitrary  diffeomorphism $f\in MS(M^2)$ the property 1) $({{\hat{\mathcal {V}_i}}},\eta_{_{\hat{\mathcal {V}_i}}})_{{\hat{\mathcal {G}_i},\mathcal T_i}}=(\hat{\mathcal V}_{i-1},\eta_{_{\hat{\mathcal V}_{i-1}}})$ is clear from the discussion about cutting and pasting operation above. Property 2) corresponds to the fact that $f$ is MS and 3) to the maximality of the linearising neighborhoods. Thus we get the following fact.

\begin{lemm}\label{lem:real-vs-abstract}
Each {diffeomorphism} $f\in MS(M^2)$ induces an abstract decomposed scheme.
\end{lemm}

\subsection{Realisability of any abstract decomposed scheme:  statement of Theorem~\ref{t.realisation}} 
\label{subsection:statementTheorem3}
For an abstract decomposed collection $(\hat{\mathcal V}_1,\eta_{_{\hat{\mathcal V}_1}},\hat{\mathcal G}_1,\mathcal T_1),\dots,(\hat{\mathcal V}_{n},\eta_{_{\hat{\mathcal V}_{n}}},\hat{\mathcal N}_{n},\mathcal T_{n})$ let 
$$\hat{\mathcal V}_{0}={(\hat{\mathcal {V}}_{1})}_{\hat{\mathcal {G}}_{1},\mathcal T_{1}},$$ 
$$\lambda_0=\sum\limits_{\hat V_j\subset\hat{\mathcal V}_0}m_{j},$$ 
$$\quad\quad\quad\quad\quad\quad\quad\quad\quad \lambda_1=\sum\limits_{i=1}^{n-1}l_i \mbox{ where } l_i=\frac12\sum\limits_{\hat \gamma_{j}\subset\hat{\mathcal G}_i}m_{j}k_j,$$  $$\lambda_2=\sum\limits_{\hat V_{j}\subset \hat{\mathcal V}_{n}}m_{j}.$$ 
 
\begin{defi} \label{opr} The collection $$S=(\hat{\mathcal V},\eta_{_{\hat{\mathcal V}}},\bigcup\limits_{i=1}^{n}\{\hat{U}_{i,\ell},\ell=1,\dots,r_i\})$$ is called a {\em realizable scheme} if there is an abstract decomposed scheme $(\hat{\mathcal V}_1,\eta_{_{\hat{\mathcal V}_1}},\hat{\mathcal G}_1,\mathcal T_1),\dots,(\hat{\mathcal V}_{n},\eta_{_{\hat{\mathcal V}_{n}}},\hat{\mathcal G}_{n},\mathcal T_{n})$ such that $\la_0-\la_1+\la_2$ is even, $\la_0-\la_1+\la_2\leq 2$ and:

1) $\hat{\mathcal V}=\hat{\mathcal V}_0$ and $\eta_{_{\hat{\mathcal V}}}=\eta_{_{\hat{\mathcal V}_0}}$;

2) $\hat{\mathcal G}_i$ consists of $r_i$ components $\hat{\Gamma}_{i,\ell},\ell=1,\dots,r_i$ such that $\hat U_{i,\ell}$ is a tubular neighborhood of $p_{_{\hat{\mathcal G}_{1},\mathcal T_{1}}}\circ\dots\circ p_{_{\hat{\mathcal G}_{i},\mathcal T_{i}}}(\hat\Gamma _{i,\ell})$ for $i\in\{1,\dots,n\}$.
\end{defi}

Notice that for the decomposed scheme of a diffeomorphism $f\in MS(M^2)$ we have that $\lambda_0$ is a number of the sink points, $\lambda_1$ is a number of the saddle points, $\lambda_2$ is a number of the source points. Thus $\lambda_0-\lambda_1+\lambda_2$ is the euler characteristic of $M^2$ equals $2-2g$, where $g$ is the genus of the surface $M^2$. Moreover, for each $i\in\{1,\dots,beh(f)\}$ the set $\hat{\mathcal U}_i$ is a  tubular neighborhood of the set $p_{_{\hat{\mathcal G}_{1},\mathcal T_{1}}}\circ\dots\circ p_{_{\hat{\mathcal G}_{i},\mathcal T_{i}}}(\hat{\mathcal G}_{i})$. Thus we have the following fact.

\begin{lemm} The decomposed scheme associated to a diffeomorphism $f\in MS(M^2)$ is an realizable in the sense above.
\end{lemm}

\begin{theo} \label{t.realisation} For any realizable decomposed scheme $S$ there is a diffeomorphism $f_{S}\in MS(M^2)$ whose scheme is equivalent to $S$ and {so that the euler characteristic of the orientable surface is equal to $\lambda_0-\lambda_1+\lambda_2$}. 
\end{theo}

{So MS orientation preserving surface diffeomorphisms can be fully classified by decomposed schemes (which automatically satisfy 
the properties of a realizable abstract decomposed scheme), and vice versa each such abstract decomposed scheme 
corresponds to a unique conjugacy class of a diffeomorphisms.}

\subsection{Examples of realizable abstract decomposed schemes and their realisations}

\begin{figure}[h!]
\centerline{\includegraphics[width=11cm,height=8.5cm]{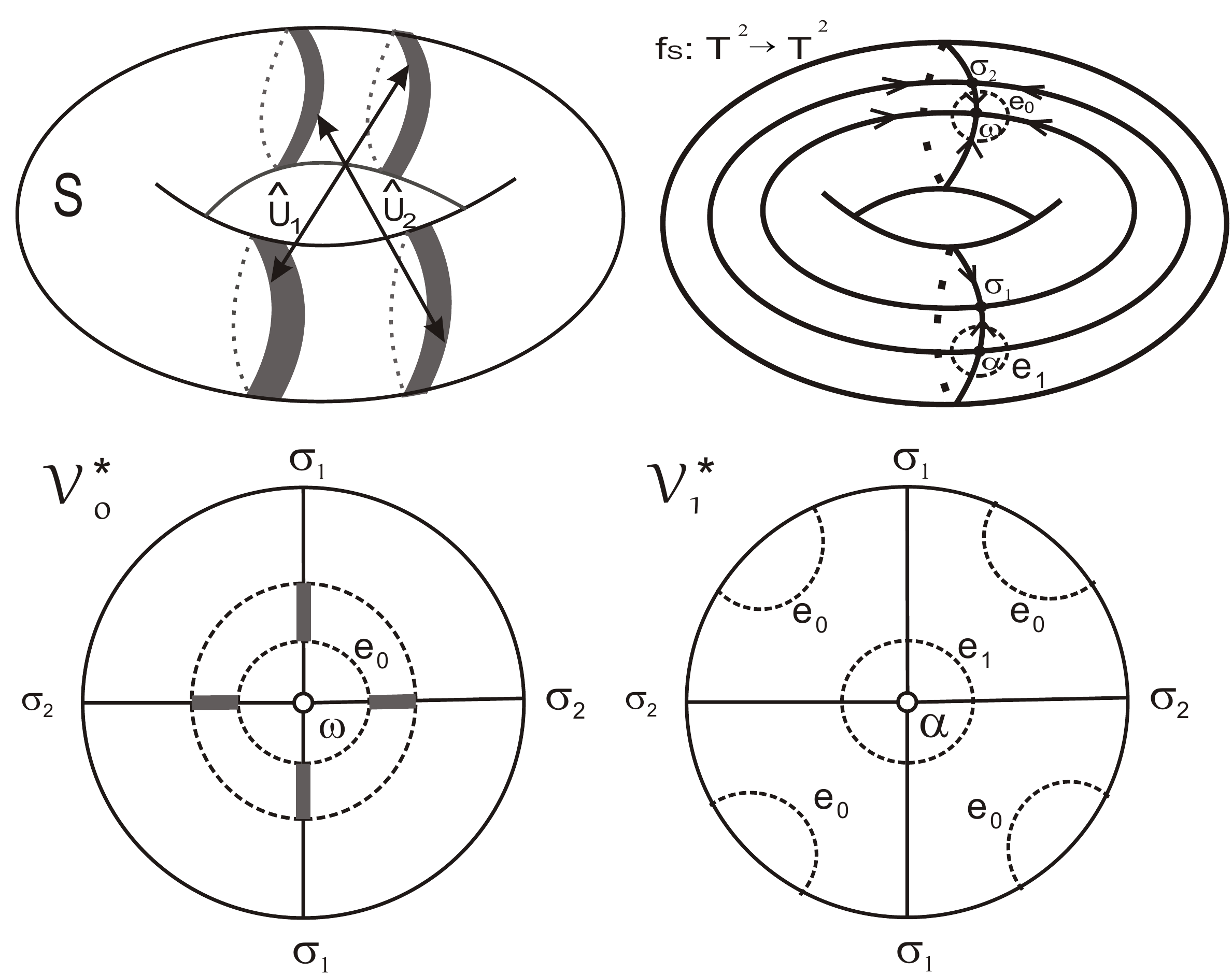}}\caption{\small Schemes $S$ and a phase portrait of the diffeomorphism $f_S\in MS(M^2)$.} \label{+schere}
\end{figure}

\begin{example}[Figure \ref{+schere}]\label{example:+schere}
Consider the scheme $S$  from Figure~\ref{+schere}. Here $\hat{\mathcal V}=\mathbb T^2$, $n=1$, 
$m_{\hat{\mathcal V}}=1$ and $\hat{U}_{i}=\hat N_i,i=1,2$ consist of two annuli on the torus $\hat{\mathcal V}$ such that $\eta_{_{\hat{\mathcal V}}}(i_{_{\hat {N}_i}*}(\pi_1(\hat{N}_i)))=\mathbb Z$. The set $(\hat{\mathcal V})_{\hat{\mathcal N}}$ consists of one torus. Thus $\lambda_0=1$, $\lambda_1=2$ and $\lambda_2=1$,  $M_S=\mathbb T^2$, $\Si_0=\{\omega\}$, $\Si_1=\{\sigma_1,\sigma_2\}$, $\nu_{\si_1}=\nu_{\si_2}=+$ and $\Si_2=\{\alpha\}$. The interior of the square corresponds to the set ${\mathcal V}^*$, where $f_S$ is the contraction to $\omega$.  The torus is obtained by identifying the boundary of the circle as suggested by the labelling of $\sigma_i$ and taking of the connected sum as suggested by the labelling of $e_i$.  
\end{example}

\begin{figure}[h!]
\centerline{\includegraphics[width=14cm,height=9cm]{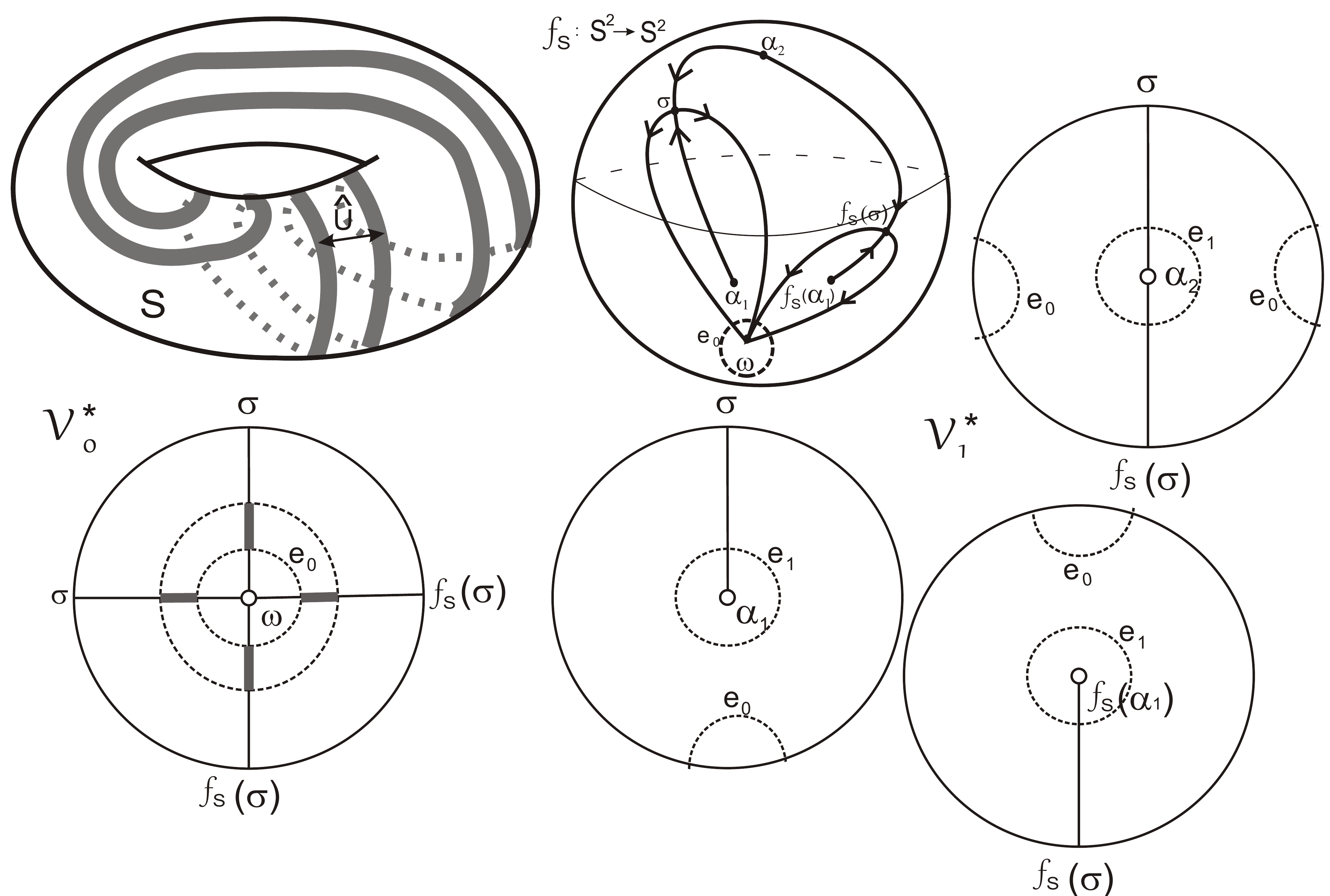}}\caption{\small Schemes $S$ and a phase portrates of the diffeomorphisms $f_S\in MS(M^2)$.} \label{s12}
\end{figure}

\begin{example}[Figure \ref{s12}]\label{example:+s12}
Consider the scheme $S$  from Figure~\ref{s12}. Here $\hat{\mathcal V}=\mathbb T^2$, $n=1$, 
$m_{\hat{\mathcal V}}=1$ and $\hat{U}=\hat N$ consist of two annuli on the torus $\hat{\mathcal V}$ such that $\eta_{_{\hat{\mathcal V}}}(i_{_{\hat {N}}*}(\pi_1(\hat{N})))=2\mathbb Z$. The set $(\hat{\mathcal V})_{\hat{N}}$ consists of two tori $\hat V_{1,1}$ and $\hat V_{1,2}$ such that $m_{\hat{V}_{1,1}}=1$ and $m_{\hat{V}_{1,2}}=2$. Thus $\lambda_0=1$, $\lambda_1=2$ and $\lambda_2=3$,  $M_S=\mathbb S^2$, $\Si_0=\{\omega\}$, $\Si_1=\{\sigma,f_S(\sigma)\}$, $\nu_{\si}=+$ and $\Si_2=\{\alpha_1,f_S(\alpha_1),\alpha_2\}$. The interior of the square corresponds to the set ${\mathcal V}^*$, where $f_S$ is the contraction to $\omega$.  The ambient sphere $\mathbb S^2$ is obtained by identifying the boundary of the circle as suggested by the labelling of $\sigma_i$ and taking of the connected sum as suggested by the labelling of $e_i$.  
\end{example}

\begin{figure}[h!]
\centerline{\includegraphics[width=14cm,height=14cm]{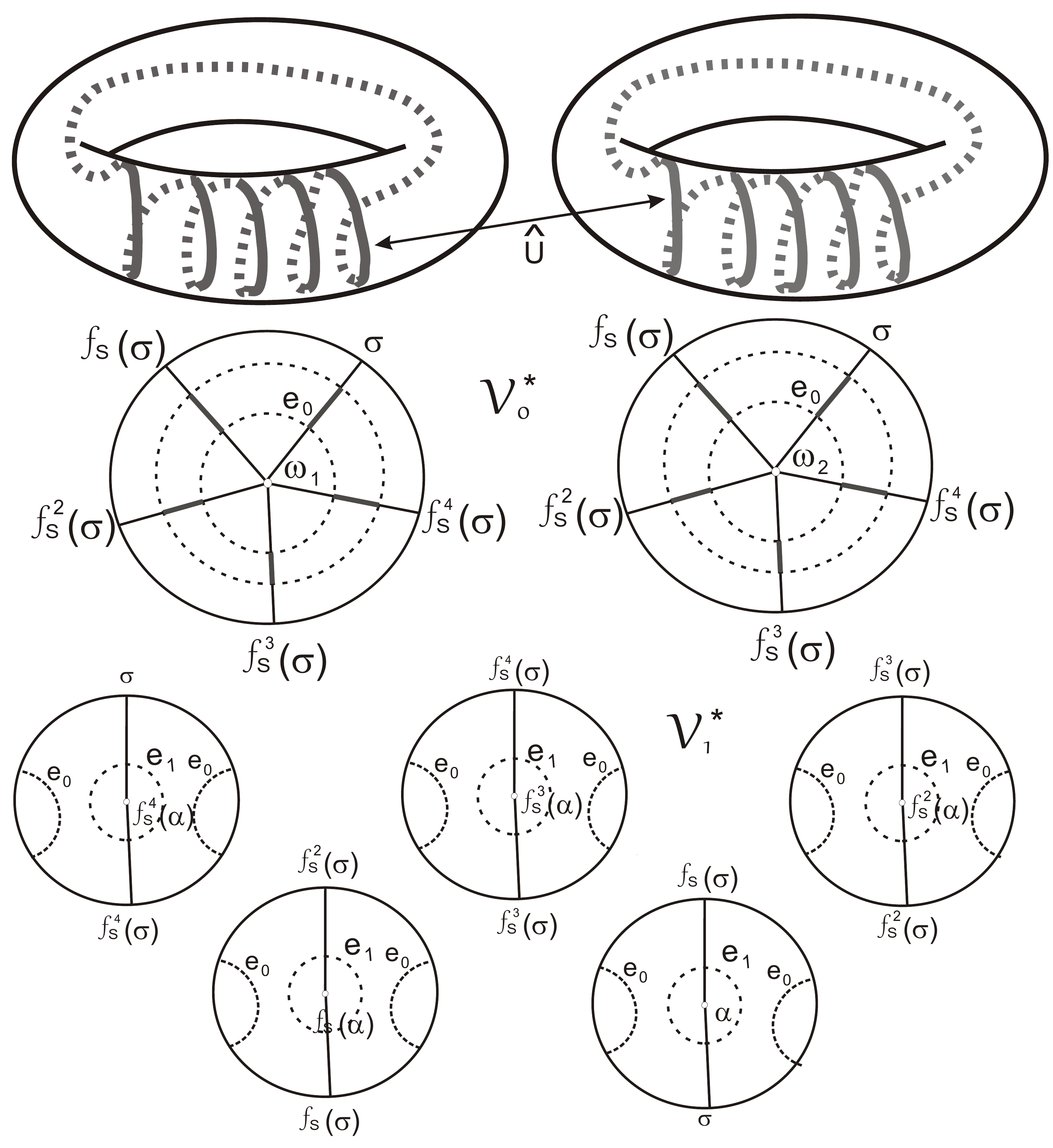}}\caption{\small Schemes $S$ and a phase portrait of the diffeomorphism $f_S\in MS(M^2)$.} \label{11}
\end{figure}

\begin{example}[Figure \ref{11}]\label{example:+11}
Consider the scheme $S$  from Figure~\ref{11}. Here $\hat{\mathcal V}=\hat V_1\cup\hat V_2$, $n=1$, 
$m_{\hat{V}_i}=1$ and $\hat{U}=\hat N$ consist of two annuli $\hat A_i,i=1,2$ on the tori $\hat{V}_i$ such that $\eta_{_{\hat{V}_i}}(i_{_{\hat {A}_i}*}(\pi_1(\hat{A}_i)))=5\mathbb Z$. The set $(\hat{\mathcal V})_{\hat{N}}$ consists of one torus $\hat V_{1,1}$ such that $m_{\hat{V}_{1,1}}=5$. Thus $\lambda_0=2$, $\lambda_1=5$ and $\lambda_2=5$,  $M_S=\mathbb S^2$, $\Si_0=\{\omega_1,\om_2\}$, $\Si_1=\{\sigma,f_S(\sigma),f^2_S(\si),f^3_S(\si),f^4_S(\si)\}$, $\nu_{\si_i}=+1,i=1,2,3,4,5$ and $\Si_2=\{\alpha,f_S(\alpha),f^2_S(\alpha),f^3_S(\al),f^4_S(\al)\}$. The interior of the punctured discs in the middle of the picture corresponds to the set ${\mathcal V}^*$, where $f_S$ is the contractions to $\omega_1,\omega_2$ in the composition with the $1/5$ part of the revolution around $\omega_1,\omega_2$, accordingly. The ambient sphere $\mathbb S^2$ is obtained by identifying the boundary of the circle as suggested by the labelling of $\sigma_i$ and taking of the connected sum as suggested by the labelling of $e_i$.  
\end{example}

\begin{figure}[h!]
\centerline{\includegraphics[width=14cm,height=12cm]{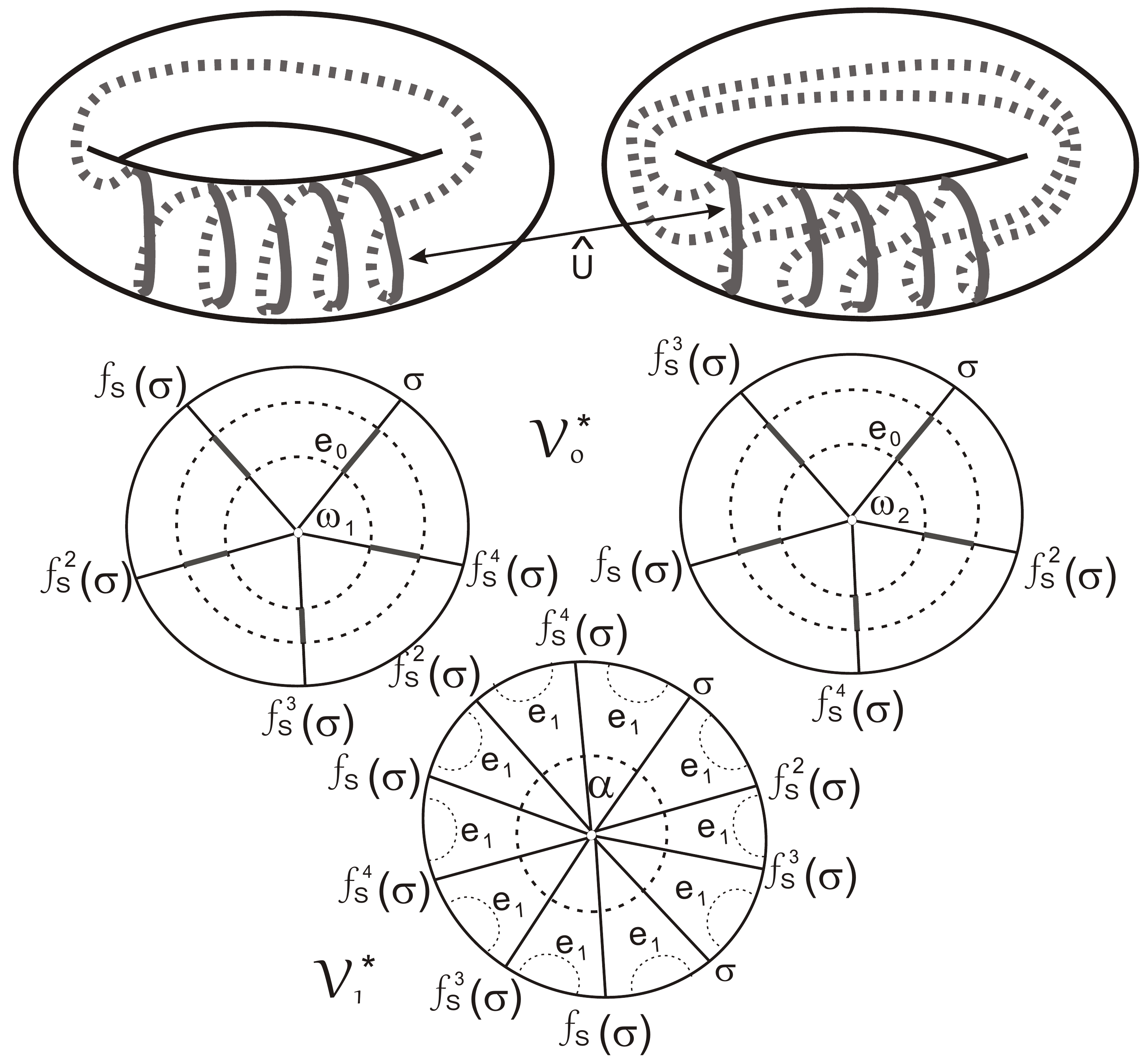}}\caption{\small Schemes $S$ and a phase portrait of the diffeomorphisms $f_S\in MS(M^2)$.   } \label{12}
\end{figure}

\begin{example}[Figure \ref{12}]\label{example:+12}
Consider the scheme $S$  from Figure~\ref{12}. Here $\hat{\mathcal V}=\hat V_1\cup\hat V_2$, $n=1$, 
$m_{\hat{V}_i}=1$ and $\hat{U}=\hat N$ consist of two annuli $\hat A_i,i=1,2$ on the tori $\hat{V}_i$ such that $\eta_{_{\hat{V}_i}}(i_{_{\hat {A}_i}*}(\pi_1(\hat{A}_i)))=5\mathbb Z$. The set $(\hat{\mathcal V})_{\hat{N}}$ consists of one torus $\hat V_{1,1}$ such that $m_{\hat{V}_{1,1}}=5$. Thus $\lambda_0=2$, $\lambda_1=5$ and $\lambda_2=1$,  $M_S$ is the surface of the genus $2$, $\Si_0=\{\omega_1,\om_2\}$, $\Si_1=\{\sigma,f_S(\sigma),f^2_S(\si),f^3_S(\si),f^4_S(\si)\}$, $\nu_{\si_i}=+1,i=1,2,3,4,5$ and $\Si_2=\{\alpha\}$. The interior of two punctured discs in the middle of picture corresponds to the set ${\mathcal V}^*$, where $f_S$ is the contractions to $\omega_1,\omega_2$ in the composition with the $1/5,2/5$ part of the revolution around $\omega_1,\omega_2$, accordingly. The ambient surface of the genus $2$ is obtained by identifying the boundary of the circle as suggested by the labelling of $\sigma_i$ and taking of the connected sum as suggested by the labelling of $e_i$.  
\end{example}

\begin{figure}[h!]
\centerline{\includegraphics[width=16cm,height=12cm]{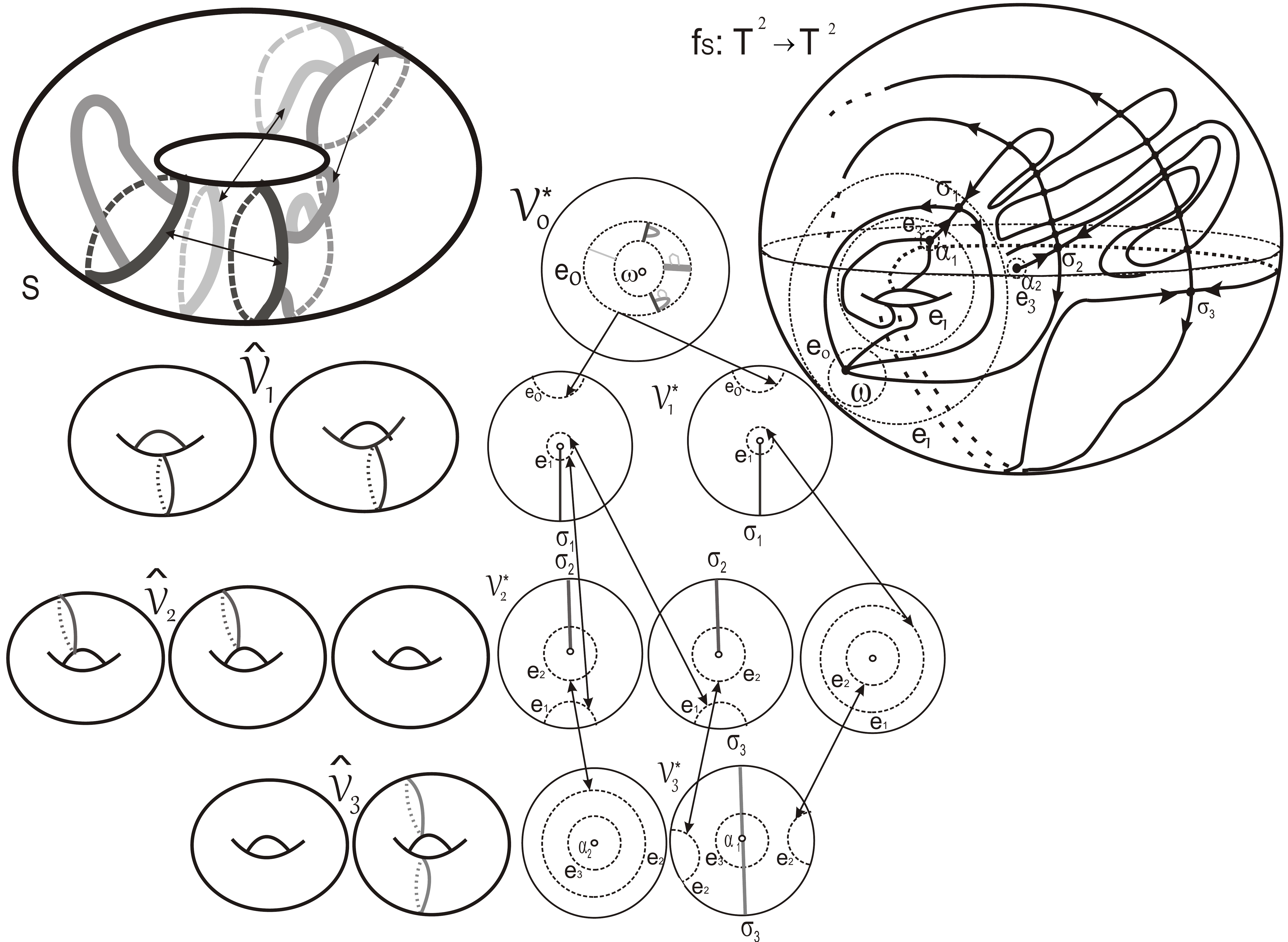}}\caption{\small Schemes $S$, the decomposed sequence and a phase portrates of the diffeomorphisms $f_S\in MS(M^2)$.} \label{s14t}
\end{figure}

\begin{example}[Figure \ref{s14t}]\label{example:+s14t}
Consider the scheme $S$ from Figure~\ref{s14t}. Here $\hat{\mathcal V}=\hat V$, $m_{\hat{V}}=1$ and $n=3$. Let us describe the decomposed sequence.
\begin{itemize}
\item The set $\hat V_3$ consists of two tori, one of which contains two knots $\hat{\mathcal G}_3$ such that $\eta_{_{\hat{V}_3}}(i_{_{\hat {\mathcal G}_{3}}*}(\pi_1(\hat{\mathcal G}_{3})))=\mathbb Z$. 
\item The set $\hat{\mathcal{V}}_2=(\hat{\mathcal V}_3)_{\hat{\mathcal G}_3}$ consists of three tori $\hat V_{1,1}$, $\hat V_{1,2}$ such that $m_{\hat{V}_{1,j}}=1$. The set $\hat{\mathcal N_1}$ consists of two annuli $\hat A_{0,i},i=1,2$ on the torus $\hat{V}_{1,2}$ such that $\eta_{_{\hat{V}_{1,2}}}(i_{_{\hat {A}_{1,i}}*}(\pi_1(\hat{A}_{1,i})))=\mathbb Z$. 
\item The set $\hat{\mathcal{V}}_2=(\hat{\mathcal V}_1)_{\hat{\mathcal N}_1}$ consists of three tori $\hat V_{2,1}$, $\hat V_{2,2}$, $\hat V_{2,3}$ such that $m_{\hat{V}_{2,j}}=1$. The set $\hat{\mathcal N_2}$ consists of two annuli $\hat A_{2,1},\hat A_{2,2}$ on the tori $\hat{V}_{2,2},\hat{V}_{2,3}$ such that $\eta_{_{\hat{V}_{2,j}}}(i_{_{\hat {A}_{2,i}}*}(\pi_1(\hat{A}_{2,i})))=\mathbb Z$. 
\item The set $\hat{\mathcal{V}}_3=(\hat{\mathcal V}_2)_{\hat{\mathcal N}_2}$ consists of two tori $\hat V_{3,1}$, $\hat V_{3,2}$ such that $m_{\hat{V}_{3,j}}=1$.
\end{itemize}
Thus $\lambda_0=1$, $\lambda_1=5$ and $\lambda_2=2$,  $M_S=\mathbb T^2$, $\Si_0=\{\omega\}$, $\Si_1=\{\sigma_1\}$, $\Si_2=\{\sigma_2\}$, $\Si_3=\{\sigma_3\}$,$\nu_{\si_i}=+,i=1,2,3$ and $\Si_4=\{\alpha_1,\alpha_2\}$. The ambient torus $\mathbb T^2$ is obtained as a connected sum three copies of 2-sphere as suggested by the labelling of saddle points $\si_i$ and the equators $e_i$.
\end{example}

\section{Proof of Theorem~\ref{add}: maximal systems of neighbourhoods}\label{maxx}
In this section we will proof Theorem \ref{add}. Let us give some notation before. Recall that we divide 
the set $\Sigma$ of the saddle points by parts $\Sigma=\Sigma_1\cup\dots\cup\Sigma_{beh(f)-1}$. For each $i\in\{0,\dots,beh(f)-1\}$ let us set $$A_i=\bigcup\limits_{j=0}^{i}W^u_j,~~R_i=\bigcup\limits_{j=i+1}^{beh(f)}W^s_j,~~\mathcal V_i=M^2\setminus(A_i\cup R_i).$$ Observe that $A_i$ is an attractor, $R_i$ is a repeller of $f$ and $f$ acts freely on $\mathcal V_i $. Set $\hat{\mathcal V}_i=\mathcal V_i/f$ and denote the natural projection by $$p_{_{i}}:\mathcal V_{i}\to\hat{\mathcal V}_{i}.$$ Notice that $\hat{\mathcal V}_0=\hat{\mathcal V}_f$ and $p_0=p_{_f}$. It is proved, for example, in book \cite{grin} that each connected component of $\hat{\mathcal V}_i$ is a closed orientable surface and $p_i$ is a cover. We introduce the following notations:

\begin{itemize} 
\item[-] for $j\in\{1,\dots,beh(f)-1\},~k\in\{0,\dots,beh(f)-1\}$ let $\hat W^s_{j,k}=p_k(W^{s}_{j}\cap{\mathcal V}_k)$, $\hat W^u_{j,k}=p_k(W^{u}_{j}\cap{\mathcal V}_k)$; 
\item[-]  $L^u=\bigcup\limits_{i=1}^{beh(f)-1}W^u_i,~L^s=\bigcup\limits_{i=1}^{beh(f)-1}W^s_i$,\,
$\hat{L}^u_i=p_i(L^u)$, $\hat{L}^s_i=p_i(L^s)$. 
\end{itemize} 

Before a proof of Theorem \ref{add} we introduce a more strong than $u$-compatibility property for a neighborhood of a saddle point. 

\begin{figure}[h!]
\centerline{\includegraphics[width=5.5cm,height=5.5cm]{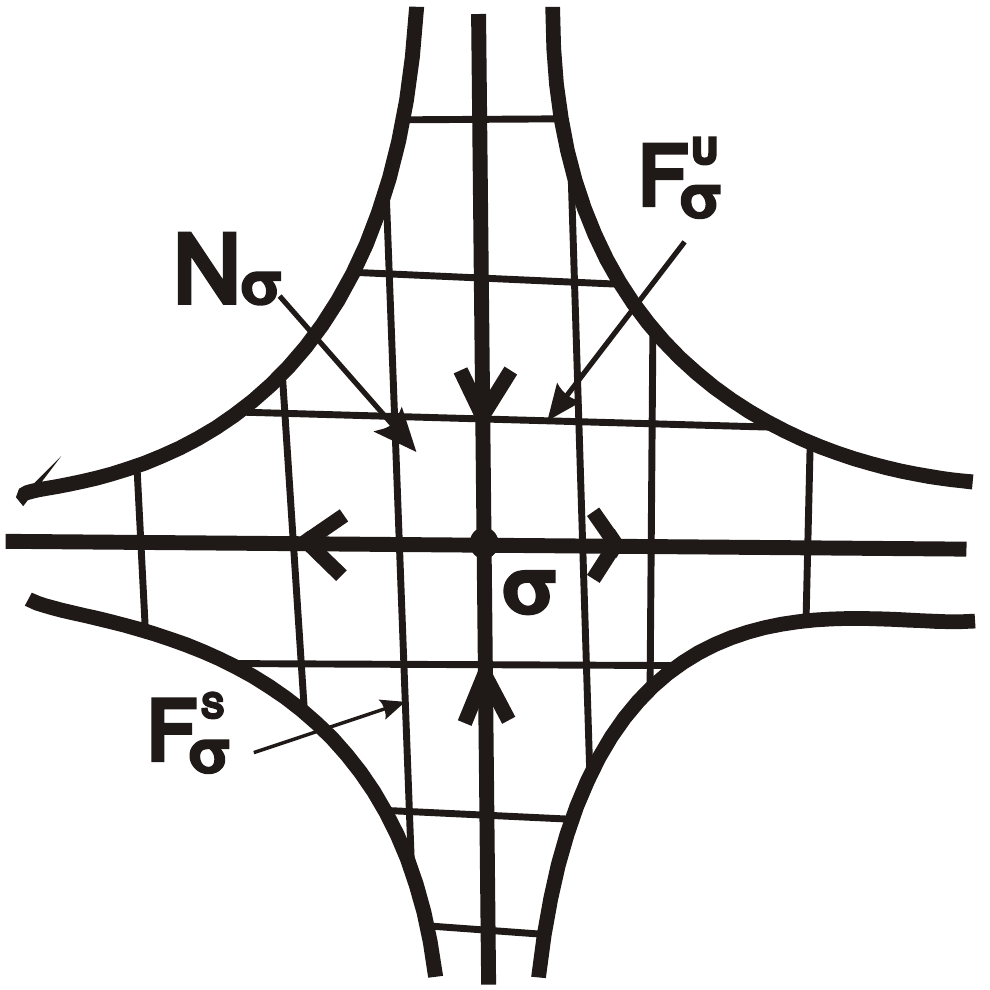}}\caption{\small Linearizing neighborhood} \label{5+}
\end{figure}

\begin{defi}[The linearizing neighborhood] \label{adus} Let $\sigma$ be a saddle periodic point for $f$. A neighborhood $\mathtt{N}_\sigma$ of the point $\sigma$ with a one-dimensional foliation ${\mathtt{F}}^u_{\sigma}$ containing $W^u_\sigma$ as a leaf and a one-dimensional foliation ${\mathtt{F}}^s_{\sigma}$ containing $W^s_\sigma$ as a leaf, is called \emph{linearizable} if there is a homeomorphism ${\mu}_\sigma:\mathtt N_\sigma\to{N}$ which conjugates the diffeomorphism $f^{k_\sigma}\vert_{\mathtt{N}_\sigma}$ to the  canonical diffeomorphism $a_{\nu_\sigma}\vert_{{N}}$ and sends leaves of the foliation ${\mathtt{F}}^u_{\sigma}$ to leaves the foliation $F^u$, also sends leaves of the foliation ${\mathtt{F}}^s_{\sigma}$ to leaves the foliation $F^s$ $($see Figure \ref{5+}$)$.
\end{defi}

For every point $x\in \mathtt{N}_{\sigma}$ denote by $\mathtt{F}^u_{\si,x},\,{\mathtt{F}}^s_{\sigma,x}$ the unique leave of the foliation  $\mathtt{F}^u_\si,\,{\mathtt{F}}^s_{\sigma}$, accordingly,  passing through the point $x$.

\begin{figure}[h!]
\centerline{\includegraphics[width=11cm,height=8cm]{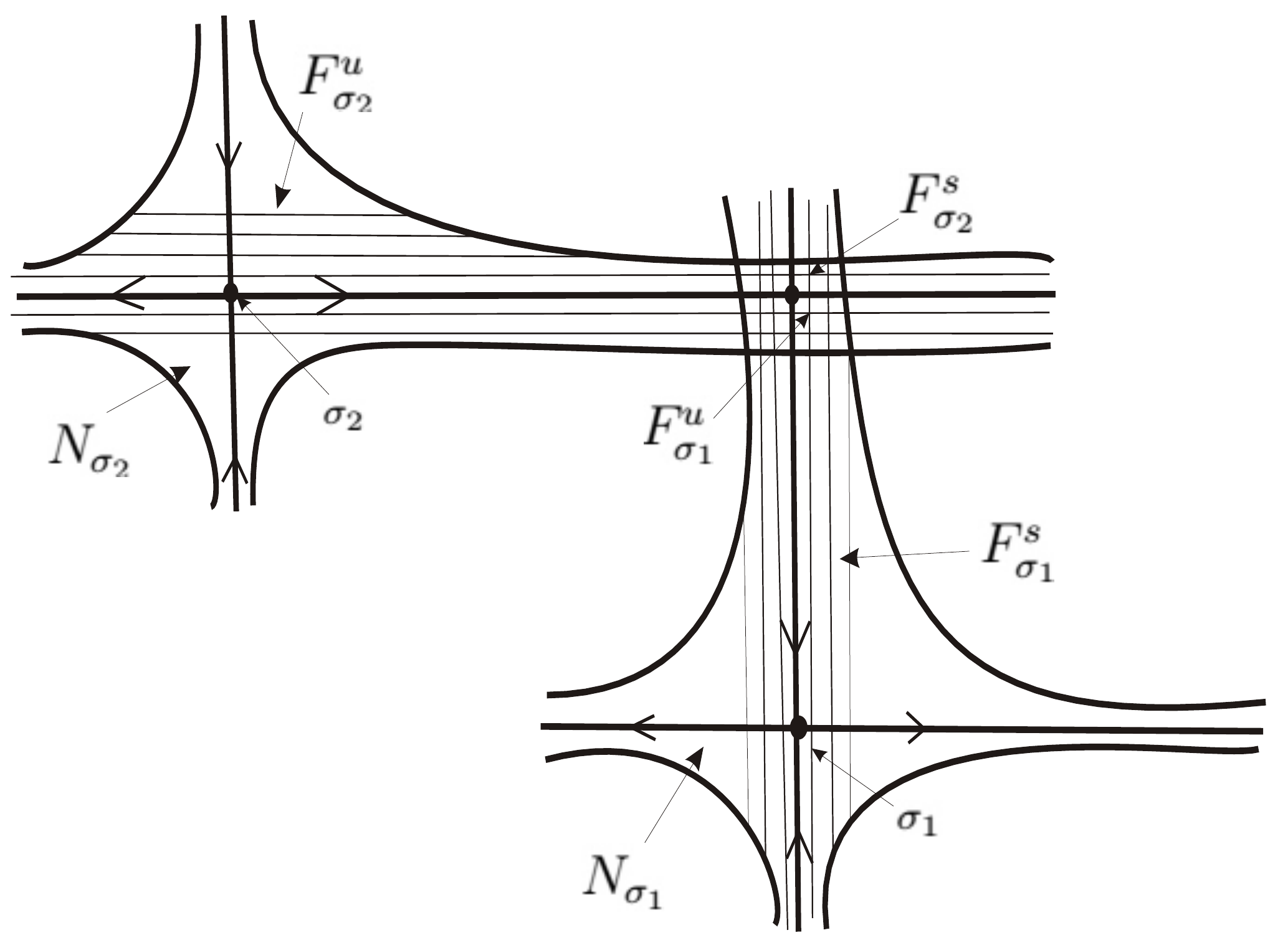}}\caption{\small A compatible system of neighbourhoods} \label{a+}
\end{figure}

\begin{defi}[The compatible system of neighbourhoods] \label{dous} An $f$-invariant collection $\mathtt{N}_f$ of linearizable neighborhoods $\mathtt{N}_\sigma$ of all saddle points $\sigma\in\Sigma$ is called \emph{compatible} if the following properties are hold:

1) $\mu_\sigma^{-1}(\partial{{N}})$ does not contain heteroclinic points for any $\sigma\in\Sigma$; 

2) if $W^{s}_{{\sigma_1}}\cap W^{u}_{{\sigma_2}}=\emptyset$ and $W^{u}_{{\sigma_1}}\cap W^{s}_{{\sigma_2}}=\emptyset$ for $\sigma_1,\sigma_2\in\Sigma$ then ${\mathtt{N}}_{{\sigma_1}}\cap{\mathtt{N}}_{{\sigma_2}}=\emptyset$;

3) if $W^s_{\sigma_1}\cap W^u_{\sigma_2}\neq\emptyset$ for $\sigma_1,\sigma_2\in\Sigma$ then $({\mathtt{F}}^u_{{\sigma_1},x}\cap {\mathtt{N}}_{{\sigma_2}})\subset{\mathtt{F}}^u_{{\sigma_2},x}$ and $({\mathtt{F}}^s_{{\sigma_2},x}\cap {\mathtt{N}}_{{\sigma_1}})\subset{\mathtt{F}}^s_{{\sigma_1},x}$ for  $x\in(\mathtt{N}_{{\sigma_1}}\cap \mathtt{N}_{{\sigma_2}})$ $($see Figure \ref{a+}$)$.
\end{defi}

\begin{lemm}\label{comp} For every diffeomorphism $f\in MS(M^2)$ there is a compatible system of  neighbourhoods.
\end{lemm}
\begin{demo} The proof consists of three steps. 

{\bf Step 1.} Here, we prove the existence of $f$-invariant neighbourhoods $Q^{s}_{1},\dots,Q^{s}_{beh(f)-1}$ of the sets $\Sigma_{1},\dots,\Sigma_{beh(f)-1}$ respectively,  equipped with  one-dimensional $f$-invariant  foliations ${\mathtt{F}}^{u}_{1},\dots,{\mathtt{F}}^{u}_{beh(f)-1}$ whose leaves are smooth
such that the following properties hold for each $i\in\{1,\dots,beh(f)-1\}$:

{\rm (i)} the unstable manifolds $W^u_i$ are leaves of the foliation ${\mathtt{F}}^{u}_{i}$ and 
each leaf of the foliation ${\mathtt{F}}^{u}_{i}$ is transverse to  $L^s$;
 
{\rm (ii)} for any  $1\leq i <k\leq beh(f)-1$ and $x\in Q^s_{i}\cap Q^s_{k}$, we have the inclusion 
$\mathtt{F}^u_{k,x}\cap {Q}^s_{i}\subset \mathtt{F}^u_{i,x}$. 

Let us prove this by a decreasing induction on  $i$ from $i=beh(f)-1$ to $i=1$.

For $i=beh(f)-1$, it follows from the definition of $\mathcal {V}_{beh(f)-1}$ that $(W^s_{beh(f)-1}\smallsetminus \Sigma_{beh(f)-1})\subset\mathcal{V}_{beh(f)-1}$. Since $f$ acts freely and properly on $W^s_{beh(f)-1}$, the quotient $\hat W^{s}_{{beh(f)-1},{beh(f)-1}}$ is a smooth submanifold of $\hat{\mathcal V}_{{beh(f)-1}}$; it consists of finite number circles. The lamination $\hat L^s_{beh(f)-1}$ accumulates on $\hat W^{s}_{{beh(f)-1},{beh(f)-1}}$. Choose a closed tubular neighbourhood ${\hat N}^{s}_{{beh(f)-1}}$ of  $\hat W^{s}_{{beh(f)-1},{beh(f)-1}}$ in $\hat{\mathcal V}_{beh(f)-1}$; denote its projection by $\pi^{u}_{beh(f)-1}:{\hat N}^{s}_{beh(f)-1}\to \hat W^{s}_{{beh(f)-1},{beh(f)-1}}$. Its fibres form a segment foliation $\{{d}^{u}_{{beh(f)-1}, x}\mid x\in \hat W^{s}_{{beh(f)-1},{beh(f)-1}}\}$. Since $\hat L^s_{beh(f)-1}$ is a $C^1$-lamination  containing $\hat W^{s}_{{beh(f)-1},{beh(f)-1}}$ as a leaf, if the tube ${\hat N}^{s}_{beh(f)-1}$ is small enough, the fibres are  transverse  to $\hat L^s_{beh(f)-1}$.
 
Set $Q^{s}_{beh(f)-1}:=p_{beh(f)-1}^{-1}({\hat N}^{s}_{beh(f)-1})\cup W^u_{beh(f)-1}$. This is a subset of $M^2$ which carries a foliation $\mathtt{F}^u_{beh(f)-1}$ defined by pullpacking the fibres of $\pi^u_{beh(f)-1} $ and by adding $W^u_{beh(f)-1}$ as extra leaves. This is the wanted foliation satisfying (i) and (ii)  for $i={beh(f)-1}$. Notice that the leaves of $\mathtt{F}^u_{beh(f)-1}$ are smooth.

For the induction we assume the construction is done for every $j>i$ and we have to construct an $f$-invariant neighborhood $Q^{s}_{i}$ of the saddle points from $\Sigma_{i}$ carrying an $f$-invariant foliation ${\mathtt{F}}^{u}_{i}$ satisfying (i) and (ii). Moreover, by genericity the boundary  $\partial Q^s_{j}$, $j>i$, is assumed to avoid all heteroclinic points.

For  $j>i$, let $\hat{Q}^{s}_{j,i}:=p_{i}(Q^{s}_{j}\cap\mathcal{V}_{i})$ and $\hat{{\mathtt{F}}}^{u}_{j,i}:=p_{_{i}}({\mathtt{F}}^{u}_{j}\cap\mathcal{V}_{i})$. For the same reason as in the case $i={beh(f)-1}$, the set $\hat W^{s}_{i,i}$ is a smooth submanifold of $\hat{\mathcal V}_{i}$ consisting of circles. Choose a tubular neighbourhood  $\hat N^{s}_{i}$  of $\hat W^{s}_{i,i}$ with a projection $\hat\pi^{u}_{i}:{\hat N}^{s}_{i}\to \hat W^{s}_{i,i}$ whose fibres are segments. Similarly, $\hat W^u_{i+1,i}$ is a compact submanifold, consisting of a finite number circles.

The set $\hat L^u_i$ is a compact lamination  and its intersection with $\hat W^{s}_{i,i}$ consists of a countable set of points which are the  projections of the heteroclinic points belonging to the stable manifolds  $W^{s}_{i}$. Actually,  there is a hierarchy in $\hat L^u_i\cap\hat W^{s}_{i,i}$  which we are going to describe in more details.
 
Set $H_{k}:= \hat W^u_{i+k,i}\cap \hat W^s_{i,i}$ for $k>0$. Since $\hat W^u_{i+1, i}$ is compact,
$H_1$ is  a finite set: $H_1= \{h^1_1, ... , h^1_{t(1)}\}$.  We are given neighbourhoods,  called {\it  boxes},  $B^1_\ell$, $\ell=1,..., t(1)$, about these points,   namely, the connected components of $\hat Q^s_{i+1,i}\cap \hat N^s_i$.  Due to the fact that $\partial \hat Q^s_{i+1,i}$ contains no heteroclinic points, 
$\partial \hat Q^s_{i+1,i}\cap W^{s}_{i,i}$  is isolated from $\hat L^u_i$. Therefore, if the tube 
$\hat N^s_i$ is small enough, $\hat L^u_i$ does not intersect $\partial \hat Q^s_{i+1,i}\cap  \hat N^s_i$.
Then, by shrinking $Q^s_j,\,j>i+1$ if necessary, we may guarantee that $\hat Q^s_{j,i}\cap \hat N^s_i$  is disjoint from $\partial\hat Q^s_{i+1,i}\cap \hat N^s_i$.
 
Since $\hat W^u_{i+2, i}$ accumulates on $W^u_{i+1,i}$, there are only finitely many points of $H_2$ outside of all boxes  $B^1_\ell$, $\ell=1,..., t(1)$.  Let $\bar H_2:= \{h^2_1,...,h^2_{t(2)}\}$ be this finite set. The open set $\hat Q^s_{i+2,i}$  is a neighbourhood of $\bar H_2$. The  connected  components of $\hat Q^s_{i+2,i}\cap \hat N^s_i$ which contain points of $\bar H_2$  will be the box $B^2_\ell$ for $\ell= 1,..., t(2)$. We argue with $B^2_\ell$ with respect to  $\hat L^u_i$ and the neighbourhoods $\hat Q^s_{j,i},\, j>i+1,$ in a similar manner as we do with $B^1_\ell$. And so on, until $\bar H_n$.
 
Due to the induction hypothesis, each above-mentioned box is foliated. Namely, $B^1_\ell$ is foliated by $\hat {\mathtt{F}}^u_{i+1, i}$; the box $B^2_\ell$ is foliated by $\hat {\mathtt{F}}^u_{i+2, i}$, and so on. But the leaves are not contained in fibres of $\hat N_i$; even more, not every leaf intersects $W^s_{i,i}$.  We have to correct this situation in order to  construct the foliation  ${\mathtt{F}}^u_i$ satisfying the wanted conditions (i) and (ii).

For every $j>i$, the foliation ${\mathtt{F}}^u_j$ may be extended to the boundary $\partial Q^s_j$ and a bit beyond.
Once this is done, if $\hat N^s_i$ is  enough shrunk, each leaf of $\hat {\mathtt{F}}^u_{i+k,i}$ though $x\in B^k_\ell$
intersects $\hat W^s_{i,i}$ (it is understood that the boxes are intersected with the shrunk tube without changing their names). So, we have a projection along the leaves $\pi_{k,\ell}: B^k_\ell\to\hat W^s_{i,i} $; 
but, the image of $\pi_{k,\ell}$ is larger than $B^k_\ell\cap\hat W^s_{i,i}$. Then, we choose a small enlargement $B'^k_\ell$ of $B^k_\ell$ such that   $B'^k_\ell\smallsetminus B^k_\ell$ is foliated  by $\hat {\mathtt{F}}^u_{i+k,i}$ and avoids  the lamination $\hat L^u_i$. 

On $B'^k_\ell\smallsetminus B^k_\ell$ we have two projections: one is $\hat\pi^u_i$ and the other one is $\pi_{k,\ell}$. We are going to interpolate between both using a partition of unity (we do it for $B^k_\ell$ but it is understood that it is done for all boxes). Let $\phi:\hat N^{s}_{i}\to[0,1]$ be a smooth function which equals 1 near   $B^k_\ell$ and whose support is contained in $B'^k_\ell$. Define a global $C^1$ retraction $\hat q: \hat N^{s}_{i}\to \hat W^s_{i,i}$ by the formula 
$$\hat q(x)=\bigl(1-\phi(x)\bigr)  \hat\pi^u_i(x) + \phi(x)\bigl(\pi_{k,\ell}(x)\bigr).$$
Here, we use an affine manifold structure on each component of $\hat W^s_{i,i}$ 
by identifying it with the circle $\mathbb S^1:= \mathbb R/\mathbb Z$. So, any positively weighted barycentric combination  makes sense for a pair of points sufficiently close. When $x\in \hat W^s_{i,i}$, we have $\hat q(x)=x$. Then, by shrinking  the tube $\hat N^{s}_{i}$ once more if necessary we make  $\hat q$ be a fibration whose fibres are transverse to the lamination $\hat L^s_i$ and we make each leaf of $\hat {\mathtt{F}}^u_{j,i}, \,j>i,$ in every  box $B^\ell_k$ be contained in a fibre of $q$. Henceforth, pullbacking  that  strip (and its fibration) by $p_i$  and adding the unstable manifold $W^u_i$ provide the wanted $Q^s_i$ and its foliation ${\mathtt{F}}^u_i$ satisfying the required properties. So, the induction is proved. 

We also have the existence of $f$-invariant neighborhoods $Q^{u}_{1},\dots,Q^{u}_{beh(f)-1}$ of  the saddle 
points from $\Sigma_{1},\dots,\Sigma_{beh(f)-1}$ respectively, equipped with one-dimensional $f$-invariant foliations ${{\mathtt{F}}}^{s}_{1},\dots,{{\mathtt{F}}}^{s}_{beh(f)-1}$ with smooth leaves such that the following properties hold for each $i\in\{1,\dots,{beh(f)-1}\}$:

{\rm (iii)} the stable manifolds $W^s_i$ are leaves of the foliation ${{\mathtt{F}}}^{s}_{i}$ and each leaf of the foliation ${{\mathtt{F}}}^{s}_{i}$ is transverse to $L^u$;
 
{\rm (iv)} for any $1\leq j<i$ and $x\in Q^u_{i}\cap Q^u_{j}$, we have the inclusion 
$({{\mathtt{F}}}^s_{j,x} \cap{Q}^u_{i})\subset
{{\mathtt{F}}}^s_{i,x}$.

The proof is done  by an increasing induction  from $i=1$; it is skipped due to similarity to the previous one.

{\bf Step 2.} We prove that for each $i=1,\dots,beh(f)-1$ there exists  an $f$-invariant neighborhood $\tilde{\mathtt N}_{i}$ of the set $\Sigma_i$ which is contained in $Q^s_{i}\cap Q^u_{i}$ and such that the restrictions of the foliations ${{\mathtt{F}}}^u_i$ and ${{\mathtt{F}}}^s_i$ to $\tilde{\mathtt N}_i$ are transverse.

For this aim, let us choose a fundamental domain   $K^s_i$ of the restriction of $f$ to $W^s_i\setminus \Sigma_i$ and take  a tubular neighborhood $N(K^s_i)$ of  $K^s_i$ whose segment fibres are contained in leaves of ${\mathtt{F}}^u_i$. Due to property (i), ${\mathtt{F}}^u_i$ is transverse to $W^s_i$. Since ${\mathtt{F}}^s_i$ is a $C^1$-foliation, if $N(K^s_i)$ is small enough, the foliations ${\mathtt{F}}^s_i$ and ${\mathtt{F}}^u_i$ have transverse intersection in $N(K^s_i)$. Set  $$\tilde{\mathtt N}_i:=W^u_i\,\bigcup_{k\in\mathbb Z}f^k\left(N(K^s_i)\right). $$ 
This is a neighborhood  of $\Sigma_i$; it satisfies condition (v) and the previous properties (i)--(iv) still hold.

{\bf Step 3.} Let us show the existence of linearizable neighborhoods 
$\mathtt{N}_i\subset\tilde{\mathtt N}_i,~i=1,\dots,beh(f)-1$, for which the required foliation are the restriction to $\mathtt N_i$ of the foliation ${{\mathtt{F}}}^u_i$.  

Let $\sigma\in\Sigma_i$ and $\tilde {\mathtt{N}}_\sigma$ be a connected component of  $\tilde{ \mathtt{N}}_i$ containing $\sigma$. There is a homeomorphism $\varphi^u_\sigma:W^u_\sigma\to W^u_{O}$ (resp. $\varphi^s_\sigma:W^s_\sigma\to W^s_{O}$) conjugating the diffeomorphisms $f^{per~\sigma}|_{W^u_\sigma}$ and $a_{\nu_\sigma}|_{W^u_{O}}$ (resp. $f^{per~\sigma}|_{W^s_\sigma}$ and $a_{\nu_\sigma}|_{W^s_{O}}$). In addition, for any point $z\in\tilde{ \mathtt{N}}_\sigma$ there is unique pair of points $z_s\in W^s_\sigma,~z_u\in W^u_\sigma$ such that $z={\mathtt{F}}^s_{i,z_u}\cap {\mathtt{F}}^u_{i,z_s}$. We define a topological embedding $\tilde{\mu}_{\sigma}:\tilde{ \mathtt{N}}_{\sigma}\to\mathbb R^2$ by the formula $\tilde{\mu}_{\sigma}(z)=(x_1,x_2)$ where  
$x_1={\varphi}^u_{\sigma}(z_u)$ and $x_2={\varphi}^s_{\sigma}(z_s)$. Choose $t_0\in(0,1]$ such that $\mathcal N^{t_0}\subset\tilde\mu_\sigma(\tilde{\mathtt{N}}_{\sigma})$ and $\partial\mathcal N^{t_0}$ does not contain images with respect $\tilde\mu_\sigma$ of heteroclinic points. Observe that $a_{\nu_\sigma}\vert_{\mathcal{N}^{t_0}}$ is conjugate to   $a_{\nu_\sigma}\vert_{\mathcal{N}}$ by the suitable homothety $h$. Set $\mathtt{N}_{\sigma}=\tilde\mu^{-1}_\sigma(\mathcal{N}^{t_0})$ and ${\mu}_{\sigma}=h \tilde\mu_\sigma:\mathtt{N}_{\sigma}\to \mathcal{N}$. Then, $\mathtt{N}_\sigma$ is the wanted neighbourhood with its linearizing homeomorphism ${\mu}_{\sigma}$. 
\end{demo}

Let us prove that for every diffeomorphism $f\in MS(M^2)$ there is a maximal $u$-compatible system of  neighbourhoods.

\begin{demo} To prove the theorem it remains to show how to do the constructed in Lemma \ref{comp} compatible system of neighborhoods $\mathtt N_f$ by maximal.  

Firstly define $N_1=\mathtt{N}_1$ and $N_{beh(f)-1}=\mathtt{N}_{beh(f)-1}$ because there is no heteroclinic rectangle connected with saddle points from $\Sigma_1$ and $\Sigma_{beh(f)-1}$. Let us describe a modification  of $\mathtt{N}_i,\,i\in\{2,\dots,beh(f)-2\}$ up to $N_i$ which is a maximal.   

For each $i\in\{1,\dots,beh(f)-1\}$, $\sigma\in\Sigma_i$ let 
$$\hat{\mathtt{N}}_\sigma=p_{_f}(\mathtt{N}_\sigma),\,\hat{\mathtt{N}}_i= \bigcup\limits_{\si\in\Si_i}\hat{\mathtt{N}}_\si,\,\hat{\mathtt{N}}_f= \bigcup\limits_{\si\in\Si_i}\hat{\mathtt{N}}_i,.$$  
If the set $\hat{\mathcal V}_f\setminus\hat{\mathtt{N}}_f$ does not contain a connected component whose boundary is a curvilinear triangle with sides $l_1\subset\hat{\mathtt N}_{i_1},\,l\subset\hat{\mathtt N}_{i},\,l_2\subset\hat{\mathtt N}_{i_2}$ for some $1\leq i_1<i<i_2\leq beh(f)-1$ then $N_i=\mathtt{N}_i$ for every $i\in\{2,\dots,beh(f)-2\}$. In the opposite case, for every such component $\hat\Delta$ there are saddle points $\sigma_1\in\Si_{i_1},\,\sigma\in\Si_{i},\,\sigma_2\in\Si_{i_2}$ such that $W^s_{\sigma_{1}}\cap W^u_\si\neq\emptyset$, $W^s_{\sigma}\cap W^u_{\si_2}\neq\emptyset$ and $W^s_{\sigma_{1}},\, W^u_\si,\, W^s_{\sigma},\,W^u_{\si_2}$ form a heteroclinic rectangle $\Pi_\si$ for which the set $\Delta=\Pi_\si\setminus (N_\si\cup N_{\si_1}\cup N_{\si_2})$ is a connected component of the preimage $p_f^{-1}(\hat\Delta)$. Let $N(\Pi_\si)$ be a neighborhood of $\Pi_\si$ bonded by $\mu_{\si_1}(\partial N^{1+\ep}),\,\mu_{\si_2}(\partial N^{1+\ep}),\,W^s_\si,W^u_\si$ for some small enough $\ep>0$. Let ${\mathtt N}^*_{\si}=\mathtt{N}_\si\cup\bigcup\limits_{k\in\mathbb Z}f^k(N(\Pi_\si))$. Due to compatibility of neighborhood $\mathtt{N}_f$ we can construct on ${\mathtt N}^*_{\si}$ a new unstable foliation $\mathtt{F}^{u*}_\si$ which is compatible with the unstable foliations in all saddle neighborhoods. Let us show that this neighborhood is a linearizable. 

For this aim it is enough to construct a homeomorphism $g:\mathtt{N}^*_\si\to\mathtt{N}_\si$ sending the foliation $\mathtt{F}^{u*}_\si$ to $\mathtt{F}^{u}_\si$ and such that $gf^{per(\si)}=f^{per(\si)}g$. Notice that each connected component of the sets $p_i(\mathtt {N}^*_\si)$ and $p_i(\mathtt{N}_\si)$ is an annulus and $p_i(\mathtt{N}^*_\si)\supset p_i(\mathtt{N}_\si)$. Let us choose a neighborhood $N(p_i(W^s_\si))$ of $p_i(W^s_\si)$ such that $N(p_i(W^s_\si))\subset p_i(\mathtt{N}_\si)$. Let us define a homeomorphism $\hat g:p_i(\mathtt{N}^*_\si)\to p_i(\mathtt{N}_\si)$ such that $\hat g|_{N(p_i(W^s_\si))}=id$,  $\hat g(p_i(\mathtt{N}^*_\si)\setminus N(p_i(W^s_\si)))=p_i(\mathtt{N}_\si)\setminus N(\hat W^s_\si))$ and $\hat g$ preserves leaves of the foliation $p_i(\mathtt{F}^{u*}_\si)$. The required homeomorphism $g$ is a lift of $\hat g$ which is identity in a neighborhood of $W^s_\si\cup W^u_\si$.   

Thus we get a new  $u$-compatible system of neighbourhoods for which the set $\hat{\mathcal V}_f\setminus(\bigcup\limits_{p\in\Om_f\setminus\si}\hat{\mathtt{N}}_p\cup\hat{\mathtt{N}}^*_\si)$ contains less by 1 connected components. When we do the same operation for all such components $\hat\Delta$ we get the desired maximal $u$-compatible system of neighborhoods.
\end{demo}

\section{Proof of Theorem \ref{t.cla}:  classifying diffeomorphisms}

In this section we prove Theorem \ref{t.cla}, i.e., we prove that two diffeomorphisms $f,f'\in MS(M^2)$ are topologically conjugate if and only if their schemes are equivalent.

\begin{demo} $ $

{\it Necessity.} Let  $f,f'\in MS(M^2)$ be two diffeomorphisms which are topologically conjugated by means of a homeomorphism $h:M^2\to M^2$. The conjugating homeomorphism sends the invariant manifolds of the periodic points of $f$ to corresponding objects of $f'$, preserving stability and periodicity,  and therefore 
$beh(f)=beh(f')$, $k_f=k_{f'}$. Thus $h$ induces a maximal system of compatible neighborhoods $\{h(N_\sigma),\sigma\in\Sigma\}$ for $f'$ which is different from $\mathcal N_{f'}$ in the general case. Then it remains to prove the following lemma.

\begin{lemm}\label{mama} The class of the equivalence of a scheme $S_f,\,f\in MS(M^2)$ does not depend on a choice of a  maximal system of compatible neighborhoods.
\end{lemm}
\begin{demo} Let $\mathcal N_f=\{N_1,\dots,N_{beh(f)-1}\}$ and $\mathsf{N}_f=\{\mathsf N_1,\dots,\mathsf N_{beh(f)-1}\}$ be different maximal systems of compatible neighborhoods with the foliations $F^u,\,\mathsf{F}^u$, accordingly. Without loss of generality we can assume that $\mathsf N_i\subset int\,N_i$ for each $i\in\{1,\dots,beh(f)-1\}$ (in the opposite case we can choose a maximal system $\mathbb{N}_{f}=\{\mathbb{N}_1,\dots,\mathbb{N}_{beh(f)-1}\}$ such that $\mathbb N_i\subset int\,N_i$ and $\mathbb N_i\subset int\,\mathsf N_i)$.

Let $G_i=\bigcup\limits_{j=1}^i N_i$, ${\mathsf G}_i=\bigcup\limits_{j=1}^i{\mathsf N}_j$, $\hat G_i=p_{_f}(G_i)$, $\hat{\mathsf G}_i=p_{_f}(\mathsf G_i)$ and $\hat G_{i,j}=p_j(G_i)$, $\hat{\mathsf G}_{i,j}=p_j(\mathsf G_i)$ for $j\in\{1,\dots,beh(f)-1\}$. By the increasing induction on $i$ from $1$ to $beh(f)-1$ let us show that $\hat G_i\setminus\hat {\mathsf G}_{i}$ is a direct product. 

For $i=1$, $\hat G_1$ and $\hat{\mathsf G}_{1}$ are tubular neighborhoods of a finite number circles $p_f(W^u_{1})$. As $\hat{\mathsf G}_{1}\subset int\,\hat G_1$ then, by the Annulus conjecture, $\hat G_1\setminus\hat {\mathsf G}_{1}$ is a direct product. 

Supposing that the statement is true for $i$ let us prove it for $i+1$. 

By assumption $\hat G_i\setminus\hat {\mathsf G}_{i}$ is a direct product. As $\hat G_{i,i}\setminus int\,\hat{\mathsf G}_{i,i}=p_i(p_f^{-1}(\hat G_i\setminus int\,\hat{\mathsf G}_i))$ then $\hat G_{i,i}\setminus int\,\hat{\mathsf G}_{i,i}$ is a direct product. Let us show that the intersections of the circles $\hat W^u_{i+1,i}$ with the sets $\hat G_{i,i}$ and $\hat{\mathsf G}_{i,i}$ consists of a finite number segments  $I_{1,i},\dots,I_{n_i,i}$ and $\mathsf{I}_{1,i},\dots,\mathsf{I}_{n_i,i}$, accordingly, such that $\mathsf I_{j,i}\subset int\, {I}_{j,i}$ for each $j\in\{1,\dots,n_i\}$. It will imply that $\hat G_{i+1}\setminus\hat {\mathsf G}_{i+1}$ is a direct product. 

Let $F^u_i$ be the corresponding foliation on $\mathcal N_i$ associated to  the foliations $F^u_\si$, $\si\in\Si_i$. Indeed, the circles $\hat W^s_{i,i}$ and $\hat W^u_{i+1,i}$ have a transversal intersection along a finite number $n^i_i$ points $z^i_{1,i},\dots,z^i_{n^i_i,i}$; $\hat N_{i,i},\,\hat{\mathsf N}_{i,i}$ are tubular neighborhoods of the circles $\hat W^s_{i,i}$ and every connected component of the intersections $\hat W^u_{i+1,i}\cap \hat N_{i,i}$, $\hat W^u_{i+1,i}\cap \hat{\mathsf N}_{i,i}$  is a leaf of the foliations $p_i(F^u_i),\,p_i(\mathsf{F}^u_i)$, accordingly.  Hence, the each intersection $\hat W^u_{i+1,i}\cap \hat N_{i,i}$ and $\hat W^u_{i+1,i}\cap \hat{\mathsf N}_{i,i}$ consist of $n^i_i$ segments $I^i_{1,i},\dots,I^i_{n^i_i,i}$ and $\mathsf{I}^i_{1,i},\dots,\mathsf{I}^i_{n^i_i,i}$, passing through the points $z^i_{1,i},\dots,z^i_{n^i_i,i}$, accordingly, such that $\mathsf I^i_{j,i}\subset int\, {I}^i_{j,i}$ for each $j\in\{1,\dots,n^i_i\}$. 

As the sets $\hat W^s_{i-1,i}\setminus \hat N_{i,i})$ and $\hat W^s_{i-1,i}\setminus \hat{\mathsf N}_{i,i}$ are compact then the intersections $(\hat W^s_{i-1,i}\setminus \hat N_{i,i})\cap\hat W^u_{i+1,i}$ and $(\hat W^s_{i-1,i}\setminus \hat{\mathsf N}_{i,i})\cap\hat W^u_{i+1,i}$ consist of a finite number points, which are the same due to properties of the maximal neighborhood, denote by $z^{i-1}_{1,i},\dots,z^{i-1}_{n^{i-1}_i,i}$ these points. Similar to the arguments above, the each intersection $\hat W^u_{i+1,i}\cap (\hat N^s_{i-1,i}\setminus \hat N_{i,i})$ and $\hat W^u_{i+1,i}\cap (\hat{\mathsf N}^s_{i-1,i}\setminus \hat{\mathsf N}_{i,i})$ consist of $n^{i-1}_i$ segments $I^{i-1}_{1,i},\dots,I^{i-1}_{n^{i-1}_i,i}$ and $\mathsf{I}^{i-1}_{1,i},\dots,\mathsf{I}^{i-1}_{n^{i-1}_i,i}$, passing through the points $z^{i-1}_{1,i},\dots,z^{i-1}_{n^{i-1}_i,i}$, accordingly, such that $\mathsf I^{i-1}_{j,i}\subset int\, {I}^{i-1}_{j,i}$ for each $j\in\{1,\dots,n^{i-1}_i\}$.  

Similarly for every $k\in\{1,\dots,i-2\}$ we have that each intersection $\hat W^u_{i+1,i}\cap (\hat N^s_{k,i}\setminus\bigcup\limits_{j=k+1}^{i} \hat N_{j,i})$ and $\hat W^u_{i+1,i}\cap (\hat {\mathsf N}^s_{k,i}\setminus\bigcup\limits_{j=k+1}^{i} \hat{\mathsf N}_{j,i})$ consist of $n^{k}_i$ segments $I^{k}_{1,i},\dots,I^{k}_{n^{k}_i,i}$ and $\mathsf{I}^{k}_{1,i},\dots,\mathsf{I}^{k}_{n^{k}_i,i}$, passing through the points $z^{k}_{1,i},\dots,z^{k}_{n^{k}_i,i}$, accordingly, such that $\mathsf I^{k}_{j,i}\subset int\, {I}^{k}_{j,i}$ for each $j\in\{1,\dots,n^{k}_i\}$. Thus we get the required statement.  

To finish the proof let us choose a tubular neighborhood $\hat U_-$ of $\partial \hat{\mathsf G}_{beh(f)-1}$ and tubular neighborhood $\hat U_+$ of $\partial \hat{G}_{beh(f)-1}$ avoiding all $p_{_f}(W^u_{\Si})$. As $\hat G_{beh(f)-1}\setminus\hat {\mathsf G}_{beh(f)-1}$ is a direct product  there is a homeomorphism $\hat \vp:\hat V_f\to\hat V_f$ which is identity on $\hat{\mathsf G}_{beh(f)-1}\setminus U_-$ and out of $\hat{G}_{beh(f)-1}\cup U_+$ and such that $\hat\vp({\mathsf G}_{beh(f)-1})={G}_{beh(f)-1}$. It is easy to modify $\vp$ such that $\hat\vp(\hat{\mathsf U}_{\sigma})=\hat U_\si$.  
\end{demo}

{\it Sufficiency.} Let us prove the sufficiency of the condition in Theorem \ref{t.cla}. Assume that the two diffeomorphisms $f,f'\in MS(M^2)$ have equivalent schemes.The schemes $S_f\mbox{ and }S_{f'}$ of diffeomorphisms  $f, f' \in MS(M^2)$, respectively, are said to be equivalent if  there exist an orientation-preserving homeomorphism $\hat\varphi:\hat{\mathcal V}_f\to\hat{\mathcal V}_{f'}$ such that:

1) $\eta_{f'}\hat\varphi_*=\eta_{f}$;

2) $\hat\varphi(\hat{\mathcal U}_f)=\hat{\mathcal U}_{f'}$, moreover for every point $\sigma\in\Sigma_i$ there is a point $\sigma'\in\Sigma'_i$ such that $\hat\varphi(\hat U_\sigma)=\hat U_{\sigma'}$ and  {$\varphi(U_\sigma)=U_{\sigma'}$}.

In a sequence of lemmas we will construct a maximal system of compatible neighborhoods $\{\mathsf N_1,\dots,\mathsf N_{beh(f)-1}\}$ for $f$ such that $\mathsf N_i\subset int(\mathcal N_i)$ for every $i\in\{1,\dots, beh(f)-1\}$ and a conjugating $f$ with $f'$ embedding $\psi$ on the union of these neighborhoods such that $\psi(\mathsf D)\subset int(\mathcal D')$, where $$\mathsf D=\bigcup\limits_{i=1}^{beh(f)-1}\mathsf N_i,\,\,\,\mathcal D=\bigcup\limits_{i=1}^{beh(f)-1}\mathcal N_i,\,\,\,\mathcal D'=\bigcup\limits_{i=1}^{beh(f)-1}\mathcal N'_i.$$ Using Lemma \ref{mama} we can interpolate $\psi$ on $\bigcup\limits_{i=1}^{beh(f)-1}\mathsf N_i$ with $\vp$ on $\bigcup\limits_{i=1}^{beh(f)-1}N_i$ and get a homeomorphism $h:M^2\setminus(\Si_0\cup\Si_{beh(f)})\to M'^2\setminus(\Si'_0\cup\Si'_{beh(f')}$ conjugating $f|_{M^2\setminus(\Si_0\cup\Si_{beh(f)})}$ with $f'|_{M\setminus(\Si'_0\cup\Si'_{beh(f')})}$. Notice that $M^2\setminus(W^s_{\Si}\cup \Si_{beh(f)})=W^s_{\Si_0}$ and $M^2\setminus(W^s_{\Si}\cup\Si'_{beh(f')})=W^s_{\Si'_0}$. Since $h(W^s_{\Si})=W^s_{\Si'}$ then $h(W^s_{\Si_0}\setminus\Si_0)=W^s_{\Si'_0}\setminus\Si'_0$. Thus for each connected component $A$ of $W^s_{\Si_0}\setminus\Si_0$, there is a sink $\Omega\in\Si_0$ such that $A=W^s_\omega\setminus\omega$. Similarly, $h(A)$ is a connected component of $W^s_{\Si'_0}\setminus\Si'_0$ such that  $h(A)=W^s_{\omega'}\setminus\omega'$ for a sink $\omega'\in\Si'_0$. Then we can continuously extend $h$ to $\Si_0$ {by defining} $h(\omega)=\omega'$ for every $\omega\in\Si_0$. A similar extension of $h$ to $\Si_{beh(f)}$ finishes the proof {of Theorem \ref{t.cla}}.  

To construct the embedding $\psi$, firstly, similar to proof of Lemma \ref{comp}, we can prove the existence in the neighborhoods $\mathcal N_1,\dots,\mathcal N_{beh(f)-1}$ of $f$-invariant one-dimensional foliations ${{F}}^{s}_{1},\dots,{{F}}^{s}_{beh(f)-1}$ whose leaves are smooth and such that the following properties hold for each $i\in\{1,\dots,beh(f)-1\}$:

- the stable manifolds $W^s_i$ are leaves of the foliation ${{F}}^{s}_{i}$ and 
the foliation ${{F}}^{s}_{i}$ is transverse to the foliation $F^u_i$;
 
- for any  $1\leq i <k\leq beh(f)-1$ and $x\in \mathcal N_{i}\cap \mathcal N_{k}$, we have the inclusion 
${F}^s_{i,x}\cap \mathcal N_{k}\subset {F}^s_{k,x}$. 

By property 2) of the equivalence of the schemes, we have that there is a one-to-one correspondence between the sets $\Sigma_i$ and $\Sigma'_i$ through the equality $\vp(U_\si)=U_{\si'}$. For $j>i$, let us denote by ${\mathcal{J}}^u_{i,j}$ the union of all connected components of ${W}^u_{j}\cap N_i$ which do not lie in {$int\,N_k$} with $i<k<j$. Let $\mathcal{J}^u_i=\bigcup\limits_{j=i+1}^{beh(f)-1}{\mathcal{J}}^u_{i,j}$. Also for $j<i$, let us denote by ${\mathcal{J}}^s_{i,j}$ the union of all connected components of ${W}^s_{j}\cap N_i$ which do not lie in {$int\,N_k$} with $j<k<i$. Let $\mathcal{J}^s_i=\bigcup\limits_{j=1}^{i-1}{\mathcal{J}}_{i,j}$. Also we introduce the same sets $\mathcal J'^u_i$ and $\mathcal J'^s_i$ for $f'$. The embedding $\vp$ gives a one-to-one correspondence between the connected components $J^u$ of $\mathcal{J}^u_{i,j}$ and the connected components $J'^u$ of $\mathcal{J}'^u_{i,j}$ by the rule: there are exactly two connected components $B_{J^u}$ of $\partial U_i\cap U_j$ such that $\partial J^u\subset B_{J^u}$ and the same for prime, then $B_{J'^u}=\vp(B_{J^u})$.  

\begin{lemm} \label{1dim} There is a homeomorphism $\psi^u:W^u_\Sigma\to W^{\prime u}_{\Sigma'}$ consisting of conjugating  homeomorphisms $\psi^u_{1}:W^u_{1}\to W^{\prime u}_{1},\dots,\psi^u_{beh(f)-1}:W^u_{beh(f)-1}\to W^{\prime u}_{beh(f)-1}$ such that:

1) $\psi^u_{i}(J^u)\subset J'^u$ for every segment $J^u\subset \mathcal J^u_{j,i},j<i$ and $$\psi^u_i(F^s_{i,y}\cap W^u_i)=F'^s_{i,\psi^u_i(y)}\cap W'^u_i,~~for~every~~y\in J^u;$$

2) $\psi^u_i(W^u_i\cap L^s\cap N_j)=W'^u_i\cap L'^s\cap N'_j,j<i$.
\end{lemm}
\begin{demo} This statement is proved step by step from $i=1$ to $i=beh(f)-1$. 
 
First, take $i=1$. Let $\si\in\Si_1$ and $\ga^u_\si$ is the unstable separatrix of $\si$. Let us choose a point  $x\in\ga^u_\si$. Let $a\in (F^s_{\si,x}\cap\partial U_\si)$, $a'=\vp(a)$. Then the point $F'^s_{\si',a'}\cap W^u_{\si'}$ belongs to the unstable separatrix of $\si'$ which we denote by $\ga_{\si'}$. Let us choose a point $x'\in\ga_{\si'}$ such that: 

(*) if $F^s_{\si,x}\cap J^u\neq\emptyset$ for $J^u\subset\mathcal J^u_{1}$ then $F'^s_{\si',x'}\cap J'^u\neq\emptyset$ also.

Denote by $I$ the segment $[x,f^{per(\ga_\si)}(x)]\subset\ga_\si$, by $I'$ the segment $[x',f^{per(\ga_{\si'})}(x')]\subset\ga_{\si'}$ and by $\psi^u_I:I\to I'$ a homeomorphism such that $\psi^u_I(x)=x'$. After that we extend it up to a homeomorphism $\psi^u_{\ga_\si}:\ga_\si\to\ga_{\si'}$ by the formula $$\psi^u_{\ga_\si}(y)=f'^{-n\cdot per(\ga_{\si'})}(\psi^u_I(f^{n\cdot per(\ga_\si)}(y))),$$ where $f^{n\cdot per(\ga_\si)}(y)\in I$. Let us do the similar construction for each unstable separatrices of $\Si_1$, satisfying the condition $\psi^u_{f(\ga_\si)}=f'\psi^u_{\ga_\si}f^{-1}$. Then we get a homomorphism on $W^u_1\setminus\Si_1$ which can be continuously extended up to the homeomorphism $\psi^u_1:W^u_1\to W'^u_1$. Due to $(*)$ we can define an embedding $\psi^u_{J^u}$ on every $J\subset\mathcal J^u_{1,j},j>1$ by the rule $y'=\psi^u_{J^u}(y),y\in J^u$, where  $$\psi^u_1(F^s_{1,y}\cap W^u_1)=F'^s_{1,y'}\cap W'^u_1.$$ 

For $i=2$ we do the similar construction for $\psi^u_2$ assuming additionally that $\psi^u_2|_{J^u}=\psi^u_{J^u}|_{J^u}$ for $J^u\in\mathcal J^u_{1,2}$.  Also we can define an embedding $\psi^u_{J^u}$ on every $J^u\subset\mathcal J^u_{2,j},j>2$ by the rule $y'=\psi^u_{J^u}(y),y\in J$, where  $$\psi^u_2(F^s_{2,y}\cap W^u_2)=F'^s_{2,y'}\cap W'^u_2.$$ Continue this process we get the desired embeddings $\psi^u_3,\dots,\psi^u_{beh(f)-1}$.
\end{demo}

By the construction the homeomorphism $\psi^u$ send the heteroclinic points of $f$ to the heteroclinic points of $f'$. Then for every connected component $J^s$ of $\mathcal J^s_{i,j}$ passing through the point $x\in W^u_i$ we will denote by $J'^s$ a connected component  of $\mathcal J'^s_{i,j}$ passing through the point $\psi^u_i(x)\in W'^u_i$. Absolutely similar to proof of Lemma \ref{1dim}, but starting from $i=beh(f)-1$ to $i=1$, we can prove the following result. 

\begin{lemm} \label{sdim} There is a homeomorphism $\psi^s:W^s_\Sigma\to W^{\prime s}_\Sigma$ consisting of conjugating  homeomorphisms $\psi^s_{1}:W^s_{1}\to W^{\prime s}_{1},\dots,\psi^s_{beh(f)-1}:W^s_{beh(f)-1}\to W^{\prime s}_{beh(f)-1}$ such that:

1) $\psi^s_{i}(J^s)\subset J'^s$ for every segment $J^s\subset \mathcal J^s_{j,i},j>i$ and $$\psi^s_i(F^u_{i,y}\cap W^s_i)=F'^u_{i,\psi^s_i(y)}\cap W'^s_i,~~for~every~~y\in J^u;$$

2) $\psi^s|_{L^s\cap L^u}=\psi^u|_{L^s\cap L^u}$.
\end{lemm}

Finitely for each $i\in\{1,\dots,beh(f)-1\}$ we construct an embedding $\psi_i:\mathcal N_i\to\mathcal N'_i$ by the formula: for a point $x=F^s_{i,x^u}\cap F^u_{i,x^s},x^u\in W^u_i,x^s\in W^s_i$ we have $\psi_i(x)=F'^s_{i,\psi^u_i(x^u)}\cap F'^u_{i,\psi^s_i(x^s)}$. It follows from Lemmas \ref{1dim} and \ref{sdim} that $\psi_i(\mathcal N_i)\subset int(\mathcal N'_i)$ and a map $\psi$ composed by $\psi_1,\dots,\psi_{beh(f)-1}$ is a homeomorphism. A choice of a maximal system of $u$-compatible neighborhoods with property $\mathsf{N}_i\subset\mathcal N_i$ finishes the proof.  
\end{demo} 

\section{Proof of  Theorem \ref{t.realisation}: realizing diffeomorphisms}\label{rea}

In this section we prove that for any abstract scheme $S\in\mathcal S$ 
there is a diffeomorphism $f\in MS(M^2)$ whose scheme is equivalent to the scheme $S$.

\begin{demo} Let  $S=(\hat{\mathcal V},\eta,\hat{\mathcal U})$ be an abstract 
scheme. Let us construct step by step a diffeomorphism $f\in MS(M^2)$ such that the 
schemes $S_{f}$ and $S$ are equivalent.\\ 
{\bf Step 1.} It follows from the definition of an abstract scheme that $\hat{\mathcal V}$ is a disjoint union of the finite number 2-tori $\hat{V}_1,\dots,\hat{V}_k$ with the map $\eta_{_{\hat{\mathcal V}}}$ composed from the non-trivial homomorphisms $\eta_{_{\hat{V}_1}},\dots,\eta_{_{\hat{V}_k}}$. Then for each $i\in\{1,\dots,k\}$ there is a number $m_{_{\hat V_i}}\in\mathbb N$ such that $\hat{V}_i=((\mathbb R^2\setminus O)\times\mathbb Z_{m_{_{\hat V_i}}})/b_i$, where $b_{i}:\mathbb R^2\times\mathbb Z_{m_{_{\hat V_i}}}\to\mathbb R^2\times\mathbb Z_{m_{_{\hat V_i}}}$ given by the formula 
$$b_{i}(x_1,x_2,\la)=\begin{cases}(\frac{x_1}{2},\frac{x_2}{2},\la+1),
\la\in\{1,\dots,m_{_{\hat V_i}}-1\};\\ (\frac{x_1}{2},\frac{x_2}{2},1),~\la=m_{_{\hat V_i}}\end{cases}$$
and the natural projection $p_i:(\mathbb R^2\setminus O)\times\mathbb Z_{m_{_{\hat V_i}}}\to\hat{V}_i$ induces the homomorphism $\eta_{_{\hat{V}_i}}$. Denote by $\mathcal W$ the disjoint union of $\mathbb R^2\times\mathbb Z_{m_{_{\hat V_1}}},\dots,\mathbb R^2\times\mathbb Z_{m_{_{\hat V_k}}}$ and by $f_{\mathcal W}$ a diffeomorphism composed by $a_1,\dots,a_r$. Also denote by $\mathcal V$ the disjoint union of $(\mathbb R^2\setminus O)\times\mathbb Z_{m_{_{\hat V_1}}},\dots,(\mathbb R^2\setminus O)\times\mathbb Z_{m_{_{\hat V_k}}}$ and by $p_{_{\mathcal V}}:\mathcal V\to\hat{\mathcal V}$ the natural projection.  

{\bf Step 2.} Let $\hat U_1=\bigcup\limits_{j=1}^{n_1}\hat U_{1,j}$. For each $j\in\{1,\dots,n_1\}$ let $U_{1,j}=p^{-1}_{\mathcal V}(\hat{U}_{1,j})$. It follows from the definition of the abstract scheme that  there are numbers $m_{1,j}\in\mathbb N,~\nu_{1,j}\in\{-1,+1\}$, a canonical neighborhood $N$ and a diffeomorphism $\mu_{1,j}:U_{1,j}\to N\times\mathbb Z_{m_{1,j}}$ which conjugate $f_{\mathcal W}|_{U_{1,j}}$ with a diffeomorphism  $a_{1,j}|_{N\times\mathbb Z_{m_{1,j}}}$ given on $N\times\mathbb Z_{m_{1,j}}$ by the formula 
$$a_{1,j}(x_1,x_2,\la)=\begin{cases}(2x_1,\frac{x_2}{2},\la+1),
\la\in\{1,\dots,m_{1,j}-1\};\\ (\nu_{1,j}\cdot 2x_1,\nu_{1,j}\cdot \frac{x_2}{2},1),~\la=m_{1,j}.\end{cases}$$

Set $A_{{1,1}}=\mathcal W\setminus int\,U_{1,1},~B_{{1,1}}={N}\times\mathbb Z_{m_{1,1}}$ and ${Q_{{1,1}}}=A_{1,1}\cup_{\mu^{-1}_{1,1}}B_{1,1}$, $\bar{Q}_{1,1}=A_{1,1}\cup B_{1,1}$. Denote by $p_{_{Q_{1,1}}}:\bar Q_{1,1}\to Q_{1,1}$ the natural projection then projections  
 $p_{_{A_{1,1}}}=p_{_{Q_{1,1}}}|_{A_{1,1}}$,  $p_{_{B_{1,1}}}=p_{_{Q_{1,1}}}|_{B_{1,1}}$ induce structure of smooth connected orientable separable $3$-manifold without boundary (possible non Hausdorff). Let us show that the space  ${Q_{1,1}}$ is Hausdorff. 

For this aim it is enough to prove that a set  
$E_{Q_{{1,1}}}=\{(x,y)\in \bar Q_{{1,1}}\times\bar Q_{{1,1}}~:~p_{_{{Q_{{1,1}}}}}(x)=p_{_{{Q_{{1,1}}}}}(y)\}$ is closed in  $\bar Q_{{1,1}}\times \bar Q_{{1,1}}$ (see, for example, a book \cite{Ko} by Kosnevski). It is equivalent that $(x,y)\in E_{Q_{{1,1}}}$ for any sequence $(x_m,y_m)\in E_{Q_{{1,1}}}$ converging in the space $\bar Q_{{1,1}}\times\bar Q_{{1,1}}$ to a point $(x,y)$. Without loss of generality we can suppose that all points of the sequence $x_m~(y_m)$ belong to the same connected component of  $\bar Q_{{1,1}}$ as  $x~(y)$ (in the opposite case it is possible to consider a subsequence  with such property). Let us consider four possibilities:
1) $x_m,y_m\in A_{{1,1}}$; 2) $x_m,y_m\in B_{{1,1}}$; 3) $x_m\in A_{{1,1}},y_m\in B_{{1,1}}$; 4) $x_m\in B_{{1,1}},y_m\in A_{{1,1}}$. 

In cases 1) and 2), $x_m=y_m$. Then  $x=y$ and, hence, $(x,y)\in E_{Q_{{1,1}}}$. In case 
3) $x_m\in A_{{1,1}}$, $y_m\in B_{{1,1}}$, hence $y_m\in\partial N_{1,1}$ and $x_m=\mu_{{1,1}}(y_m)$. As $\partial N_{1,1}$ is closed in $\bar Q_{{1,1}}$ then, using continuation of the map $\mu_{{1,1}}$, we get next series of equalities: $x=\lim\limits_{m\to\infty}x_m=\lim\limits_{m\to\infty}\mu_{{1,1}}(y_m)=
\mu_{{1,1}}(\lim\limits_{m\to\infty}(y_m))=\mu_{{1,1}}(y)$. Thus, 
$(x,y)\in E_{Q_{{1,1}}}$. In case 4) similar to above it is possible to prove that $(x,y)\in E_{Q_{{1,1}}}$.
 
Thus $Q_{{1,1}}$ is smooth connected orientable $2$-manifold without boundary. 
Set  $f_{A_{{1,1}}}=p_{_{A_{{1,1}}}}f_{\mathcal W} p_{_{A_{{1,1}}}}^{-1}:p_{_{A_{{1,1}}}}(A_{{1,1}})\to p_{_{A_{{1,1}}}}(A_{{1,1}})$ and $f_{B_{{1,1}}}=p_{_{B_{{1,1}}}}a_{{1,1}}p_{_{B_{{1,1}}}}^{-1}:p_{_{B_{{1,1}}}}(B_{{1,1}})\to p_{_{B_{{1,1}}}}(B_{{1,1}})$. By the construction the maps $f_{A_{{1,1}}}$ and
$f_{B_{{1,1}}}$ are diffeomorphisms coinciding on set 
$p_{_{A_{{1,1}}}}(A_{{1,1}})\cap p_{_{B_{{1,1}}}}(B_{{1,1}})$. 
Then a map $f_{{Q_{{1,1}}}}:{Q_{{1,1}}}\to {Q_{{1,1}}}$ given by formula\\
$$f_{{Q_{{1,1}}}}(x)=\begin{cases}f_{_{A_{{1,1}}}}(x),
x\in p_{_{A_{{1,1}}}}(A_{{1,1}});
\\ f_{_{B_{{1,1}}}}(x),x\in
p_{_{B_{{1,1}}}}(B_{{1,1}})\end{cases}$$  
is a diffeomorphism of the manifold ${Q_{{1,1}}}$. 
By the construction non-wandering set of $f_{{Q_{{1,1}}}}$ 
consists of unique saddle periodic orbit and $r$ sink periodic orbits. 

We will do the same operation with the all connected components of $U_{1,j}\subset U_1$ to get a smooth connected orientable $2$-manifold without boundary $Q_{U_{1}}$ and a diffeomorphism 
$f_{{Q_{U_{1}}}}:Q_{U_{1}}\to Q_{U_{1}}$ with a finite set $\Sigma_1$ of the saddle periodic points and a finite set $\Sigma_0$ of the sink periodic points. 

Set $\check V=(Q_{{U}_{1}}\setminus W^u_{\Sigma_1\cup\Sigma_0})/_{f_{Q_{{U}_{1}}}}$ 
and denote by $p_{_{\check V}}:Q_{{U}_{1}}\setminus W^u_{\Sigma_1\cup\Sigma_0}\to\check V$ 
the natural projection. Then the manifold $\check V$ is equivalent to the manifold 
$\hat V_{\hat U_1}=\hat{\mathcal V}_2$. 

Continuing this process we get a smooth connected orientable noncompact 2-manifold $Q$ without
boundary and a diffeomorphism $f_{Q}:Q\to Q$ whose non-wandering set consists of finite set $\Sigma$ of the saddle periodic hyperbolic orbits.

{\bf Step 3.} Set $C={Q}\setminus W^u_{\Omega_{f_{Q}}}$.
Denote by $\hat{C}$ the space of orbit of action $f_{{Q}}$ on $C$ 
and by  $p_{_{\hat C}}:\hat C\to\hat C$ the natural projection. 
By the construction $C$ is obtained by surgery of $\hat{\mathcal V}_n$ along $\hat{\mathcal N}_n$ and, 
hence, due to Proposition \ref{rM}, it is homeomorhic to finite number (denote it $k_u$) of 
copies of the torus. Similar to Step 1 for each connected component $\hat c$ of $\hat C$ there is a number $m_{\hat c}\in\mathbb N$ and a diffeomorphism $b_{\hat c}:\mathbb R^2\times\mathbb Z_{m_{\hat c}}\to\mathbb R^2\times\mathbb Z_{m_{\hat c}}$ given by the formula 
$$b_{\hat c}(x_1,x_2,\la)=\begin{cases}(2x_1,2x_2,\la+1),
\la\in\{1,\dots,m_{c}-1\};
\\ (2x_1,2x_2,1),~\la=m_{\hat c}\end{cases}$$
such that $f_Q|_{p_{_{\hat C}}(\hat c)}$ topologically conjugated with $b_{\hat c}|_{(\mathbb R^2\setminus O)\times\mathbb Z_{m_{\hat c}}}$ by means a diffeomorphism $\mu_{\hat c}$. 

Denote by $D_{u}$ a set composed by $\mathbb R^2\times\mathbb Z_{m_{\hat c}},~\hat c\subset \hat C$, by  $D'_{u}$ a set composed by $(\mathbb R^2\setminus O)\times\mathbb Z_{m_{\hat c}},~\hat c\subset\hat C$ , by $\mu_C:C\to D'_u$ the map composed by the homeomorphisms $\mu_{\hat c},~\hat c\subset \hat C$ and by $b_C:D_u\to D_u$ a map composed by $b_{\hat c},~\hat c\subset\hat C$.
Set $M^2=Q\cup_{\mu_C}D_u$, $\bar M^2=Q\cup D_u$ and denote by  
$p_{_{M^2}}:\bar M^2\to M^2$ the natural projection. Like above a proof 
that the topological space $M^2$ is smooth connected orientable 2-manifold 
without boundary reduces to checking that set $E_{{M^2}}=\{(x,y)\in 
\bar {M^2}\times \bar {M^2}~:~p_{_{{M^2}}}(x)=p_{_{{M^2}}}(y)\}$ is closed in 
$\bar {M^2}\times \bar {M^2}$, that is if a sequence  $(x_m,y_m)\in E_{{M^2}}$ 
converges in $\bar {M^2}\times\bar {M^2}$ to a point $(x,y)$ then the point 
$(x,y)$ belongs to $E_{{M^2}}$.

Consider four cases:
1) $x_m,y_m\in Q$; 2) $x_m,y_m\in D_u$; 3) $x_m\in Q,y_m\in D_u$; 4) 
$x_m\in D_u,y_m\in Q$. 

In cases 1) and 2), $x_m=y_m$. Then  $x=y$ and, hence, $(x,y)\in E_{{M^2}}$. 
In case 3) $x_m\in Q$, $y_m\in D^{\prime}_{u}$,  $y_m=\mu_{_C}(x_m)$ and 
there are two subcases: 3a) $y\in D^{\prime}_{u}$; 3b) 
$y=O$. In subcase 3a), like to above, $x=\mu^{-1}_{C}(y)$ and, hence, 
$(x,y)\in E_{{M^2}}$. Show that case 3b) is impossible.

As $y_m\in D^{\prime}_{u}$ and $y=O$ then the sequence $x_m=\mu^{-1}_{C}(y_m)$ converges to $x\in W^u_{\Sigma_1}\cup\Sigma_0$. Then there is a sequence $k_m\to+\infty$ such that $f_Q^{-k_m}(x_m)\to z\in C$. Thus $b_C^{-k_m}(y_m)\to\mu_C^{-1}(z)$. It is contradiction because $b_C^{-k_m}(y_m)\to O$. 

The case 4) can be proved similarly to case 3). 

Set $p_{_{Q}}=p_{_{{M^2}}}|_{Q}$,  $p_{_{D_u}}=p_{_{{M^2}}}|_{D_u}$ and 
$f_{D_u}=p_{_{D_u}}d_{s}p_{_{D_u}}^{-1}:p_{_{D_u}}(D_u)\to p_{_{D_u}}(D_u)$. 
Similar to above we can prove that a map $f:{M^2}\to {M^2}$ given by formula\\
$f(x)=\begin{cases}f_{Q}(x), x\in p_{_{Q}}(Q);\\ f_{D_u}(x), x\in p_{_{D_u}}(D_u)\end{cases}$ is a diffeomorphism of the manifold ${M^2}$ whose non-wandering set consists of $k$ saddle periodic hyperbolic orbits and of $k_s$ sink periodic hyperbolic orbits and of $k_u$ source periodic hyperbolic orbits.

{\bf Step 4.} In this step we show that the manifold $M^2$ is compact and, hence, the 
diffeomorphism $f$ belongs to the class $MS(M^2)$ and its scheme by the construction 
is equivalent to the abstract scheme $S$.

For proof of compactness of $M^2$ it is enough to show that any sequence $\{x_n\}\in M^2$ 
has converging subsequence. If infinitely many members of $\{x_n\}$ belong to $\Omega_{f}$ the 
fact is obvious. Consider opposite case. By the construction  
$M^2=\bigcup\limits_{p\in\Omega_f}W^s_p=\bigcup\limits_{p\in\Omega_f}W^u_p$. 
Up to consider a  subsequence there is a point $p_1\in\Omega_f$ such that 
$\{x_n\}\subset(W^s_{p_1}\setminus p_1)$. Denote by $K$ fundamental domain of the restriction 
of $f$ to $W^s_{p_1}\setminus p_1$. Then for each member $x_n$ of the sequence $\{x_n\}$ 
there is an integer $k_n$ such that  $y_n=f^{k_n}(x_n)\in K$. Without loss of generality we 
can suppose that sequence $\{y_n\}=\{f^{k_n}(x_n)\}$ converges to a point $y\in K$ 
(in opposite case we can consider subsequence with such property). For the sequence 
$\{k_n\}$ there are two possibilities:

1)  $\{k_n\}$ is bounded;

2)  $\{k_n\}$ is not bounded.

In case 1), up to consider a  subsequence, the sequence $\{k_n\}$ converges to an integer 
$k$. Then  $\lim\limits_{n\to\infty}x_{n}=\lim\limits_{n\to\infty}f^{-k_n}(y_n)=f^{-k}(y)$. 
Thus a subsequence of $\{x_{n}\}$ converges to $f^{-k}(y)\in W^s_{p_1}$.

In case 2), up to consider a  subsequence,  $\{k_n\}$ converges to  $ + \infty$ or $-\infty $. 
In case $k_n\to -\infty $ a subsequence of  $\{x_{n}=f^{-k_n}(y_n)\}$ converges to $p_1$. In 
case $k_n\to +\infty $, up to consider a  subsequence, there is a point $p_2\in\Omega_f$ 
such that $\{x_n\}\subset(W^u_{p_2}\setminus p_2)$ and, hence, a subsequence of  
$\{x_{n}=f^{-k_n}(y_n)\}$ converges to $p_2$. 
\end{demo}

\end{document}